\documentclass[3p,times]{elsarticle}

\usepackage[T1]{fontenc}
\usepackage{lipsum}
\usepackage{amsfonts}
\usepackage{graphicx}
\usepackage{epstopdf}
\usepackage{algorithmic}
\usepackage{algorithm}
\usepackage{graphicx}
\usepackage{epstopdf}
\usepackage{url}
\usepackage{amsmath}
\usepackage{mathtools}
\usepackage{subcaption}
\usepackage{booktabs}
\usepackage{multirow}
\usepackage{color}

\usepackage{amsopn}

\newtheorem{remark}{Remark}

\newproof{proof}{Proof}

\definecolor{ao}{rgb}{0.0, 0.5, 0.0}

\makeatletter
\def\ps@pprintTitle{%
 \let\@oddhead\@empty
 \let\@evenhead\@empty
 \def\@oddfoot{}%
 \let\@evenfoot\@oddfoot}
\makeatother

\begin{document}

\begin{frontmatter}



\title{{ Anomalous Nonlinear Dynamics Behavior of Fractional Viscoelastic Structures}}

\author[1,2]{Jorge Suzuki}
\author[1,2]{Pegah Varghaei}
\author[4]{Ehsan Kharazmi}
\author[1,3]{Mohsen Zayernouri\corref{cor1}}
\cortext[cor1]{Corresponding author: zayern@msu.edu, Tel.(517) 432-0464}

\address[1]{Department of Mechanical Engineering, Michigan State University, MI 
48824, USA}
\address[2]{Department of Computational Mathematics, Science, and Engineering 
(CMSE), Michigan State University, MI 48824, USA}
\address[3]{Department of Statistics and Probability, Michigan State 
University, MI 48824, USA}
\address[4]{Division of Applied Mathematics, Brown University, RI, 02912, USA} 
\begin{abstract}
{Fractional models and their parameters are sensitive to changes in the intrinsic micro-structures of anomalous materials. We investigate how such physics-informed models propagate the evolving anomalous rheology to the nonlinear dynamics of mechanical systems. In particular, we analyze the} vibration of a fractional, {geometrically nonlinear} viscoelastic cantilever beam, under base excitation and free vibration, where the {viscoelastic response is general through a} distributed-order fractional model. We {employ} Hamilton's principle to {obtain the corresponding} equation of motion with the choice of specific material distribution functions that {recover a} fractional Kelvin-Voigt viscoelastic model {of order $\alpha$}. {Through} spectral decomposition in space, the resulting time-fractional partial differential equation reduces to a nonlinear time-fractional ordinary differential equation, {in which the linear counterpart is numerically integrated by employing a direct L1-difference scheme.} We further develop a semi-analytical scheme to solve the nonlinear system {through a method of multiple scales,} {which yields a} cubic algebraic equation {in terms of the frequency}. {Our numerical results suggest a set of $\alpha$-dependent anomalous dynamic qualities, such as far-from-equilibrium power-law amplitude decay rates, super-sensitivity of amplitude response at free vibration,} and bifurcation in steady-state amplitude at primary resonance.
\end{abstract}


\end{frontmatter}


\section{Introduction}

{Nonlinearities are inherent characteristics in many real physical systems arising from a variety of sources, such as significant changes in geometry, material properties (\textit{e.g.}, ageing), and boundary effects (\textit{e.g.}, development of boundary layers and shock). In this work we focus on the analysis of nonlinear systems subject to anomalous dynamics that arise from nonlocal/history effects. Despite the existence of ``nearly-pure" systems, in which standard features evolve to anomalous qualities, \textit{e.g.}, laminar-to-turbulent flows \cite{Sagaut2018, AkhavanSafaei2020} and dislocation pile-up in localized plastic yielding \cite{habtour2016detection}, in our study, the source of anomalies is due to the employment of extraordinary materials.} 

{Power-law rheology is a constitutive behavior observed in a wide range of anomalous materials. Such complex rheology exhibits macroscopic \textit{memory-effects} by means of single-to-multiple power-law relaxation/creep \cite{Kapnistos2008}, and dynamic storage/dissipation visco-elasto-plasticity \cite{McKinley2013}. These power-law characteristics are multi-scale fingerprints of spatial/temporal sub-diffusive processes \cite{Metzler2000} of heterogeneous, fractal-like micro-structures, where the mean squared displacement of constituents/defects follows a non-linear scaling in time of the form $\langle \Delta r \rangle^2 \propto t^\alpha$ \cite{Nnetu2013, Wong2004}. As anomalous materials undergo cyclic loads, they endure micro-structural changes, \textit{i.e.} rearrangement/unfolding of polymer networks/chains \cite{Kapnistos2008}, plastic stretching/buckling of micro-fibers \cite{Bonadkar2016}, formation, arresting, relaxation of dislocations \cite{Richeton2005}, among others. Such multi-scale physics change the characteristic fractal and spectral dimensions of the microstructure, affecting the small-scale diffusion mechanisms and therefore the micro/macro-rheological properties.}

{Classical (integer-order) viscoelastic models, \textit{i.e.}, Maxwell, Kelvin-Voigt \cite{christensen2012theory}, provide accurate fits for exponential-like relaxation data with a limited number of relaxation times
\cite{pipkin2012lectures, christensen2012theory}. However, complex arrangements and a large number of Hookean/Newtonian springs/dashpots is required to simply estimate the complex hereditary behavior observed for a broad class of anomalous/non-standard (non-exponential/power law) materials. This leads to high-dimensional parameter spaces, adversely affecting the conditioning of ill-posed inverse problems of parameter estimation \cite{aster2018parameter}. In addition, multi-exponential approximations merely represent a truncated power-law relaxation \cite{bagley1989power}, providing satisfactory representations only for short observation times, therefore lacking predictability and requiring recalibration for multiple time-scales \cite{Jaishankar2013}.}

{Fractional differential equations (FDEs) allow excellent predictability of material responses across multiple time-scales for anomalous materials. Nutting and Gemant \cite{nutting1921new,gemant1936method} demonstrated that power law kernels are more descriptive for creep and relaxation. Later on, Bagley and Torvik 
\cite{bagley1983theoretical} proposed a link of fractional viscoelasticity with molecular theories of polymers dynamics through frequency-dependent moduli. The basic building block of fractional viscoelasticity is the \textit{so-called} Scott-Blair (SB) element with fractional order $0 < \alpha < 1$, which provides a constitutive interpolation between Hookean springs $(\alpha \to 0)$ and Newtonian dashpots $(\alpha \to 1)$. Distinct mechanical arrangements of SB elements allow the modeling multiple experimentally observed power-laws through corresponding multi-term FDEs. Such flexible and compact mathematical tools allowed researchers to develop and employ fractional rheological models in diverse fields, such as bio-engineering\cite{naghibolhosseini2015estimation,naghibolhosseini2018fractional}, visco-elasto-plastic modeling for power-law strain hardening \cite{suzuki2016fractional}, among others \cite{shitikova2017interaction,rossikhin1997applications,Samiee2019turbulence}. The most general forms of viscoelastic constitutive laws are represented through distributed order differential equations (DODEs) \cite{lorenzo2002variable,Atanackovic2009distributional}, where the fractional-order distributions (and therefore distributions of SB elements) code the heterogeneous multi-scale material properties arising from the evolving material microstructure. In addition, specific choices of material distributions recover known discrete fractional models. DODEs are used in \cite{caputo1995m,caputo2001m} to generalize the stress-strain 
relation of inelastic media and Fick's law. Their connection with diffusion-like equations was established in \cite{chechkin2002,	chechkin2008generalized} and their applications are also discussed in time domain analysis of control, filtering and signal processing \cite{li2011distributed,li2014lyapunov}, vibration \cite{duan2018steady}, frequency domain analysis \cite{bagley2000}, and uncertainty quantification \cite{kharazmi2017petrov, kharazmi2018fractional}.}

{Regarding the dynamics of fractional visco-elastic beams, {\L}ab\c{e}dzki \textit{et al.} \cite{labkedzki2018transverse} investigated the resonant characteristics of an Euler-Bernoulli piezoelectric cantilever beam by replacing the usual sum of stiffness and damping terms in the strong form of the equation of motion by a single fractional derivative operator, and solved the system using a Rayleigh-Ritz method.  Ansari \textit{et al.} \cite{Ansari2016Nanobeam} analyzed the free vibration response of a fractional Kelvin-Voigt viscoelastic Euler-Bernoulli nanobeam with nonlocal elastic response. Their work employed a direct Ritz method for space discretization and a time-fractional Adams-Moulton scheme for the resulting time-fractional ODEs, and observed higher damping for higher fractional order values. Utilizing the same model, Faraji Oskouie \textit{et al.} \cite{FarajiOskouie2017} incorporated the effects of surface stresses using the Gurtin-Murdoch theory in simply supported and cantilever beams. More recently, Eyebe \textit{et al.} \cite{Eyebe2018} analyzed the nonlinear vibration of a nanobeam resting on a fractional order Winkler-Pasternak foundation, utilizing the D'Alembert principle to obtain the governing equations and a method of multiple scales to approximate the resulting nonlinear problem. In \cite{lewandowski2017nonlinear}, Lewandowski \textit{et al.} analyzed the non-linear, steady state vibration of viscoelastic beams using a fractional Zener model, where the amplitude equations were obtained using the finite element method together with the harmonic balance method, and solved using the continuation method, followed by a stability analysis.}

The sophistication of numerical methods for FDEs allowed increasing applications of fractional models in the last two decades. Here we outline some spectral methods for spatial/temporal discretization of FDEs \cite{samiee2019unified1, samiee2019unified2} and DODEs \cite{samiee2018petrov}. Among different schemes for time-fractional integration of FDEs \cite{lubich1986discretized, zayernouri2015tempered, suzuki2018automated, zayernouri2016fractionalAdams}, for simplicity, we are particularly interested in the direct L1 finite-difference (FD) scheme by Lin and Xu \cite{lin2007finite} and {we} refer the readers to \cite{Zhou2019IMEX} for a brief review of numerical methods {for time-fractional ODEs. Despite the developed works on nonlinear vibration of fractional viscoelastic beams, they} employed direct Ritz discretizations in the developed strong forms of the governing equations, which requires more smoothness to the employed basis functions. The application of spectral methods for nonlinear fractional beam models, where proper finite dimensional function spaces accounting for fractional operators are still lacking in the literature. Furthermore, from the rheology standpoint, studying the emergence of anomalous dynamics from evolving properties of {extraordinary materials}, as well as their sensitivity also require more attention. Such view is a fundamental step for physics- and mathematically-informed learning of constitutive laws from available data or desired mechanical response of the system. 

In this work, we analyze how evolving anomalous constitutive laws leads to (counter-intuitive) anomalous dynamics of mechanical systems. Our approximation of such systems is done through free- and forced- vibration response of a geometrically nonlinear Euler-Bernoulli cantilever beam with a fractional Kelvin-Voigt viscoelastic model, where:
\begin{itemize}

	\item The fractional Kelvin-Voigt model is obtained both from the Boltzmann superposition principle and a general distributed-order viscoelastic form through the choice of specific distributions of fractional orders.
	
	\item {Our framework is motivated by the influence of evolving fractal microstructures on the macroscopic dynamics of the material. {Therefore,} we study the effects of fractional orders on the response of the continuum system.}
	
	\item We {employ} Hamilton's principle to avoid {the non-trivial} decomposition of conservative (elastic) and non-conservative (viscous) {parts} of fractional constitutive laws.
	
	\item {The} weak form of {the} governing equation {is derived,} and {a single-mode approximation in space is employed, reducing the original} system to a non-linear fractional ODE.
	
	\item We perform a perturbation analysis of the resulting nonlinear fractional ODE {through the} method of multiple scales for different boundary and forcing cases.
	
	\item {Finally, we} perform a sensitivity analysis of amplitude decay rates with respect to the fractional order $\alpha$.
\end{itemize}

Our several numerical and semi-analytical experiments demonstrate a series of anomalous responses linked to far-from-equilibrium fractional dynamics, such as $\alpha$-dependent hardening-like drifts in linear amplitude-frequency behavior, long-term power law response, and super sensitivity of amplitude response with respect to fractional order values. A softening-like behavior is observed until a critical fractional order value, followed by a hardening-like response, both justified from the constitutive model standpoint. We also observe a bifurcation behavior under steady-state amplitude at primary resonance. {Such anomalous $\alpha-$dependent softening/hardening behavior motivates the notion of evolving anomalous effects, where the changing fractal material microstructure drives the fractional operator form through $\alpha$ \cite{Mashayekhi2019Fractal}. Furthermore, the observed sensitivity of the amplitude with respect to $\alpha$ could also be potentially related to the identification of early damage precursors, before the onset of macroscopic plasticity and cracks \cite{habtour2016detection}}.

This work is organized as follows. In Section \ref{Sec: Mathematical Formulation}, we derive the governing equation for the nonlinear in-plane 
vibration of a viscoelastic cantilever beam. Afterwards, we define the concepts of fractional viscoelasticity used in this work, and employ the extended 
Hamilton's principle to derive the equation of motion with external 
forces as base excitation. Then, we obtain the weak formulation of the 
problem and use assumed modes in space to reduce the problem to system of fractional ODEs. In Section  \ref{Sec: Direct Numerical Time Integration}, we obtain the corresponding linearized equation of motion and comment on the employed time-fractional integration scheme. We perform a perturbation analysis in Section \ref{Sec: Perturbatin Analysis} to 
solve the resulting nonlinear fractional ODE. A series of numerical results are presented, followed by the conclusions in Section \ref{Sec: Summery}.

\section{Mathematical Formulation}
\label{Sec: Mathematical Formulation}
We formulate our anomalous physical system and discuss the main assumptions {utilized} to derive the corresponding equation of motion. 

\subsection{Nonlinear In-Plane Vibration of a {Viscoelastic} Cantilever Beam}
\label{Sec: Fractional NL Beam}

\begin{figure}[t]
	\centering
	\begin{subfigure}{0.72\textwidth}
		\centering
		\includegraphics[width=1\linewidth]{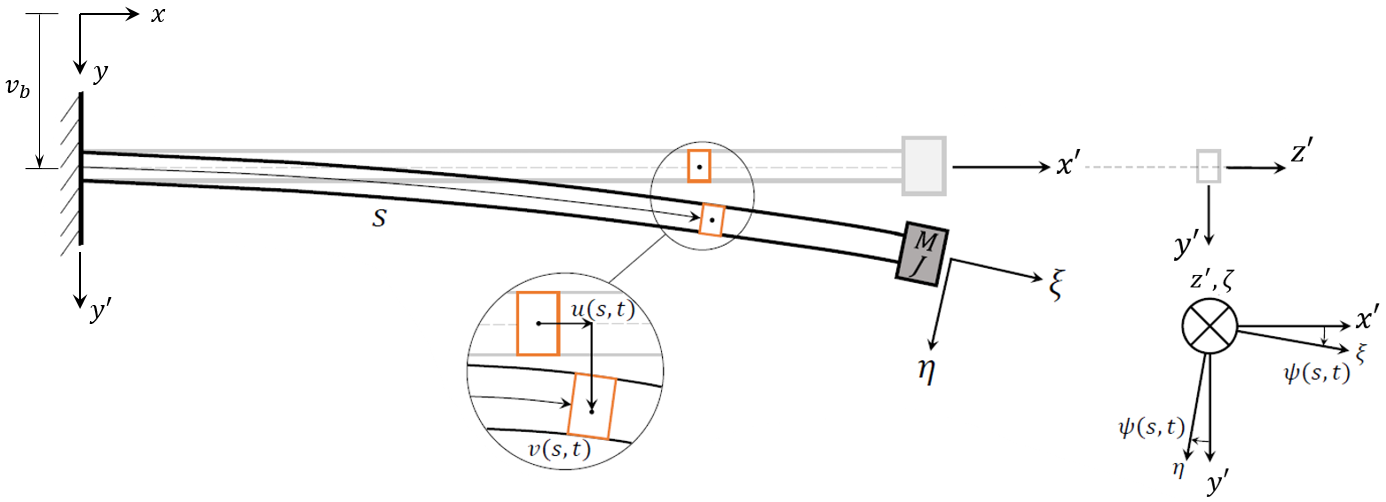}
		\caption{\label{Fig: Cantilever Beam}}
	\end{subfigure}
	\begin{subfigure}{0.27\textwidth}
		\centering
		\includegraphics[width=1\linewidth]{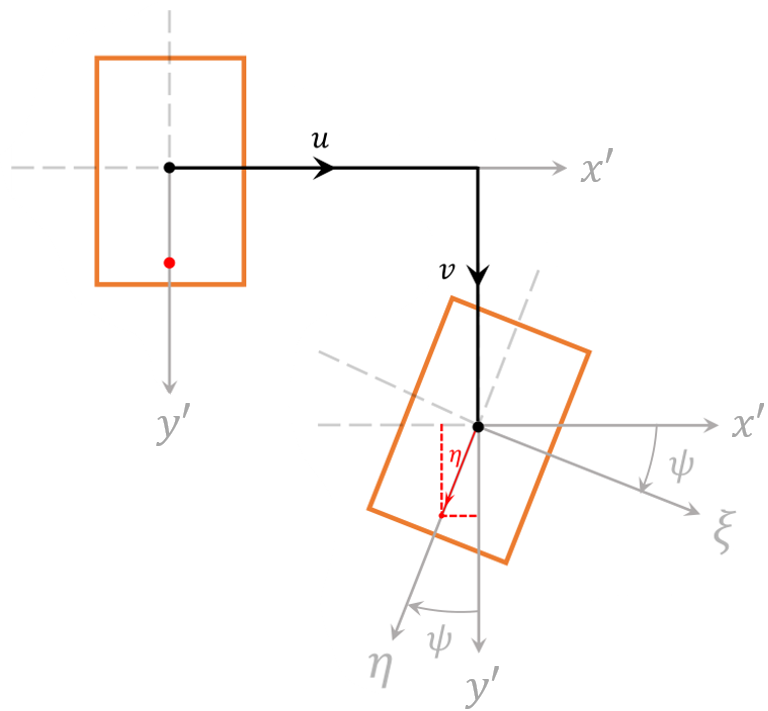}
		\caption{\label{Fig: Cantilever Beam deformation}}
	\end{subfigure}
	\caption{\textit{(a)} In-plane {kinematics of the cantilever beam subject to a base displacement $v_b(t)$ with respect to an inertial coordinate system $(x,y,z)$}. The terms $u(s,t)$ and $v(s,t)$ denote, respectively, the axial and {vertical} displacements {with respect to a $(x^\prime,y^\prime,z^\prime)$ coordinate system attached to the base}, and $\psi(s,t)$ is the rotation angle about {the} $z^\prime$-axis. \textit{(b)} The total deformation of an arbitrary {(red) point, composed of an axial displacement $u$ and vertical displacement $v$, as well as the displacement due to rotation $\psi$.} }
\end{figure}

Let the nonlinear response of a slender viscoelastic 
cantilever beam {with symmetric cross-section}, subject to harmonic {vertical} 
base excitation {denoted by $v_b$} (see Figures \ref{Fig: Cantilever Beam} and 
\ref{Fig: Cantilever Beam deformation}). We {employ} the nonlinear Euler-Bernoulli 
beam theory, where the geometric 
nonlinearities {are taken into account} in 
the equations of motion. We consider the following {kinematic and 
geometric assumptions:}
\begin{itemize}
	\item The beam is {inextensional}, i.e., {the strech along the neutral axis is 	      negligible}. The effects of warping and shear 
	deformation are ignored. Therefore, the strain {states in the cross section 
	are only due to bending}.
	
	\item The beam is slender with {symmetric} cross section{, and} undergoes purely planar flexural vibration.
	
	\item The length $L$, cross section area $A$, mass per unit length $\rho$, mass $M$ and rotatory inertia $J$ of the lumped mass at 
	the tip of beam are constant.
	
	\item The axial displacement along length of beam and the lateral 
	displacement are respectively denoted by $u(s,t)$ and $v(s,t)$.
	
	\item We consider the in-plane {vertical} vibration of the beam and reduce the problem to 1-dimension.
	
\end{itemize}
Figures \ref{Fig: Cantilever Beam} and \ref{Fig: Cantilever Beam deformation} {illustrates the kinematics of the cantilever beam under consideration. Let $(x,y,z)$ be an inertial coordinate system and $(x^\prime,y^\prime,z^\prime)$ be a moving coordinate system attached to the base of the beam, such that $(x^\prime_0,y^\prime_0,z^\prime_0) = (0,v_b,0)$}. {We note that both systems coincide when the base displacement is zero, \textit{i.e.}, $v_b(t) = 0$. Furthermore, a differential element of the beam rotates about the $z^\prime$-axis with an angle $\psi(s,t)$} to the coordinate system 
$(\xi,\eta,\zeta)$, where
{
\begin{align*}
	\begin{bmatrix} \textbf{e}_{\xi} \\ \textbf{e}_{\eta} \\ \textbf{e}_{\zeta} \end{bmatrix}
	= \begin{pmatrix}cos(\psi) & sin(\psi) & 0 \\ -sin(\psi) & cos(\psi) & 0 \\ 0 & 0 & 1\end{pmatrix}
	\begin{bmatrix} \textbf{e}_{x^\prime} \\ \textbf{e}_{y^\prime} \\ \textbf{e}_{z^\prime} \end{bmatrix},
\end{align*}}
and $\textbf{e}_{i}$ is the unit vector along the $i^{\text{\,th}}$ coordinate. 
The angular velocity and curvature at any {point $s$ along the length of 
the beam at time $t$} can be written, respectively, as
{
\begin{align}
	\label{Eq: angular velocity curvature}
	\boldsymbol{\omega}(s,t) = \frac {\partial{\psi}}{\partial{t}} \, \textbf{e}_{z^\prime}, \quad 
	\boldsymbol{\rho}(s,t) =  \frac {\partial{\psi}}{\partial{s}} \, \textbf{e}_{z^\prime}
\end{align}}
%
\noindent{Therefore, the total displacement and velocity of an arbitrary point of the beam with respect to the inertial coordinate system takes the form:}
\begin{align}
	\label{Eq: displacement}
	\textbf{r}  & = \left (u - \eta \,  \sin(\psi) \right) \, \textbf{e}_x + \left (v + {{v_b}} + \eta \, \cos(\psi)\right) \, \textbf{e}_y ,
	\\ \label{Eq: velocity}
	\frac {\partial{\textbf{r}}}{\partial{t}} & = \left( \frac {\partial{u}}{\partial{t}} - \eta \,  \frac {\partial{\psi}}{\partial{t}} \, \cos(\psi) \right) \, \textbf{e}_x 
	 + \left ( \frac {\partial{v}}{\partial{t}}  +  \frac {\partial{{{v_b}}}}{\partial{t}} - \eta \,  \frac {\partial{\psi}}{\partial{t}} \, \sin(\psi)\right) \, \textbf{e}_y .
\end{align}

\begin{figure}[t]
	\centering
	\includegraphics[width=0.55\linewidth]{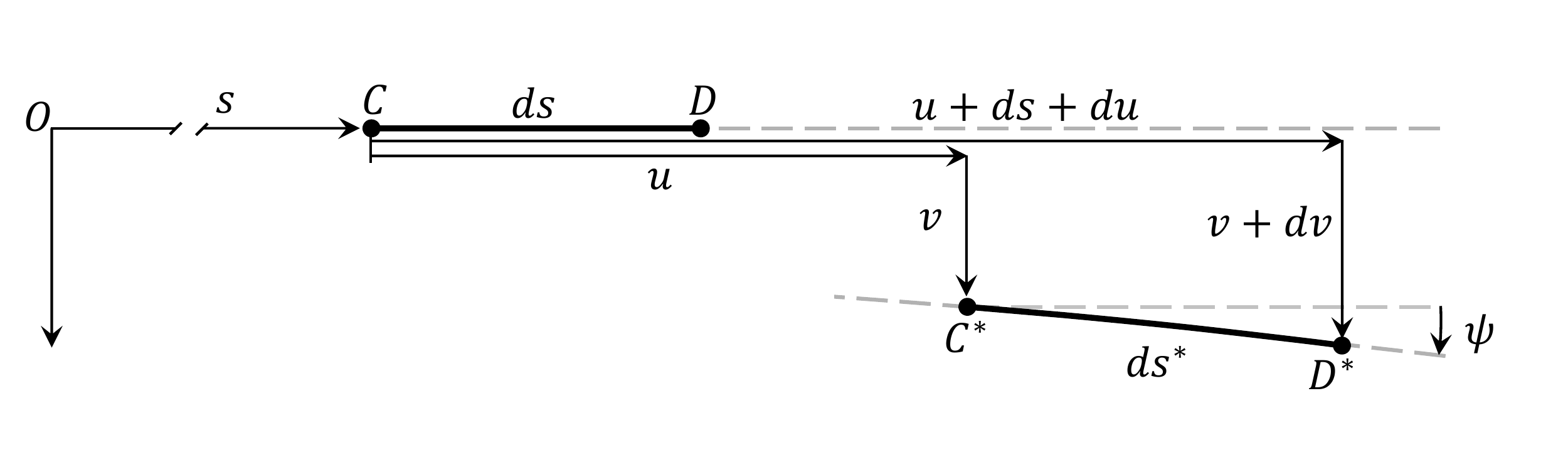}
	\caption{Deformation of an arbitrary element of the beam. The initial 
		configuration $CD$ {translates, rotates and elongates to an updated} configuration 
		$C^*D^*$.}
	\label{Fig: inextensibility}
\end{figure}

\noindent We also let an arbitrary element $CD$ {with initial length $ds$ on the neutral axis, located at a distance $s$ from the origin $O$ of the moving system $(x^\prime,y^\prime,z^\prime)$, to deformation to an updated configuration} $C^* D^*$ \textit{(see Fig. \ref{Fig: inextensibility})}. The displacement components of points $C$ and $D$ 
are denoted by the pairs $(u,v)$ and $(u+du,v+dv)$, respectively. The axial strain $e(s,t)$ at {point $C$  is} given by
\begin{align}
	\label{Eq: strain}
	e&=\frac{ds^*-ds}{ds}=\frac{\sqrt{(ds+du)^2+dv^2} - ds}{ds} =  \sqrt{(1+  \frac {\partial{u}}{\partial{s}})^2+ ( \frac {\partial{v}}{\partial{s}})^2 } - 1.
\end{align}
Applying the inextensionality {assumption}, i.e. $e = 0$, (\ref{Eq: strain}) 
becomes
\begin{align}
	\label{Eq: inextensionality}
	1 +  \frac {\partial{u}}{\partial{s}} = \left(1 - ( \frac {\partial{v}}{\partial{s}})^2\right)^{1/2}.
	%
\end{align}
Moreover, based on the assumption of {negligible vertical shear strains,} and using \eqref{Eq: inextensionality}, we have {the following expression for the rotation:}
\begin{align}
\label{Eq: no-shear}
\psi = \text{tan}^{-1}\frac{ \frac {\partial{v}}{\partial{s}}}{1+ \frac {\partial{u}}{\partial{s}}} = \text{tan}^{-1}\frac{ \frac {\partial{v}}{\partial{s}}}{\left(1 - ( \frac {\partial{v}}{\partial{s}})^2\right)^{1/2}} .
\end{align}
Using the expansion $\tan^{-1}(x) = x -\frac{1}{3}x^3 + \cdots$, the curvature can be approximated up to {third-order terms} as
\begin{align}
\label{Eq: rotation}
\psi 
& =  \frac {\partial{v}}{\partial{s}}(1 - ( \frac {\partial{v}}{\partial{s}})^2)^{-1/2} - \frac{1}{3} ( \frac {\partial{v}}{\partial{s}})^3(1 - ( \frac {\partial{v}}{\partial{s}})^2)^{-3/2} + \cdots \nonumber \\ 
& \simeq  \frac {\partial{v}}{\partial{s}}(1 + \frac{1}{2} ( \frac {\partial{v}}{\partial{s}})^2) - \frac{1}{3} ( \frac {\partial{v}}{\partial{s}})^3
\simeq  \frac {\partial{v}}{\partial{s}}+ \frac{1}{6} ( \frac {\partial{v}}{\partial{s}})^3\
\end{align}
%
\\
Therefore, the angular velocity and curvature of the beam, i.e. $ \frac {\partial{\psi}}{\partial{t}}$ and $ \frac {\partial{\psi}}{\partial{s}}$, respectively, can be approximated as:
\begin{align}
\label{Eq: angular velocity}
& \frac {\partial{\psi}}{\partial{t}} \simeq \frac {\partial^2{v}}{\partial{t}\partial{s}} + \frac{1}{2} \, \frac {\partial^2{v}}{\partial{t}\partial{s}} \, (\frac {\partial{v}}{\partial{s}} )^2 \simeq  \frac {\partial^2{v}}{\partial{t}\partial{s}}( 1 + \frac{1}{2} (\frac {\partial{v}}{\partial{s}})^2) ,
\\ \label{Eq: curvature}
&\frac {\partial{\psi}}{\partial{s}}   \simeq \frac {\partial^2{v}}{\partial{s}^2} + \frac{1}{2}\frac {\partial^2{v}}{\partial{s}^2} (\frac {\partial{v}}{\partial{s}})^2 \simeq \frac {\partial^2{v}}{\partial{s}^2}(1 + \frac{1}{2} (\frac {\partial{v}}{\partial{s}})^2) .
\end{align}

By the Euler-Bernoulli beam assumptions a slender {beam without vertical shear strains}, the strain-curvature {relationship} takes the form
\begin{align}
\label{Eq: strain-curvature relation}
\varepsilon(s,t) = 
-\eta \,\frac {\partial{\psi}(s,t)}{\partial{s}}
\end{align}
%
\subsection{Linear Viscoelasticity: Boltzmann Superposition Principle}
\label{Sec: Viscoelasticity}
In this section, we start with a \textit{bottom-up} derivation of our rheological 
building block, \textit{i.e.}, the Scott-Blair model through the Boltzmann 
superposition principle. Then, in a \textit{top-bottom} fashion, we demonstrate 
how the fractional Kelvin-Voigt model is obtained from a general 
distributed-order form. Assuming linear viscoelasticity, and 
applying a small step strain increase, denoted by $\delta \varepsilon(t)$, at a 
given time $t=\tau_1$, the resulting stress in the material is given by: 
\begin{align}
\label{Eq: relaxation function}
\sigma(t) = G(t-\tau_1) \delta\varepsilon(\tau_1), \quad t > 
\tau_1,
\end{align}
where $G(t)$ denotes the relaxation function. The Boltzmann superposition 
principle states that resulting stresses from distinct applied small strains 
are additive. Therefore, the total tensile stress of the specimen at time $t$ 
is obtained from the superposition of infinitesimal changes in strain at some 
prior time $\tau_j$, given as $G(t-\tau_j) \delta \varepsilon(\tau_j)$. 
Therefore,
\begin{align}
\label{Eq: tensile stress - 1}
\sigma ( t ) 
= \sum_{\tau_j<t} G(t-\tau_j) \frac{\delta \varepsilon(\tau_j)}{\delta 
\tau_j}\delta\tau_j,
\end{align}
where the limiting case $\delta \tau_j \rightarrow 0$ yields the following 
integral form:
\begin{align}
\label{Eq: tensile stress - 2}
\sigma ( t ) = \int_{-\infty}^{t} G(t - \tau) \, \dot{\varepsilon}(\tau)\, 
d\tau,
\end{align}
where $\dot{\varepsilon}$ denotes the strain rate.

\subsubsection{Exponential Relaxation (Classical Models) vs. Power-Law Relaxation (Fractional Models)}

\begin{figure}[t]
	\centering
	\includegraphics[width=0.45\linewidth]{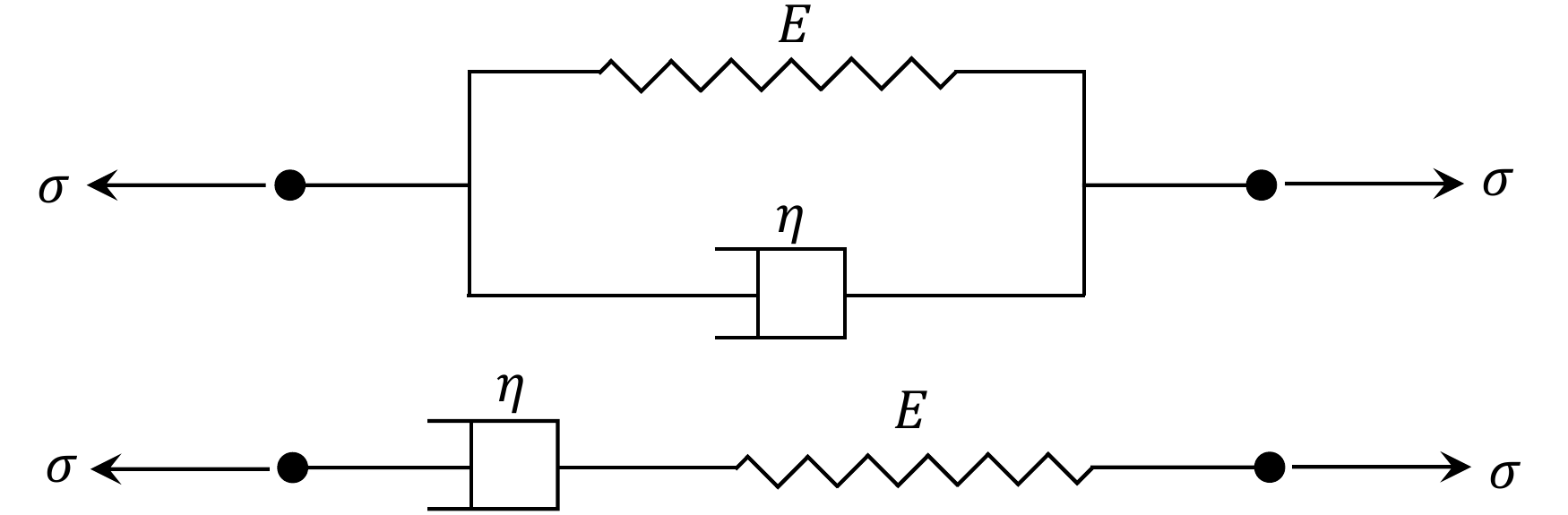}
	\caption{Classical viscoelastic models as a combination of spring (purely elastic) and dash-pot (purely viscous) elements. Kelvin-Voigt (top) and Maxwell (bottom) rheological models.}
	\label{Fig: Maxwell}
\end{figure}

The relaxation function $G(t)$ is traditionally expressed as the summation of 
exponential functions with different exponents and constants, which yields the 
\textit{so-called} generalized Maxwell form as: 
\begin{align}
\label{Eq: exponential relaxation}
G(t) = \sum C_i e^{- t/\tau_i}.
\end{align}
For the simple case of a single exponential term (a single Maxwell branch), we 
have $G(t) = E e^{- t/\tau}$. Therefore, in the case of zero 
initial strain 
$(\varepsilon(0) = 0)$, we have:
\begin{align}
\sigma ( t ) 
= E \int_{0}^{t}  e^{- (t - \tilde{t})/\tau} \, \dot{\varepsilon}(\tilde{t})
\, d\tilde{t},
\end{align}
which solves the integer-order differential equation $ 
\frac{\partial{\varepsilon}}{\partial{t}} = \frac{1}{E}  
\frac{\partial{\sigma}}{\partial{t}} + \frac{1}{\eta} \sigma$, where the 
relaxation time constant $\tau = \eta / E$ is obtained from experimental 
observations. The Maxwell model is in fact a combination of purely elastic and 
purely viscous elements in series, as illustrated in Fig. \ref{Fig: Maxwell}.

By letting the relaxation function (kernel) 
in \eqref{Eq: tensile stress - 2} have a modulated power-law form 
$G(t) = 
E_\alpha\, g(\alpha) ( \, t-\tau \, )^{-\alpha}$, equation (\ref{Eq: tensile 
stress - 2}) for the stress takes the following form,
\begin{align}
\label{Eq: fractional stress}
\sigma ( t ) 
&= E_\alpha \, g(\alpha) \,  \int_{-\infty}^{t}  
\frac{\dot{\varepsilon}(\tau)}{( \, t-\tau \, )^{\alpha}} \, d\tau.
\end{align}
where $E_\alpha$ denotes a \textit{pseudo-constant} with units $[Pa.s^\alpha]$.
If we choose the modulation $g(\alpha) = \frac{1}{\Gamma(1-\alpha)}$, then the 
integro-differential operator \eqref{Eq: fractional stress} gives the 
Liouville-Weyl fractional derivative \cite{mainardi2010fractional}. Although 
the lower integration limit of (\ref{Eq: fractional stress}) is taken as 
$-\infty$, under hypothesis of causal histories, which states that the 
viscoelastic body is quiescent for all time prior to some starting point $t=0$, 
\eqref{Eq: fractional stress} can be re-written as 
\begin{align}
\label{Eq: fractional stress - 2}
\sigma ( t ) 
& = \varepsilon(0^{+}) \, \frac{E_\alpha \, g(\alpha)}{t ^{\alpha}} + E_\alpha 
\, g(\alpha) \,  \int_{0}^{t}  \frac{\dot{\varepsilon}}{( \, t-\tau \, 
)^{\alpha}} \, d\tau,
 \nonumber \\
& = \varepsilon(0^{+}) \, \frac{E_\alpha \, g(\alpha)}{t ^{\alpha}} + E_\alpha 
\, \prescript{C}{0}{\mathcal D}_{t}^{\alpha} \, \varepsilon,
 \nonumber \\
& = E_\alpha\, \prescript{RL}{0}{\mathcal D}_{t}^{\alpha} \, \varepsilon,
\end{align}
where $\prescript{C}{0}{\mathcal D}_{t}^{\alpha}$ and 
$\prescript{RL}{0}{\mathcal D}_{t}^{\alpha}$ denote, respectively, the Caputo 
and Riemann-Liouville fractional derivatives \cite{mainardi2010fractional}. 
Both definitions are equivalent here due to homogeneous initial conditions for 
the strain. 

\begin{remark}	\label{Rem: Scott Blair}
	The constitutive equation \eqref{Eq: fractional stress - 2} can be thought of as an interpolation between a pure elastic (spring) and a pure viscous (dash-pot) elements, i.e., the Scott Blair element \cite{rogosin2014george, mainardi2008time, mainardi2011creep, suzuki2016fractional}. It should be noted that in the limiting cases of $\alpha \rightarrow 0$ and $\alpha \rightarrow 1$, the relation \eqref{Eq: fractional stress - 2} recovers the corresponding equations for spring and dash-pot, respectively.
\end{remark}

\subsubsection{Multi-Scale Power-Laws, Distributed-Order Models}

In the most general sense, materials intrinsically possess a spectrum of 
power-law relaxations, and therefore we need a distributed-order representation 
for the stress-strain relationship. Consequently, the relaxation function 
$G(t)$ in \eqref{Eq: tensile stress - 2} does not only contain a single power-law 
as in \eqref{Eq: fractional stress}, but rather a distribution over a 
range of 
values. Considering nonlinear 
viscoelasticity with material heterogeneities, the distributed order 
constitutive equations over $t > 0$ with orders $\alpha \in [\alpha_{min}, 
\alpha_{max}]$ and $\beta \in [\beta_{min}, \beta_{max}]$ can be expressed in 
the general form as
\begin{align}
\label{Eq: Dis const law}
&\int_{\beta_{min}}^{\beta_{max}} \Phi(\beta; x, t, \sigma) 
\prescript{*}{0}{\mathcal D}_{t}^{\beta} \sigma(t) \, d\beta = \int_{\alpha_{min}}^{\alpha_{max}}  \Psi(\alpha; x, t, \varepsilon)
\prescript{*}{0}{\mathcal D}_{t}^{\alpha} \varepsilon(t) d\alpha ,
\end{align}
in which the prescript $*$ stands for any type of fractional derivative, {and initial conditions also depending on such definitions.} The 
functions $\Phi(\beta; x, t, \sigma)$ and $\Psi(\alpha; x, t, \varepsilon)$ can 
be thought of as distribution functions, where $\alpha \mapsto \Psi(\alpha; x, 
t, \varepsilon)$ 
and $\beta \mapsto \Phi(\beta; x, t, \sigma)$  are continuous mappings in 
$[\alpha_{min} , 
\alpha_{max}]$ and $[\beta_{min} , \beta_{max}]$. Furthermore, the dependence 
of the distributions on the (thermodynamically) conjugate pair $(\sigma, 
\varepsilon)$ introduces the notion of nonlinear viscoelasticity, and the 
dependence on a material coordinate $x$ induces material heterogeneties in 
space.

\begin{remark}
	The pairs ($\alpha_{min}$, $\alpha_{max}$) and ($\beta_{min}$, $\beta_{max}$) 
	are only the theoretical lower and upper terminals in the definition of 
	distributed order models. In general, the distribution function $\Phi(\beta; x, 
	t, \sigma)$ and $\Psi(\alpha; x, t, \varepsilon)$ can arbitrarily confine the 
	domain of integration in each realization of practical rheological problems and 
	material design. If we let the distribution be summation of some delta 
	functions, then, the distributed order model becomes the following multi-term 
	model: 
	\begin{align*}
	\left( 1+ \sum_{k=1}^{p_{\sigma}} \, a_k \, \prescript{}{0}{\mathcal D}_{t}^{\beta_k} \right) \, \sigma(t)
	=
	\left( c + \sum_{k=1}^{p_{\varepsilon}} \, b_k \, \prescript{}{0}{\mathcal D}_{t}^{\alpha_k} \right) \, \varepsilon(t).
	\end{align*}
\end{remark}

In order to obtain the fractional Kelvin-Voigt model, we let $\Phi(\beta) = 
\delta(\beta)$ and $\Psi(\alpha) = E_{\infty} \delta(\alpha) + E_{\alpha} 
\delta(\alpha-\alpha_0)$ in \eqref{Eq: Dis const law}, and therefore,
\begin{align}
\label{Eq: frac KV}
\sigma(t)
=
E_{\infty} \, \varepsilon(t) + E_{\alpha} \, \prescript{RL}{0}{\mathcal D}_{t}^{\alpha} \, \varepsilon(t), \quad  { \alpha \in (0,1).}
\end{align}

\subsection{Extended Hamilton's Principle}
%
We derive the equations of motion by employing the extended Hamilton's principle $$\int_{t_1}^{t_2} \left (\delta T - \delta W\right) \, dt = 0,$$ where $\delta T$ and $\delta W$ {denote} the variations of kinetic energy and total work \cite{meirovitch2010fundamentals}. The only source of external input to our system of interest is the base excitation, {which is linearly superposed to the beam's vertical displacement $v(t)$, and therefore contributes to the kinetic energy taken in the inertial (Lagrangian) coordinate system}. Hence, the total work only {contains} the internal work done by {the stress state, with the variation expressed} as \cite{bonet1997nonlinear}
\begin{align}
\label{Eq: work var}
\delta W = \int_{\mathbb{V}} \sigma \, \delta \varepsilon \, dv,
\end{align}
where the integral is taken over the whole system volume $\mathbb{V}$. 

\begin{remark}
	It is remarked in (\ref{Rem: Scott Blair}) that the fractional {Scott-Blair elements} exhibit 
	both elasticity and viscosity {behaviors}. There have been {attempts in the literature 
	to separate the conservative (elastic) and non-conservative (viscous) parts of 
	fractional constitutive equations at the free-energy level} \cite{lion1997thermodynamics}. {However,} we note this separation in the time domain is not trivial for sophisticated fractional constitutive equations, {and therefore we choose to formulate our problem in terms of the total work in order to avoid such additional complexities.}
\end{remark}

The full derivation of the governing equation using the extended Hamilton's principle is given in \ref{Sec: App. extended Hamilton}.\ We recall that $M$ and $J$ are the mass and rotatory inertia of the lumped mass at the tip of beam, $\rho$ is the mass per unit length of the beam, $I = \int_{A} \eta^2 \, dA$, and let $m = \frac{\rho}{ E_{\infty} \, I}$ and $E_r = \frac{E_{\alpha}}{E_{\infty}}$. We approximate the nonlinear terms up to third order and use the following dimensionless variables
\begin{align*}
s^{*} = \frac{s}{L}, \,\, 
v^{*} &= \frac{v}{L}, \,\,
t^{*} = t \left(\frac{1}{m L^4}\right)^{1/2}, \,\, 
E_r^{*} = E_r \left(\frac{1}{m L^4}\right)^{\alpha/2}, \,\, J^{*} = \frac{J}{\rho L^3}, \,\,
M^{*} = \frac{M}{\rho L}, \,\,
{{v_b}}^* =  \frac{{{v_b}}}{L},
\end{align*}
and derive the strong form of the equation of motion. Therefore, by choosing a proper function space $V$, the problem reads as: find $v \in V$ such that
\begin{align}
\label{Eq: eqn of motion - dimless - body}
&\frac{\partial^2{v}}{\partial{t}^2} 
+ 
\frac{\partial^2}{\partial{s}^2} \left(  
\frac{\partial^2{v}}{\partial{s}^2}  + \frac{\partial^2{v}}{\partial{s}^2}  (\frac{\partial{v}}{\partial{s}} )^2 \right. 
\left.+ E_r \prescript{RL}{0}{\mathcal D}_{t}^{\alpha} \, \Big[\frac{\partial^2{v}}{\partial{s}^2}(1 + \frac{1}{2} (\frac{\partial{v}}{\partial{s}} )^2)\Big]
+ \frac{1}{2}  E_r (\frac{\partial{v}}{\partial{s}} )^2 \, \prescript{RL}{0}{\mathcal D}_{t}^{\alpha} \, \frac{\partial^2{v}}{\partial{s}^2} 
\right) \nonumber\\
& - 
\frac{\partial}{\partial{s}}\left( 
\frac{\partial{v}}{\partial{s}} \, (\frac{\partial^2{v}}{\partial{s}^2})^2  
+ E_r \frac{\partial{v}}{\partial{s}} \, \frac{\partial^2{v}}{\partial{s}^2} \,  \prescript{RL}{0}{\mathcal D}_{t}^{\alpha} \,  \frac{\partial^2{v}}{\partial{s}^2}
\right)
= - \ddot{{{v_b}}},
\end{align}
which is subject to the following boundary conditions:

\begin{align}
\label{Eq: bc - dimless - body}
& v \, \Big|_{s=0} = \frac{\partial{v}}{\partial{s}} \, \Big|_{s=0} = 0 ,
 \nonumber\\
& 
J \left(  \frac{\partial^3{v}}{\partial{t}^2\partial{s}} ( 1 + ( \frac{\partial{v}}{\partial{s}})^2) +  \frac{\partial{v}}{\partial{s}}  (\frac{\partial^2{v}}{\partial{s}\partial{t}})^2  \right) \nonumber \\
& +\left( \frac{\partial^2{v}}{\partial{s}^2}+  \frac{\partial^2{v}}{\partial{s}^2} ( \frac{\partial{v}}{\partial{s}})^2
+ E_r \, \prescript{RL}{0}{\mathcal D}_{t}^{\alpha} \,  \Big[\frac{\partial^2{v}}{\partial{s}^2}(1 + \frac{1}{2}  (\frac{\partial{v}}{\partial{s}})^2) \Big] \right.
\left. +\frac{1}{2}  E_r \,  (\frac{\partial{v}}{\partial{s}})^2 \, \prescript{RL}{0}{\mathcal D}_{t}^{\alpha} \,  \frac{\partial^2{v}}{\partial{s}^2}
\right) \,\, \Bigg|_{s=1} = 0 ,
\nonumber\\
& 
M ( \frac{\partial^2{v}}{\partial{t}^2}\!+\! \ddot{{{v_b}}})
\!-\! \frac{\partial{v}}{\partial{s}}\left(
\frac{\partial^2{v}}{\partial{s}^2}\!+\! \frac{\partial^2{v}}{\partial{s}^2} (\frac{\partial{v}}{\partial{s}})^2
 +E_r \, \prescript{RL}{0}{\mathcal D}_{t}^{\alpha} \, \Big[\frac{\partial^2{v}}{\partial{s}^2}(1 \!+\! \frac{1}{2} (\frac{\partial{v}}{\partial{s}})^2)\Big]
\!+\! \frac{1}{2}  E_r \, (\frac{\partial{v}}{\partial{s}})^2 \, \prescript{RL}{0}{\mathcal D}_{t}^{\alpha} \, \frac{\partial^2{v}}{\partial{s}^2}
\right)
\nonumber \\
& 
+\left( 
\frac{\partial{v}}{\partial{s}} \, (\frac{\partial^2{v}}{\partial{s}^2})^2  + E_r \, \frac{\partial{v}}{\partial{s}} \, \frac{\partial^2{v}}{\partial{s}^2}  \, \prescript{RL}{0}{\mathcal D}_{t}^{\alpha} \, \frac{\partial^2{v}}{\partial{s}^2} 
\right) \,\, \Bigg|_{s=1}  = 0.
\end{align}
%
\subsection{Weak Formulation}
\label{Sec: Weak Formulation}
%
The common practice in analysis of numerical methods for PDEs are mostly concerned with linear equations. {Analyses} for linear PDEs are well-developed and well-defined, however {they are still scarce for nonlinear PDEs}. The linear theories are usually applicable to nonlinear problems if the solution is sufficiently smooth \cite{tadmor2012review}. We do not intend to investigate/develop analysis for our proposed nonlinear model. Instead, by assuming smooth solution, we employ linear theories in our analysis. Let $v: \mathbb{R}^{1+1}\rightarrow\mathbb{R}$ for $\alpha \in (0,1)$ and $\Omega = [0,T]\times [0,L]$. Here, we construct the solution space, $\mathcal{B}^{\alpha} \, (\Omega)$, endowed with proper norms \cite{samiee2019unified2}, in which the corresponding weak form of \eqref{Eq: eqn of motion - dimless - body} can be formulated. If we recall the equation \eqref{Eq: eqn of motion - dimless - body} as E, then:
\begin{align}
\mathcal{B}^{\alpha} \, (\Omega):=\Big\{v\in\prescript{l}{0}{H}^{\alpha}(\Omega) \, \Big|  \int_{\Omega} E \,d\Omega< \infty \Big\}
\end{align}
where
\begin{align*}
\prescript{l}{0}{H}^{\alpha}(\Omega) =\prescript{l}{0}{H}^{\alpha}\,\Big(I;L^2\,(\Omega)\Big) \cap L^2(I ;\prescript{}{0}{H}^{2}(\Omega)), \qquad 
\prescript{}{0}{H}^{2}(\Omega)=\Big\{v\in{H}^{2}(\Omega) \, \Big| \, \,v \, \Big|_{s=0} = \frac{\partial{v}}{\partial{s}} \, \Big|_{s=0} = 0  \Big\}
\end{align*}

We obtain the weak form of the problem by multiplying the strong form \eqref{Eq: eqn of motion - dimless - body} with proper test functions $\tilde{v}(s) \in \mathcal {B}^{\alpha} (\Omega)$ and integrating over the dimensionless spatial computational domain { $\Omega_s = [0,1]$. The test functions satisfy the} boundary conditions, i.e. $\tilde{v}(0) = \frac{\partial{\tilde{v}}}{\partial{s}}(0) = 0 $. Therefore, { we obtain:}
\begin{align}
\label{Eq: weak form - 1}
&\! \int_{0}^{1}\! \frac{\partial^2{v}}{\partial{t}^2}  \tilde{v} ds  
\!+\! \!\int_{0}^{1}\! 
\frac{\partial^2}{\partial{s}^2}\left(  
\frac{\partial^2{v}}{\partial{s}^2} \!+\! \frac{\partial^2{v}}{\partial{s}^2}(\frac{\partial{v}}{\partial{s}})^2 
\!+\! E_r \prescript{RL}{0}{\mathcal D}_{t}^{\alpha}  \Big[\frac{\partial^2{v}}{\partial{s}^2}(1 \!+\! \frac{1}{2} (\frac{\partial{v}}{\partial{s}})^2)\Big]\!+\! \frac{1}{2}  E_r 
(\frac{\partial{v}}{\partial{s}})^2  \prescript{RL}{0}{\mathcal D}_{t}^{\alpha}  \frac{\partial^2{v}}{\partial{s}^2}\right)
\tilde{v}  ds
\nonumber \\
&
- \int_{0}^{1} 
\frac{\partial}{\partial{s}}\left( 
\frac{\partial{v}}{\partial{s}} \,  (\frac{\partial^2{v}}{\partial{s}^2})^2  
+ E_r
\frac{\partial{v}}{\partial{s}} \, \frac{\partial^2{v}}{\partial{s}^2} \,  \prescript{RL}{0}{\mathcal D}_{t}^{\alpha} \, \frac{\partial^2{v}}{\partial{s}^2}
\right)
\, \tilde{v} \, ds  = - \int_{0}^{1} \ddot{{{v_b}}} \, \tilde{v} \, ds.
\end{align}
{ Integrating the above equation by parts, we obtain:}
\begin{align}
\label{Eq: weak form - 2}
& \frac{\partial^2}{\partial {t}^2}  \int_{0}^{1} v \tilde{v}  ds  
\!+\! 
\int_{0}^{1} 
\left(  
\frac{\partial^2{v}}{\partial{s}^2}+ \frac{\partial^2{v}}{\partial{s}^2} (\frac{\partial{v}}{\partial{s}})^2 
\!+\! E_r
\prescript{RL}{0}{\mathcal D}_{t}^{\alpha}  \Big[\frac{\partial^2{v}}{\partial{s}^2}(1 \!+\! \frac{1}{2} (\frac{\partial{v}}{\partial{s}})^2)\Big]  
\!+\! \frac{1}{2}  E_r 
(\frac{\partial{v}}{\partial{s}})^2  \prescript{RL}{0}{\mathcal D}_{t}^{\alpha}  \frac{\partial^2{v}}{\partial{s}^2}
\right) \frac{\partial^2{\tilde{v}}}{\partial{s}^2}  ds
\nonumber \\
&
+
\int_{0}^{1} 
\left( 
\frac{\partial{v}}{\partial{s}} \,  (\frac{\partial^2{v}}{\partial{s}^2})^2  
+ E_r
\frac{\partial{v}}{\partial{s}} \, \frac{\partial^2{v}}{\partial{s}^2}\,  \prescript{RL}{0}{\mathcal D}_{t}^{\alpha} \, \frac{\partial^2{v}}{\partial{s}^2}
\right)
\, \frac{\partial{\tilde{v}}}{\partial{s}} \, ds
+ M (\frac{\partial^2{v}}{\partial{t}^2} + \ddot{{{v_b}}}) \,\tilde{v} \, \Big|_{s=1}
+ J \left(  \frac{\partial^3{v}}{\partial{t}^2\partial{s}}( 1 + (\frac{\partial{v}}{\partial{s}})^2) + \frac{\partial{v}}{\partial{s}} (\frac{\partial^2{v}}{\partial{t}\partial{s}})^2  \right) \frac{\partial{\tilde{v}}}{\partial{s}} \, \Big|_{s=1}  \nonumber \\
&
= - \ddot{{{v_b}}} \, \int_{0}^{1} \tilde{v} \, ds .
\end{align}
By rearranging the terms, we get
\begin{align}
\label{Eq: weak form - 3}
& \frac{\partial^2}{\partial {t}^2} \, \left( \int_{0}^{1} v \, \tilde{v} \, ds +  M \, v \, \tilde{v} \, \Big|_{s=1} + J \, \frac{\partial{v}}{\partial{s}} \, \frac{\partial{\tilde{v}}}{\partial{s}} \, \Big|_{s=1}  \right) 
+
J \left(  \frac{\partial^3{v}}{\partial{t}^2\partial{s}} (\frac{\partial{v}}{\partial{s}})^2 + \frac{\partial{v}}{\partial{s}} \ (\frac{\partial^2{v}}{\partial{t}\partial{s}})^2 \right) \frac{\partial{\tilde{v}}}{\partial{s}} \, \Big|_{s=1}
+\int_{0}^{1} \frac{\partial^2{v}}{\partial{s}^2}\, \frac{\partial^2{\tilde{v}}}{\partial{s}^2} \, ds
 \nonumber \\
&
+ E_r 
\int_{0}^{1} \prescript{RL}{0}{\mathcal D}_{t}^{\alpha} \, \Big[\frac{\partial^2{v}}{\partial{s}^2}\Big]\,\, \frac{\partial^2{\tilde{v}}}{\partial{s}^2} \, ds+ 
\int_{0}^{1} \frac{\partial^2{v}}{\partial{s}^2}(\frac{\partial{v}}{\partial{s}})^2 \,\, \frac{\partial^2{\tilde{v}}}{\partial{s}^2} \, ds
+
\int_{0}^{1} \frac{\partial{v}}{\partial{s}} \,  (\frac{\partial^2{v}}{\partial{s}^2})^2  \,\, \frac{\partial{\tilde{v}}}{\partial{s}} \, ds
+ \frac{E_r}{2}
\int_{0}^{1} \prescript{RL}{0}{\mathcal D}_{t}^{\alpha} \, \Big[\frac{\partial^2{v}}{\partial{s}^2}(\frac{\partial{v}}{\partial{s}})^2\Big] \,\, \frac{\partial^2{\tilde{v}}}{\partial{s}^2} \, ds
 \nonumber \\
& 
\!+\! \frac{E_r}{2}  
\int_{0}^{1} (\frac{\partial{v}}{\partial{s}})^2  \prescript{RL}{0}{\mathcal D}_{t}^{\alpha} \Big[\frac{\partial^2{v}}{\partial{s}^2}\Big] \frac{\partial^2{\tilde{v}}}{\partial{s}^2}  ds
\!+\! E_r
\int_{0}^{1} \frac{\partial{v}}{\partial{s}} \frac{\partial^2{v}}{\partial{s}^2} \prescript{RL}{0}{\mathcal D}_{t}^{\alpha}  \Big[\frac{\partial^2{v}}{\partial{s}^2}\Big] \frac{\partial{\tilde{v}}}{\partial{s}}  ds
= - \ddot{{{v_b}}} \left( \int_{0}^{1} \tilde{v}  ds \!+\! M  \tilde{v}  \Big|_{s=1} \right).
\end{align}
%
\subsection{Assumed Mode: {A Spectral} Approximation in Space}
\label{Sec: Spectral Galerkin Method}
%
We employ the following modal discretization to obtain a reduced-order model of the beam. Therefore,
\begin{align}
\label{Eq: assumed mode}
v(s,t) \simeq v_N(s,t) =  \sum_{n=1}^{N} q_n(t) \, \phi_n(s),
\end{align}
where the spatial functions $\phi_n(s), \,\, n=1,2,\cdots,N$ are assumed \textit{a priori} and the temporal functions $q_n(t), \,\, n=1,2,\cdots,N$ are the unknown modal coordinates. The assumed modes $\phi_n(s)$ in \eqref{Eq: assumed mode} are obtained in \ref{Sec: App. Eigenvalue Problem of Linear Model},  by solving the corresponding linear eigenvalue problem of our nonlinear model. Subsequently, we construct the proper finite dimensional spaces of basis/test functions as:
\begin{align}
\label{Eq: Solution/Test Space}
V_N = \tilde{V}_N =  span \, \Big\{ \,\, \phi_n(x) : n = 1,2, \cdots, N \,\, \Big\}.
\end{align}
Since $V_N = \tilde{V}_N \subset V = \tilde{V} $, problem \eqref{Eq: weak form - 3} read as: find $v_N \in V_N$ such that
\begin{align}
\label{Eq: weak form - discrete}
& \frac{\partial^2}{\partial {t}^2}  \left( \int_{0}^{1} v_N  \tilde{v}_N  ds \!+ \! M  v_N  \tilde{v}_N  \Big|_{s=1} \!+\! J  \frac{\partial{{v}_N}}{\partial{s}}  \frac{\partial{\tilde{v}_N}}{\partial{s}} \, \Big|_{s=1}  \right) 
\!+\! 
J \left(  \frac{\partial^3{{v}_N}}{\partial{t}^2\partial{s}} (\frac{\partial{{v}_N}}{\partial{s}})^2 \!+\! \frac{\partial{{v}_N}}{\partial{s}}(\frac{\partial^2{{v}_N}}{\partial{t}\partial{s}})^2  \right) \frac{\partial{\tilde{v}_N}}{\partial{s}}  \Big|_{s=1}
+ 
\int_{0}^{1} \frac{\partial^2{{v}_N}}{\partial{s}^2} \, \frac{\partial^2{\tilde{v}_N}}{\partial{s}^2} \, ds
 \nonumber\\
&
+ E_r 
\int_{0}^{1} \prescript{RL}{0}{\mathcal D}_{t}^{\alpha} \, \Big[\frac{\partial^2{{v}_N}}{\partial{s}^2}\Big] \,\,  \frac{\partial^2{\tilde{v}_N}}{\partial{s}^2} \, ds
+ 
\int_{0}^{1} \frac{\partial^2{{v}_N}}{\partial{s}^2} (\frac{\partial{{v}_N}}{\partial{s}})^2 \,\,  \frac{\partial^2{\tilde{v}_N}}{\partial{s}^2} \, ds
+
\int_{0}^{1} \frac{\partial{{v}_N}}{\partial{s}} \,  (\frac{\partial^2{{v}_N}}{\partial{s}})^2  \,\, \frac{\partial{\tilde{v}_N}}{\partial{s}}  \, ds
\nonumber \\
&
+ \frac{E_r}{2}
\int_{0}^{1} \prescript{RL}{0}{\mathcal D}_{t}^{\alpha} \, \Big[\frac{\partial^2{{v}_N}}{\partial{s}^2} (\frac{\partial{v_N}}{\partial{s}} )^2\Big] \,\, \frac{\partial^2{\tilde{v}_N}}{\partial{s}^2}\, ds
+ \frac{E_r}{2} 
\int_{0}^{1} (\frac{\partial{{v}_N}}{\partial{s}})^2 \, \prescript{RL}{0}{\mathcal D}_{t}^{\alpha} \, \Big[\frac{\partial^2{{v}_N}}{\partial{s}^2}\Big] \,\, \frac{\partial^2{\tilde{v}_N}}{\partial{s}^2}  \, ds 
\nonumber \\
&
+ E_r
\int_{0}^{1} \frac{\partial{{v}_N}}{\partial{s}} \, \frac{\partial^2{v_{N}}}{\partial{s}^2} \,  \prescript{RL}{0}{\mathcal D}_{t}^{\alpha} \, \Big[\frac{\partial^2{{v}_N}}{\partial{s}^2}\Big] \,\, \frac{\partial{\tilde{v}_N}}{\partial{s}}  \, ds
= -\ddot{{{v_b}}} \left( \int_{0}^{1} \tilde{v}_N \, ds + M \, \tilde{v}_N \, \Big|_{s=1} \right),
\end{align}
for all $ \tilde{v}_N \in \tilde{V}_N$. 

%
\subsection{Single Mode Approximation} \label{subsection}
%
In general, the modal discretization \eqref{Eq: assumed mode} in \eqref{Eq: weak form - discrete} leads to {a} coupled nonlinear system of fractional ordinary differential equations. We note that while the fractional operators already impose numerical challenges, {these are increased by the presence of nonlinearities, leading to failure of existing numerical schemes}. However, without loss of generality, we can assume that only one mode (primary mode) of motion is involved in the dynamics of system of interest. 

\subsubsection{Why is single-mode approximation useful?}

Although single-mode approximations are simplistic in nature, they encapsulate the most fundamental dynamics and the highest energy mode in the motion of nonlinear systems. Furthermore, as shown by numerous studies below, such approximation also proved capable of capturing the complex behavior of structures.

Azrar et al. \cite{azrar1999semi} {demonstrated sufficient approximations of single-} and multi-modal representation {for} the nonlinear forced vibration of a simply supported beam under a uniform harmonic distributed force. Tseng and Dugundji \cite{tseng1971nonlinear} showed {similar results between} single and two mode approximations for nonlinear vibrations of clamped-clamped beams {far} from the crossover region. Loutridis et al. \cite{loutridis2005forced} implemented a crack detection {method for} beams using a single-degree-of-freedom system with time varying stiffness. In \cite{hamdan1997large}, the effects of base stiffness and attached mass on the nonlinear{, planar flexural free vibrations of beams were studied}. Lestari and Hanagud \cite{lestari2001nonlinear} studied the nonlinear free vibrations of buckled beams with elastic end constraints, {where the single-mode assumption led to a closed-form solution in terms of elliptic functions}. 
	
{Of particular interest, Habtour \textit{et al.} \cite{habtour2016detection} detected and validated the response of a nonlinear cantilever beam subject to softening due to local stress-induced, early fatigue damage precursors \textit{prior} to crack formation. Their findings demonstrate that the pragmatism of a single-mode approximation provides sufficient sensitivity of the amplitude response with respect to the nonlinear stiffness, making their framework an effective practical tool for early fatigue detection. We also refer the reader to \cite{eisley1964nonlinear, hsu1960application, pillai1992nonlinear, evensen1968nonlinear} for additional applications.} 

Therefore, we let the anomalous dynamics of our system be driven by the fractional-order $\alpha$, and following the aforementioned studies, we replace \eqref{Eq: assumed mode} with the one-mode discretization $v_N = q(t) \, \phi(s)$ (where we let $N=1$ and drop the subscript $1$ for simplicity). Upon substituting in \eqref{Eq: weak form - discrete}, we obtain the unimodal governing equation of motion as (see \ref{Sec: App. Single Mode Decomposition}),
\begin{align}
\label{Eq: weak form - discrete 2}
&\mathcal{M}  \ddot{q} + \mathcal{J} (\ddot q  q^2 \!+\! q  {\dot q}^2)
+ \mathcal{K}_l  q + E_r  \mathcal{C}_l  \prescript{RL}{0}{\mathcal D}_{t}^{\alpha} q  
+ 2 \mathcal{K}_{nl} \, q^3 
\!+\!\frac{E_r  \mathcal{C}_{nl}}{2}  \left( \prescript{RL}{0}{\mathcal D}_{t}^{\alpha} q^3 \!+\! 3  q^2  \prescript{RL}{0}{\mathcal D}_{t}^{\alpha} q  \right) = -\mathcal{M}_b  \ddot{{{v_b}}},
\end{align}
in which
\begin{align}
\label{Eq: coeff unimodal}
& \mathcal{M} = \int_{0}^{1} \phi^2 \, ds +  M \, \phi^2(1) + J \, {\phi^{\prime}}^2(1),  \quad
\mathcal{J} =  J \, {\phi^{\prime}}^4(1),
 \nonumber\\
& \mathcal{K}_l = \mathcal{C}_l = \int_{0}^{1}  {\phi^{\prime\prime}}^2 \,\, ds, \quad
\mathcal{K}_{nl} = \mathcal{C}_{nl} = \int_{0}^{1}  {\phi^{\prime}}^2 \, {\phi^{\prime\prime}}^2 \,\, ds, \,\,
\mathcal{M}_b = \int_{0}^{1} \phi \, ds + M \, \phi(1) .
\end{align}

\begin{remark}
	We note that one can isolate any mode of vibration $\phi_n(s), \,\, n = 1,2, \cdots, N$ (and not necessarily the primary mode) by assuming that $\phi_n(s)$ is the only active one, and thus, end up with similar equation of motion as \eqref{Eq: weak form - discrete 2}, where the coefficients in \eqref{Eq: coeff unimodal} are obtained based on the active mode $\phi_n(s)$. Therefore, we can also make sense of \eqref{Eq: weak form - discrete 2} as a decoupled equation of motion associated with mode $\phi_n(s)$, in which the interaction with other inactive modes is absent.
\end{remark}

\section{Linearized Equation: Direct Numerical Time Integration}
\label{Sec: Direct Numerical Time Integration}
%
{We linearize our equation of motion for the cantilever beam by assuming small motions \textit{(see \ref{Sec: App. Linearization})}, and obtain the following form:}
\begin{align}
\label{Eq: weak form linear - discrete 2}
\ddot{q} + E_r \, c_l \, \prescript{RL}{0}{\mathcal D}_{t}^{\alpha} q  + k_l \, q  = { -m_b \ddot{v}_b}, \quad q(0) = \frac{v(L,0)}{\phi(L)}, \quad \dot{q}(0) = \frac{\dot{v}(L,0)}{\phi(L)}
\end{align}
with the coefficients
\begin{equation}\label{eq:linear_coefs}
c_l = \frac{\mathcal{C}_l}{\mathcal{M}}, \qquad k_l = \frac{\mathcal{K}_l}{\mathcal{M}},\qquad m_b = \frac{\mathcal{M}_b}{\mathcal{M}}.    
\end{equation}
The linearized, unimodal form \eqref{Eq: weak form linear - discrete 2} {is equivalent to the vibration of a lumped fractional Kelvin-Voigt rheological element, and can be thought of as a fractional oscillator, shown schematically in Fig.\ref{Fig:Frac_oscillator}. 
\begin{figure}[h]
	\centering
	\includegraphics[width=0.4\linewidth]{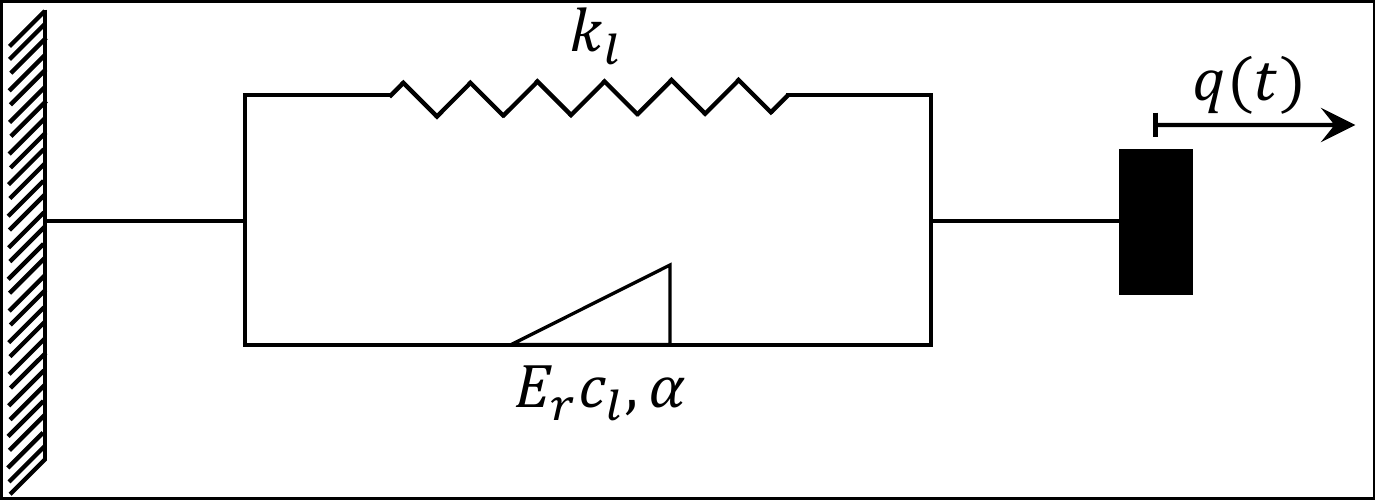}
	\caption{Lumped representation of the system as a  fractional damper, with constants $E_r \, c_l$ and fractional order $\alpha$.}
	\label{Fig:Frac_oscillator}
\end{figure}
Let a uniform time-grid with $N$ time-steps of size $\Delta t$, such that $t_n = n \Delta t$, $n = 0,\,1,\,\dots,\,N$. We employ the following equivalence relationship between the Riemann-Liouville and Caputo definitions:
\begin{equation}\label{eq:eqv_RL_C}
    \prescript{RL}{0}{\mathcal D}_{t}^{\alpha} q(t) = \prescript{C}{0}{\mathcal D}_{t}^{\alpha} q(t) + \frac{q(0)}{\Gamma(1-\alpha) t^\alpha}.
\end{equation}
Substituting (\ref{eq:eqv_RL_C}) into (\ref{Eq: weak form linear - discrete 2}), evaluating both sides implicitly at $t = t_{n+1}$, and approximating the time-fractional Caputo derivative through an L1-difference scheme \cite{lin2007finite}, we obtain:
\begin{equation}\label{eq:linear_L1}
    \ddot{q}_{n+1} + E_r c_l \left[\frac{1}{\Delta t^\alpha \Gamma(2-\alpha)} \left( q_{n+1} - q_n + \mathcal{H}^\alpha q_{n+1} \right) + \frac{q_0}{\Gamma(1-\alpha) t^\alpha_{n+1}} \right] + k_l q_{n+1} = -m_b \ddot{v}_{b,n+1},
\end{equation}
where $\mathcal{H}^{\alpha}_{k+1} = \sum^{n-1}_{j=0} = b_j (q_{n-j+1} - q_{n-j}) $ represents the discretized history term, with $\alpha$-dependent convolution coefficients $b_j = (j+1)^\alpha - j^\alpha$. We approximate the acceleration $\ddot{q}_{n+1}$ and velocity $\dot{q}_{n+1}$ through a Newmark-$\beta$ method as follows:
\begin{equation}\label{eq:Newmark}
    \ddot{q}_{n+1} = a_1 \left(q_{n+1} - q_n\right) - a_2 \dot{q}_n - a_3 \ddot{q}_n,
\end{equation}
\begin{equation}
    \dot{q}_{n+1} = a_4 \left(q_{n+1} - q_n\right) + a_5 \dot{q}_n + a_6 \ddot{q}_n,
\end{equation}
with approximation coefficients given by
\begin{equation*}
		a_1 = \frac{1}{\beta \Delta t^2}, \quad a_2 = \frac{1}{\beta \Delta t}, \quad a_3 = \frac{1-2 \beta}{2 \beta}, \quad
		a_4 = \frac{\gamma}{\beta \Delta t^2}, \quad a_5 = \left( 1- \frac{\gamma}{\beta}\right), \quad a_6 = \left( 1 - \frac{\gamma}{2 \beta} \right) \Delta t,
\end{equation*}
where we choose $\beta = 0.5$, $\gamma = 0.25$ for unconditional stability. Inserting (\ref{eq:Newmark}) into (\ref{eq:linear_L1}), we obtain the following closed form for $q_{n+1}$:
\begin{equation}\label{eq:linear_final_displacement}
       q_{n+1} = \frac{\left(a_1 + E^*\right) q_n +a_2 \dot{q}_n + a_3 \ddot{q}_n -m_b \ddot{v}_{b,n+1} - E^*\left[ \mathcal{H}^\alpha_{n+1} + q_0 \frac{(1-\alpha)}{(n+1)^\alpha}  \right]}{a_1 + E^* + k_l}
\end{equation}
with $E^* = (E_r\, c_l) / (\Delta t^\alpha \Gamma(2-\alpha))$. We observe that since the Newmark method is second-order accurate with respect to $\Delta t$, the overall accuracy is dominated by the accuracy of the $L1$ scheme, which is of $\mathcal{O}(\Delta t^{2-\alpha})$. We also observe that a discretization of a Caputo-variant of the FDE (\ref{Eq: weak form linear - discrete 2}) is recovered if we remove the term $q_0 (1-\alpha)/(n+1)^\alpha$ from (\ref{eq:linear_final_displacement}).

{We consider two numerical tests. In the first one, we solve the above system under harmonic base excitation, and in the second one, we consider a free-vibration response. For both tests, we set $E_r = 1$ and consider the lumped mass at the tip, with $M = J = 1$, that is, we utilize (\ref{eq:eig_mass}) for $\phi(s)$, which yields the coefficients $c_l = k_l = 1.24$. }

{
\subsection{Harmonic Base Excitation}
We solve  (\ref{Eq: weak form linear - discrete 2}) in the presence of base excitation in the harmonic form  $v_b = a_b sin(\omega_b t)$, where $\omega_b \in [0.5,3.5]$ and $a_b=0.01$ denote, respectively, the base frequency and displacement amplitude. The coefficient $m_b = -0.042$ is calculated through (\ref{eq:linear_coefs}) and (\ref{Eq: coeff unimodal}). We employ homogeneous initial conditions, \textit{i.e.}, $q(0) = 0$, $\dot{q}(0) = 0$, and set the time $t \in (0,100]$, with step size $\Delta t = 10^{-3}$. The maximum displacement amplitude after reaching the steady state response of the system is evaluated. Figure \ref{Fig: maxamp} illustrates the amplitude \textit{vs} base frequency response with respect to varying fractional orders $\alpha$. We observe the existence of a critical point at $\omega_b = 1$ that changes the dissipation nature of the fractional order parameter. Regarding the maximum observed amplitudes, increasing the fractional order in the range $\alpha \in [0.1,0.4]$, decreases and slightly shifts the amplitude peaks to higher (right) frequencies (an anomalous quality). On the other hand, as the fractional order is increased in the range $\alpha \in[0.5,0.6]$, the peak amplitudes slightly shift towards the lower (left) frequencies, which is also observed in standard systems with the increase of modal damping values.

\begin{figure}[h]
	\centering
	\includegraphics[width=0.75\linewidth]{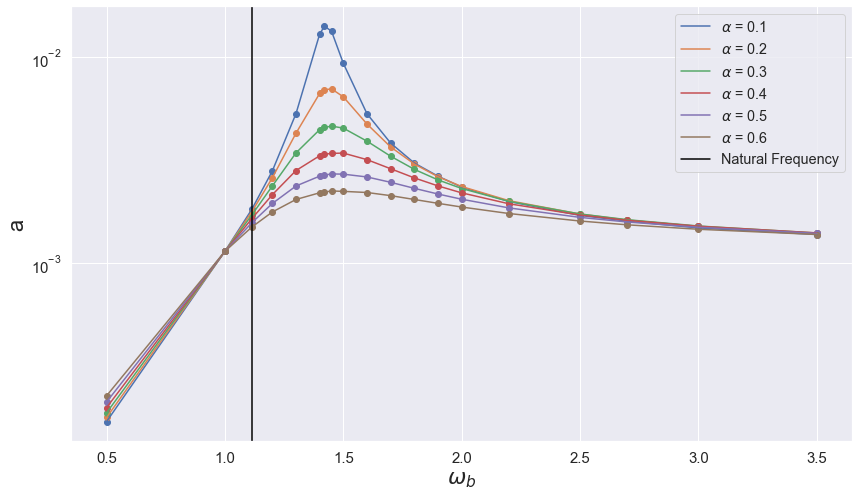}
	\caption{Anomalous change of the maximum amplitude vs frequency at the tip of the beam in presence of harmonic base excitation for different fractional order $\alpha$. The solid black line represents the standard, undamped natural frequency of the system.}
	\label{Fig: maxamp}
\end{figure}
}

%

{
\subsection{Free Vibration}

Following the observed anomalous amplitude \textit{vs.} base frequency behaviors and presence of a critical point nearby the standard natural frequency of the system illustrated in Figure \ref{Fig: maxamp}, we solve (\ref{Eq: weak form linear - discrete 2}) in a free-vibration setting employing Riemann-Liouville and Caputo definitions, where we set $\ddot{v}_b = 0$, and $q(0) = 0.01$, $\dot{q}(0) = 0$.} Figure \ref{Fig: fracOscil} (left) illustrates {the obtained results for $q(t)$ for varying fractional orders using a Riemann-Liouville definition. We observe an $\alpha$-dependent amplitude decay, which converges to a classical integer-oder oscillator as $\alpha \to 1$. Furthermore, an anomalous transient region is observed at the short time-scale $t \in [0.1,\,1]$. On the other hand, in Figure \ref{Fig: fracOscil} (right), anomalies are present at large time-scales through a (far-from-equilibrium) power-law relaxation, while the short-time behavior is ``standard-like". Such contrast between the obtained results provides interesting insights towards modeling desired anomalous ranges in such power-law materials. By replacing the fractional 
damper with a classical integer-order one \textit{(see Fig.\ref{Fig: freevib_classical})}, we notice that neither anomalous dynamics are present.} The obtained results are in agreement with the power-law and exponential relaxation kernels described in Sec.\ref{Sec: Viscoelasticity}. We note that since the fractional element provides a constitutive interpolation between spring and dash-pot elements (see Sec.\ref{Sec: Viscoelasticity} for more discussion and references), it contributes to both effective stiffness and damping ratio of the system, and therefore increasing values of $\alpha$ (decreasing stiffness), yield a reduction in the frequency response. 

Fractional linear oscillators are also considered in 
\cite{svenkeson2016spectral} {for systems with memory}, where their 
interaction with a fluctuating environment causes the time evolution of the 
system to be intermittent. The authors in \cite{svenkeson2016spectral} apply 
the Koopman operator theory to the corresponding integer order system and then 
make a L$\grave{\text{e}}$vy transformation in time to recover long-term memory 
effects; they observe a power-law behavior in the amplitude decay of the 
system's response. Such an anomalous decay rate has also been investigated in 
\cite{shoshani2017anomalous} for an extended theory of decay of classical 
vibrational models brought into nonlinear resonances. The authors report a 
``non-exponential" decay in variables describing the dynamics of the system in 
the presence of dissipation and also a sharp change in the decay rate close to 
resonance. 

%
\begin{figure}[t]
	\centering
	\begin{subfigure}{0.49\textwidth}
		\centering
		\includegraphics[width=1\linewidth]{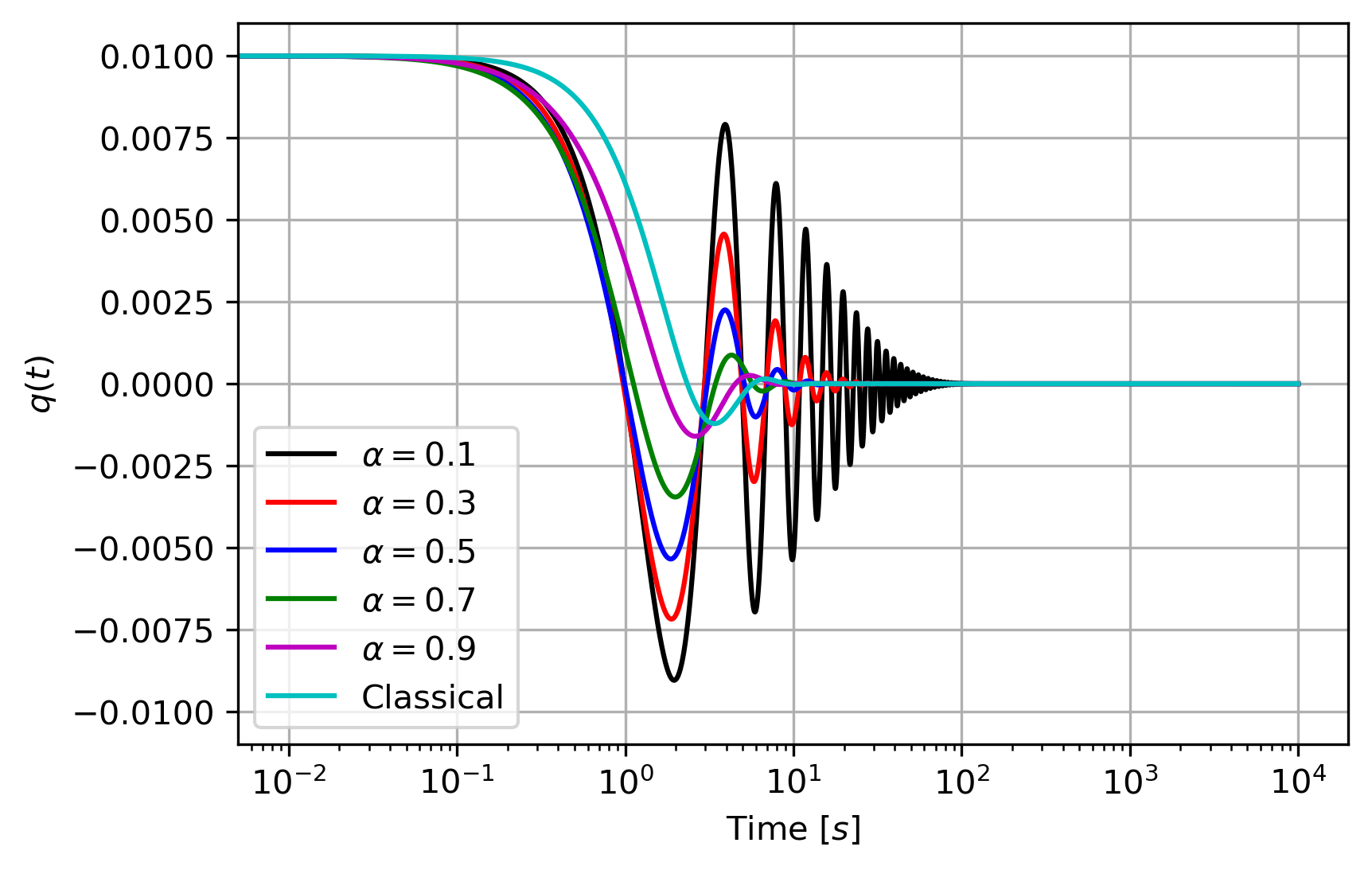}
	\end{subfigure}
	\begin{subfigure}{0.49\textwidth}
		\centering
		\includegraphics[width=1\linewidth]{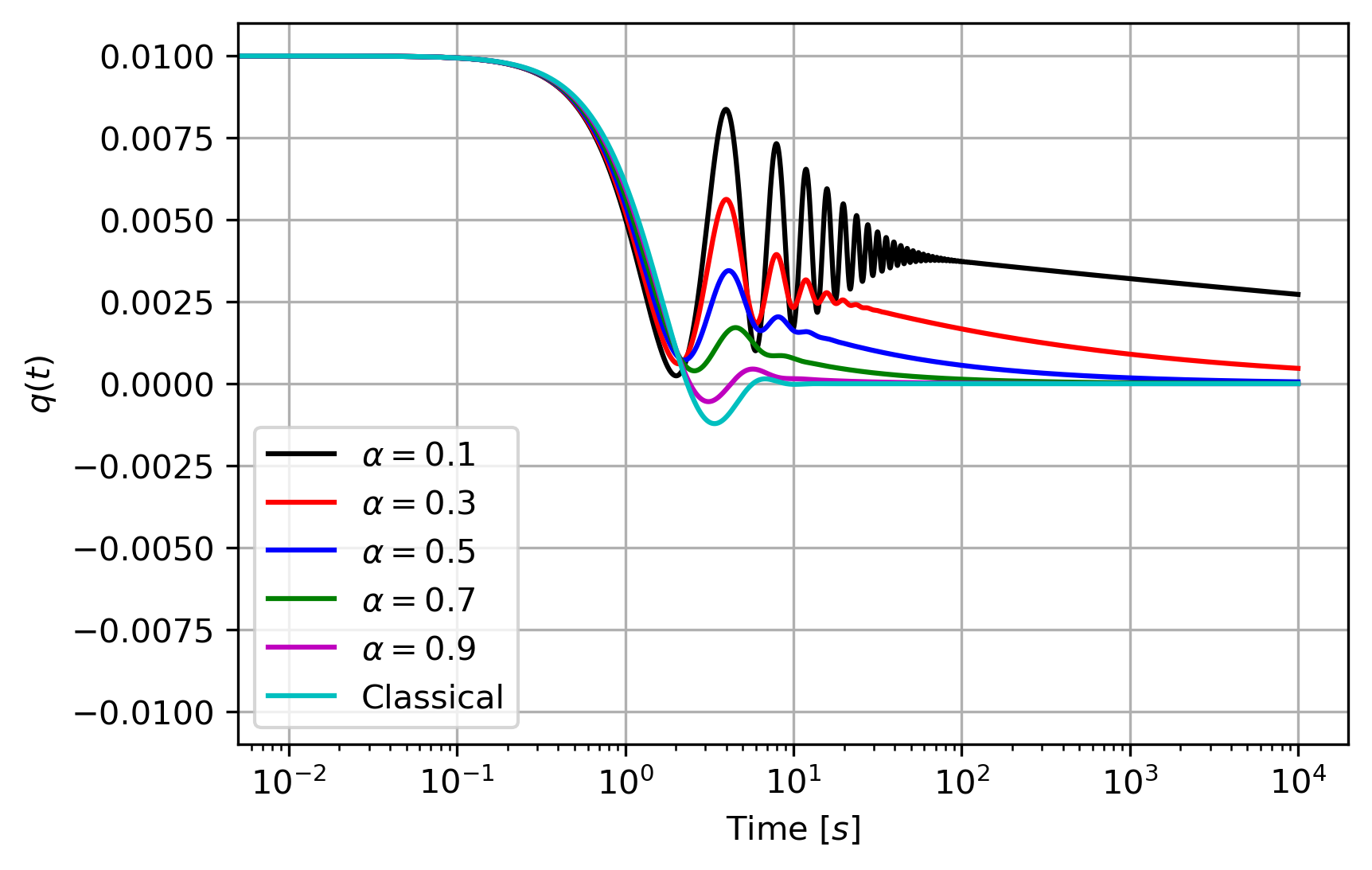}
	\end{subfigure}
	\caption{Anomalous linear free vibration modal displacement \textit{vs.} time. \textit{(Left)} Riemman-Liouville definition with short-time anomalies. \textit{(Right)} Caputo definition with long-time anomalies. }
	\label{Fig: fracOscil}
\end{figure}
%

\begin{figure}[h]
	\centering
	\includegraphics[width=0.49\linewidth]{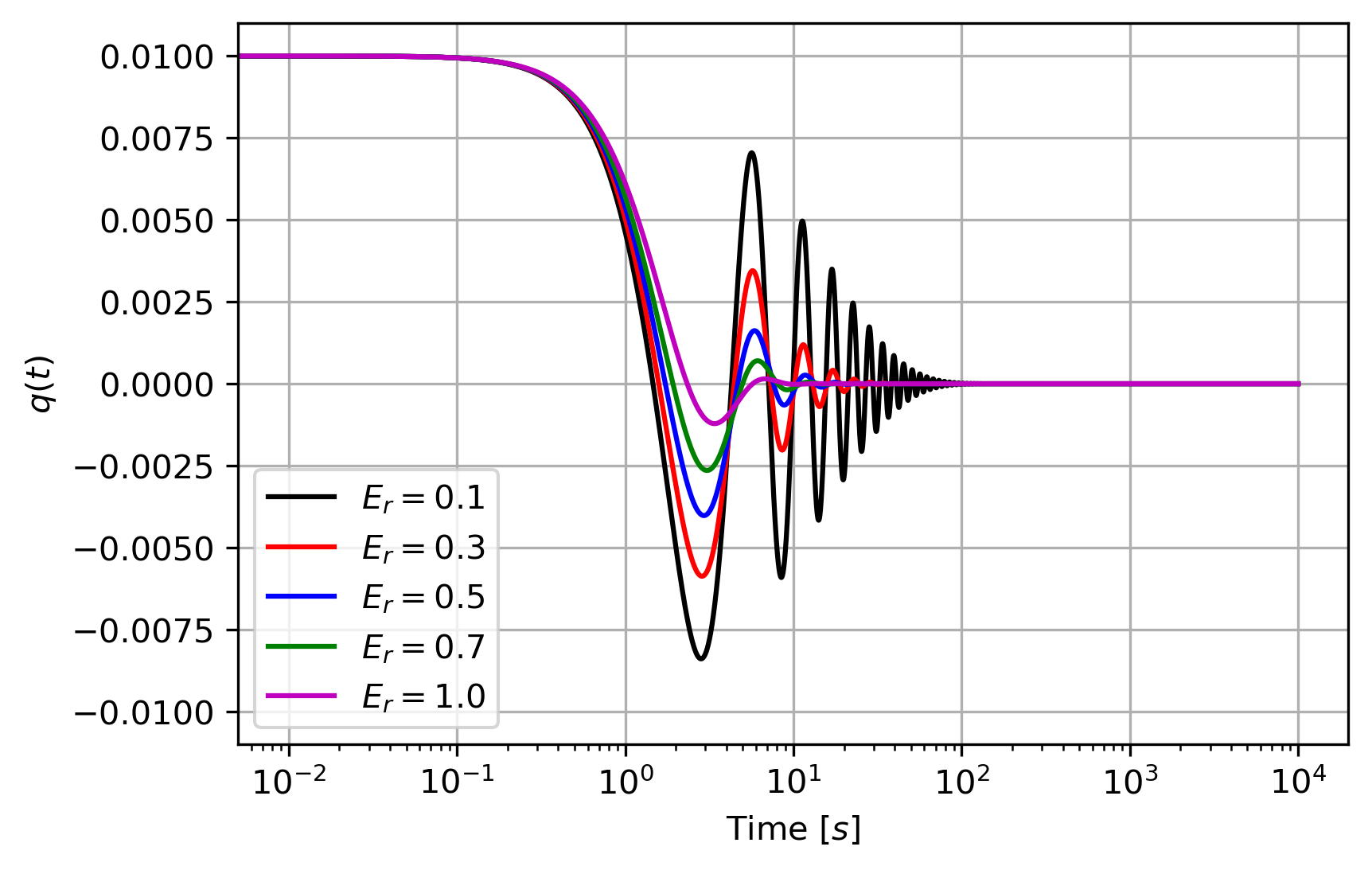}
	\caption{Classical linear free vibration modal displacement \textit{vs.} time under varying values of $E_r$.}
	\label{Fig: freevib_classical}
\end{figure}
}

%
\section{Perturbation Analysis of Nonlinear Equation}
\label{Sec: Perturbatin Analysis}
%
We use perturbation analysis to investigate the behavior of a nonlinear system, where we reduce a nonlinear fractional differential equation to an algebraic equation to solve for the steady state amplitude and phase of vibration.
%
\subsection{Method of Multiple Scales}
%
To investigate the dynamics of the system described by \eqref{Eq: weak form - discrete 2}, we use the method of multiple scales \cite{nayfeh2008nonlinear,rossikhin2010application}. The new independent time scales and the integer-order derivative with respect to them are defined as
\begin{equation}
\label{Eq: time scales}
T_m = \epsilon^{m}\, t, \qquad D_m = \frac{\partial}{\partial T_m} , \qquad m = 0, 1, 2, \cdots.
\end{equation}
It is also convenient to utilize another representation of the fractional derivative {\textit{(see \cite{samko1993fractional}, Equation 5.82)},} which according to the Rieman-Liouville fractional derivative, is equivalent to the fractional power of the operator of conventional time-derivative, i.e. $\prescript{RL}{0}{\mathcal D}_{t}^{\alpha} = (\frac{d}{dt})^{\alpha}$. Therefore,  
\begin{align}
\label{Eq: time derivatives}
\frac{d}{dt} &= D_0 + \epsilon D_1 + \cdots , \qquad 
\frac{d^2}{dt^2} &= D^2_0 + 2 \epsilon D_0 D_1 + \cdots, \qquad 
\prescript{RL}{0}{\mathcal D}_{t}^{\alpha} &= (\frac{d}{dt})^{\alpha} = D^{\alpha}_0 + \epsilon \alpha D^{\alpha-1}_0 D_1 + \cdots
\end{align}
The solution $q(t)$ can then be represented in terms of series {expansion:}
\begin{align}
\label{Eq: MMS expansion}
q(T_0,T_1,\cdots) &= q_0(T_0,T_1,\cdots) + \epsilon q_1(T_0,T_1,\cdots) + \epsilon^2 q_2(T_0,T_1,\cdots) + \cdots
\end{align}
We assume that the coefficients in the equation of motion {have} the following scaling
\begin{align}
\label{Eq: coeff scaling}
&\frac{\mathcal{J}}{\mathcal{M}} = \epsilon \, m_{nl} , \quad
\frac{\mathcal{K}_l}{\mathcal{M}} = k_l = \omega^2_0  , \quad
\frac{\mathcal{C}_l}{\mathcal{M}} = \epsilon \, c_l , \quad 
\frac{\mathcal{K}_{nl}}{\mathcal{M}} = \epsilon \, k_{nl} , \quad
\frac{\mathcal{C}_{nl}}{\mathcal{M}} = \epsilon \, c_{nl} , 
\end{align}
and the base excitation $- \frac{\mathcal{M}_b}{\mathcal{M}} \, \ddot{{{v_b}}} $ is a harmonic function {in the} form $\epsilon \, F \cos(\Omega \, t)$. Thus, \eqref{Eq: weak form - discrete 2} can be expanded as
\begin{align}
\label{Eq: governing eqn unimodal - 4}
&(D^2_0 + 2 \epsilon D_0 D_1 + \cdots)(q_0 + \epsilon q_1 + \cdots) 
+ \epsilon \, m_{nl} (D^2_0 + 2 \epsilon D_0 D_1 + \cdots)(q_0 + \epsilon q_1 + \cdots) \nonumber \,\times\, (q_0 + \epsilon q_1 + \cdots)^2 
 \nonumber\\
+ & \epsilon \, m_{nl} (q_0 + \epsilon q_1 + \cdots)
\times\, \left( (D_0 + \epsilon D_1 + \cdots)(q_0 + \epsilon q_1 + \cdots) \right)^2 
+ \omega^2_0 \, (q_0 + \epsilon q_1 + \cdots) 
\nonumber\\
+ & \epsilon \, E_r \, c_l \, (D^{\alpha}_0 + \epsilon \alpha D^{\alpha - 1}_0 D_1 + \cdots)(q_0 + \epsilon q_1 + \cdots) 
+ 2 \epsilon \, k_{nl} \, (q_0 + \epsilon q_1 + \cdots)^3 
 \nonumber\\
+& \frac{1}{2} \epsilon \, E_r \, c_{nl} \, (D^{\alpha}_0 + \epsilon \alpha D^{\alpha - 1}_0 D_1 + \cdots)(q_0 + \epsilon q_1 + \cdots)^3 
+ \frac{3}{2} \epsilon \, E_r \, c_{nl}(q_0 + \epsilon q_1 + \cdots)^2 \nonumber\\
&   \left[ (D^{\alpha}_0 + \epsilon \alpha D^{\alpha - 1}_0 D_1 + \cdots)(q_0 + \epsilon q_1 + \cdots) \right]
= \epsilon \, F \cos(\Omega \, T_0).
\end{align}
By collecting similar coefficients of zero-th and first orders of $\epsilon$, we obtain the following equations
\begin{alignat}{2}
\mathcal{O}(\epsilon^0) :  D_0^2 q_{0} + \omega_0^2 q_{0}\,\,=&\,\, 0 , \label{Eq: order zero} \\
\mathcal{O}(\epsilon^{1}) :  D_0^2 q_{1} + \omega_0^2 q_{1} \,\,=&\,\, -2 D_0 D_1 q_{0}
- m_{nl} \, \left( q_0^2 D_0^2 q_0 + q_0 (D_0 q_0)^2 \right) - E_r \, c_l \, D^{\alpha}_0 q_{0} - 2 \, k_{nl} \, q^3_{0} 
\nonumber \\
&- \frac{1}{2} \, E_r \, c_{nl} \, D^{\alpha}_0 q^3_{0}- \frac{3}{2} \, E_r \, c_{nl} \, q^2_0 D^{\alpha}_0 q_{0}
+F \cos(\Omega \, T_0).         \label{Eq: first order}
\end{alignat}
The solution to \eqref{Eq: order zero} is of the form
\begin{align}
\label{Eq: order zero solution}
q_0(T_0,T_1) = A(T_1) \, e^{i \, \omega_0 \, T_0} + c.c 
%
\end{align}
where ``c.c" denotes the complex conjugate. By substituting \eqref{Eq: order zero solution} into the right-hand-side of \eqref{Eq: first order}, we observe that different resonance cases are possible. In each case, we obtain the corresponding solvability conditions by removing the secular terms, i.e. the terms that grow in time unbounded. Then, {we utilize the} polar form $A = \frac{1}{2} a \, e^{i \, \varphi}$, where the real valued functions $a$ and $\varphi$ are the amplitude and phase lag of time response, respectively. Thus, the solution $q(t)$ becomes
\begin{align}
\label{Eq: order zero solution - 2}
q(t) = a(\epsilon \, t) \, \text{cos}(\omega_0 \, t + \varphi(\epsilon \, t)) + \mathcal{O}(\epsilon),
\end{align}
where the governing equations of $a$ and $\varphi$ are obtained by separating the real and imaginary parts. 

%
\subsubsection{Case 1: No Lumped Mass At The Tip}
\label{Sec: case 1}
\begin{figure}[h]
	\centering
	\begin{subfigure}{0.4\textwidth}
		\centering
		\includegraphics[width=0.9\linewidth]{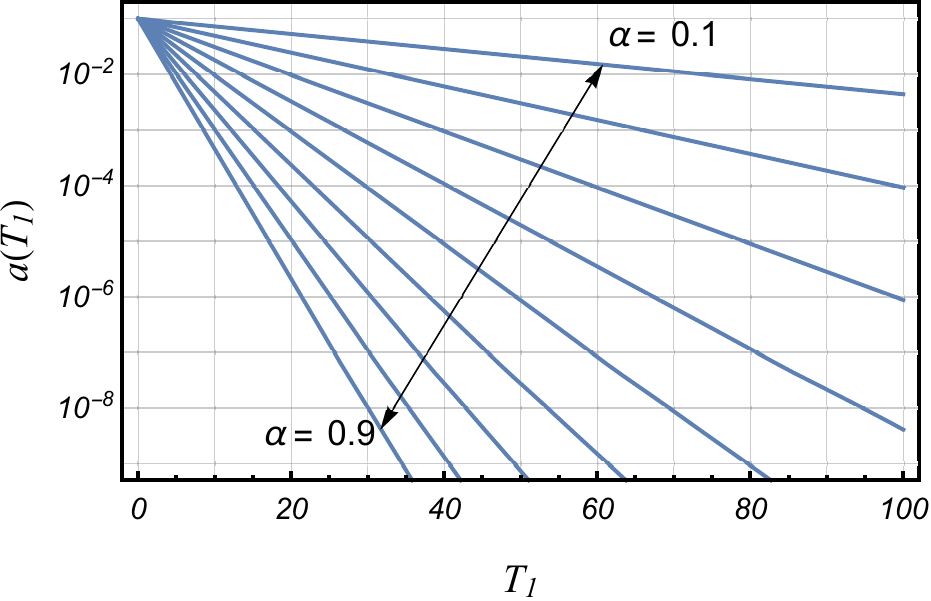}
	\end{subfigure}
	\begin{subfigure}{0.4\textwidth}
		\centering
		\includegraphics[width=0.9\linewidth]{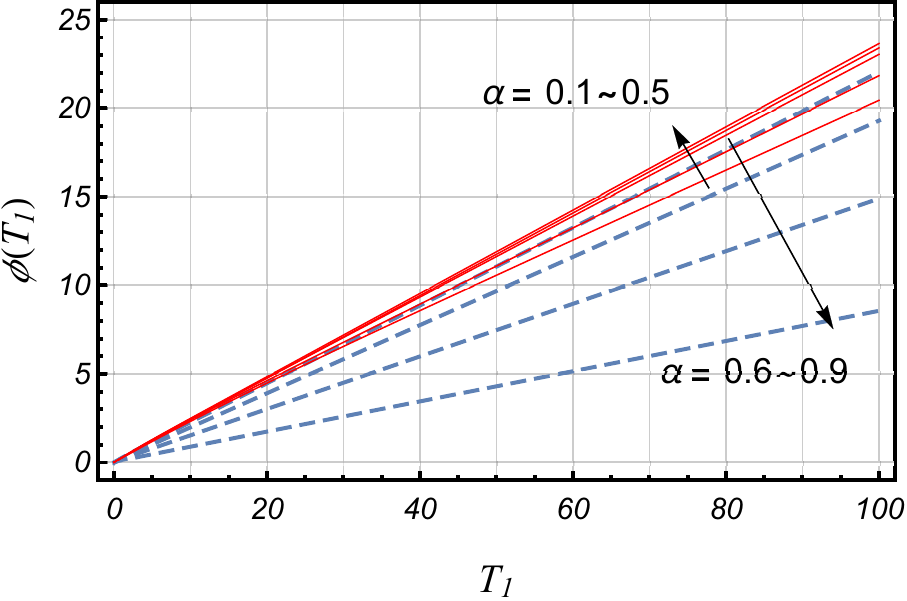}
	\end{subfigure}
	\caption{{Nonlinear anomalous free vibration of a viscoelastic cantilever beam with no lumped mass at the tip and $E_r = 0.1$.} The rate of amplitude {decay strongly depends on the fractional order $\alpha$, where a rapid decay is observed as $\alpha$ is increased \textit{(left)}. On the other hand, for increasing $\alpha$, the phase lag $\varphi(\epsilon t)$, increases in the lower range of $\alpha$, and decreases in the higher range of $\alpha$}}
	\label{Fig: perturbation free vib}
\end{figure}

%
In this case, $M = J = 0$, and thus, given the {form (\ref{eq:eig_no_mass}) for the eigenfunctions $\phi(s)$ in \ref{Sec: App. Eigenvalue Problem of Linear Model}}, the coefficients are computed as $\mathcal{M} = 1$, $\mathcal{K}_l = \mathcal{C}_l =12.3624$, $\mathcal{M}_b = 0.782992$, and $\mathcal{K}_{nl} = \mathcal{C}_{nl} = 20.2203$. We consider the following cases:

\vspace{0.2 in}
\noindent$\bullet$ Free Vibration, $F = 0$: Super Sensitivity to $\alpha$\\
In this case, the beam is not externally excited and thus, $F = 0$. By removing the secular terms that are the coefficients of $e^{i \, \omega_0 \, T_0}$ in the solvability condition, we find the governing equations of solution amplitude and phase as 
%
\begin{align}
\label{Eq: free amp}
\frac{d a}{d T_1}&= 
- E_r \, \omega_0^{\alpha-1} \sin (\alpha \frac{\pi }{2}) 
\left( \frac{1}{2} \, c_l \, a +  \frac{3}{8} \, c_{nl} \, a^3 \right),
\\ \label{Eq: free phase}
\frac{d \varphi}{d T_1} &=
\frac{1}{2} c_l \, E_r \, \omega _0^{\alpha-1} \, \cos \left(\frac{\pi  \alpha }{2}\right)  
+ \frac{3}{4} c_{nl} \, E_r \, \omega _0^{\alpha-1} \, \cos \left(\frac{\pi  \alpha }{2}\right) \, a^2 
+\frac{3}{4} \, \omega _0^{-1} \, k_{nl} \, a^2 .
\end{align}
We can see from the first equation \eqref{Eq: free amp} that the amplitude of free vibration decays out, where the decay rate $\tau_d = c_l \, E_r \, \omega_0^{\alpha-1} \sin (\alpha \frac{\pi }{2})$ directly depends on values of the fractional derivative $\alpha$ and the coefficients $E_r$ (see Fig. \ref{Fig: perturbation free vib}). 
%
%
We introduce the sensitivity index $S_{\tau_d, \alpha}$ as the partial derivative of decay rate with respect to $\alpha$, i.e. 
%
\begin{figure}[h]
	\centering
	\includegraphics[width=0.4\linewidth]{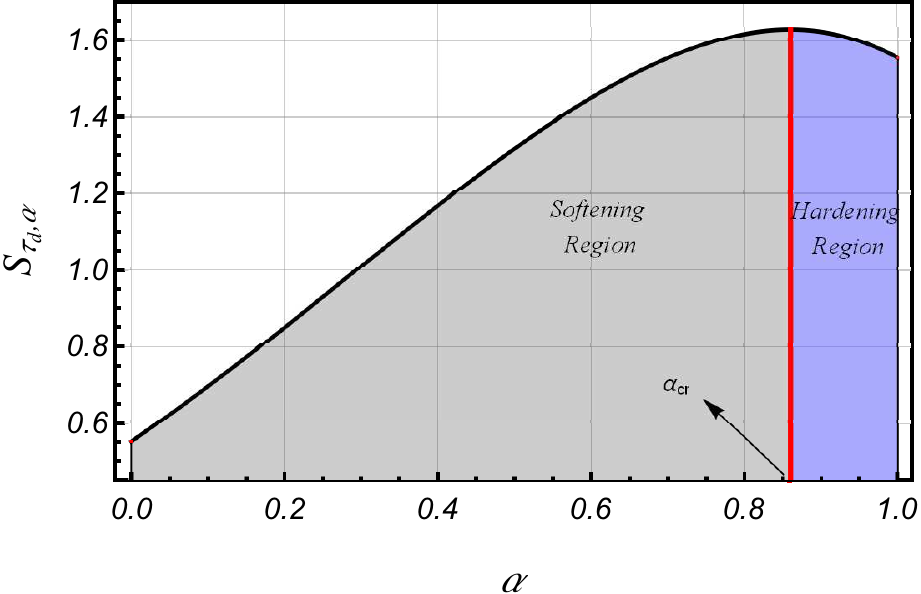}
	\caption{{Anomalous super-sensitivity of the decay rate $\tau_d$ with respect to $\alpha$ under free vibration}. Increasing $\alpha$ when $\alpha < \alpha_{cr}$ leads to higher dissipation and decay rate. The reverse effect is observed when $\alpha > \alpha_{cr}$. {Here, the notions of softening/hardening are associated to high/lower decay rates as $\alpha$ is increased (introducing extra viscosity)}.}
	\label{Fig: perturbation free vib Sen}
\end{figure}
%
\begin{align}
\label{Eq: Free Vib Amp Sen}
S_{\tau_d, \alpha} = \frac{d \tau_d}{d \alpha} = &\frac{\pi}{2} c_l \, E_r \, \omega_0^{\alpha -1} \, \cos(\alpha \frac{\pi }{2}) + c_l \, E_r \, \omega _0^{\alpha -1} \sin(\alpha \frac{\pi }{2})  \log(\omega _0).
\end{align}

The sensitivity index is computed and plotted in Fig. \ref{Fig: perturbation free vib Sen} for the same set of parameters as in Fig. \ref{Fig: perturbation free vib}. There exists a critical value
\begin{align}
\label{Eq: Free Vib Amp Sen critical}
%
\alpha_{cr}  = 
-\frac{2}{\pi}\tan^{-1} \left(\frac{\pi}{2 \, \log(\omega _0)} \right)  ,
\end{align}
where $(dS_{\tau_d, \alpha} / d\alpha) = 0 $. We observe in Fig.\ref{Fig: 
	perturbation free vib Sen} that by increasing $\alpha$ when $\alpha < 
\alpha_{cr}$, i.e. introducing more viscosity to the system, the dissipation rate, and thus decay rate, increases; this can be {interpreted} as a softening (stiffness-decreasing) region. Further increasing $\alpha$ when $\alpha > \alpha_{cr}$, will reversely results in decrease of decay rate; this can be {interpreted} as a hardening (more stiffening) region. We also note that $\alpha_{cr}$ solely depends on value of $\omega_0$, given in \eqref{Eq: coeff 	scaling}, and even though the value of $E_r$ affects decay rate, it does not change the value of $\alpha_{cr}$. {Therefore, the region of super-sensitivity, where the anomalous transition between softening/hardening regimes takes place only depends on the \textit{standard} natural frequency of the system.}

Although the observed hardening response after a critical value of $\alpha$ in Fig.\ref{Fig: perturbation free vib Sen} might seem counter-intuitive at first, we remark that here the notions of softening and hardening have a mixed nature regarding energy dissipation and time-scale dependent material stress response, which have anomalous nature for fractional viscoelasticity. {Similar anomalous dynamics were also observed in ballistic, strain-driven yield stress responses of fractional visco-elasto-plastic truss structures \cite{suzuki2016fractional}. In the following,} we demonstrate two numerical tests by purely utilizing the constitutive response of the fractional Kelvin-Voigt model (\ref{Eq: frac KV}) to justify the observed behavior in Fig.\ref{Fig: perturbation free vib Sen} by employing the \textit{tangent loss} and the stress-strain response under monotone loads/relaxation.

\begin{figure}[h]
	\centering
	\begin{subfigure}{0.33\textwidth}
		\centering
		\includegraphics[width=1\linewidth]{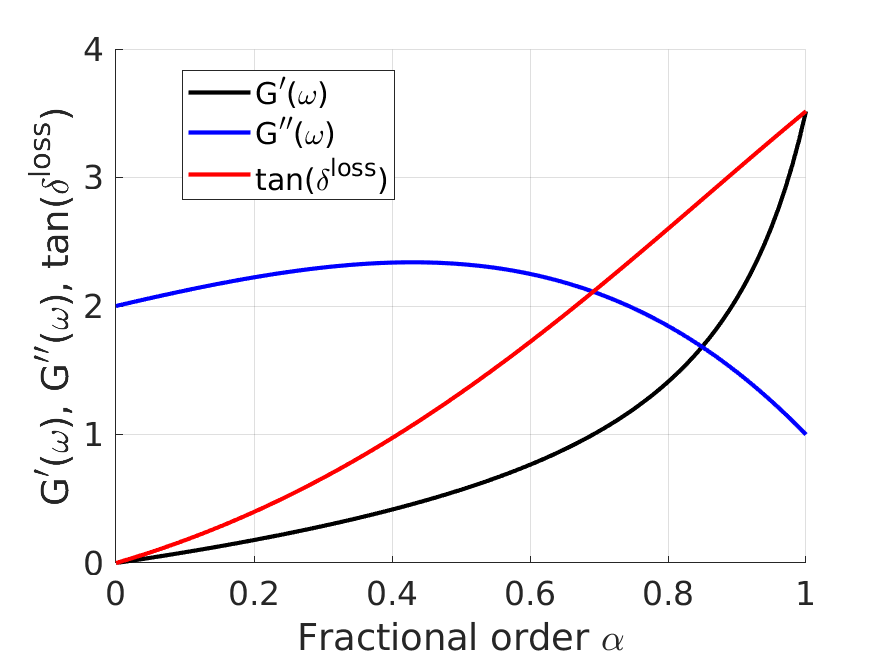}
	\end{subfigure}
	\begin{subfigure}{0.33\textwidth}
		\centering
		\includegraphics[width=1\linewidth]{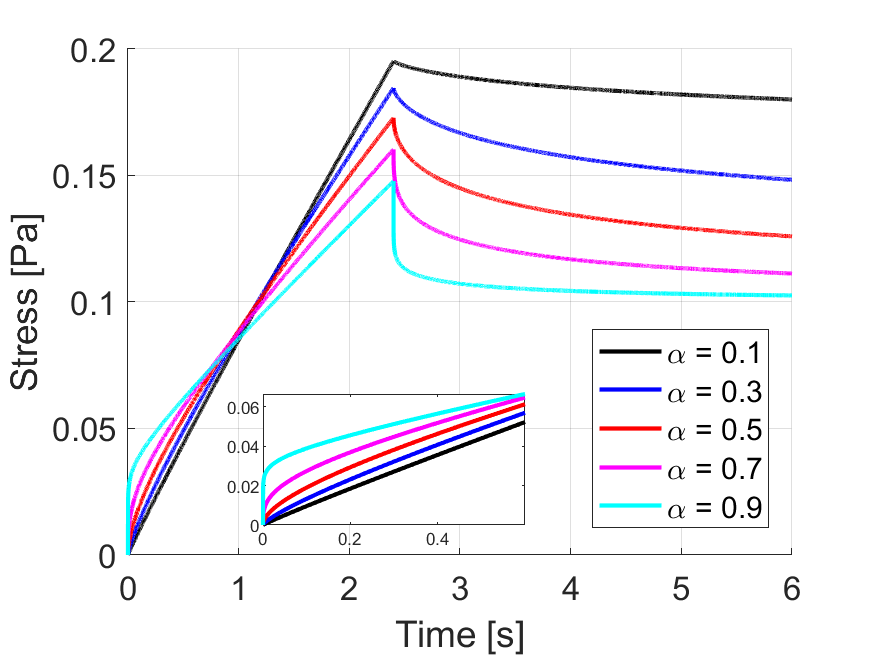}
	\end{subfigure}
	\caption{\textit{(Left)} Storage and loss moduli, and tangent loss for the fractional Kelvin-Voigt model at $\omega_0$ with varying fractional-orders and $E_r = 1$. \textit{(Right)} Stress-time response under a monotone load with constant strain rate undergoing ballistic hardening response for short-time and higher $\alpha$, followed by a stress relaxation.}
	\label{Fig:Constitutive_response}
\end{figure}

\noindent\textbf{Dissipation \textit{via} tangent loss}: By taking the Fourier transform of (\ref{Eq: frac KV}), we obtain the so-called complex modulus $G^*$ \cite{mainardi2010fractional}, which is given by:
\begin{equation}
G^*(\omega) = E_\infty + E_\alpha \omega^\alpha \left( \cos\left(\alpha \frac{\pi}{2}\right) + i \sin\left(\alpha \frac{\pi}{2}\right) \right),
\end{equation}
from which the real and imaginary parts yield, respectively, the storage and loss moduli, as follows:
\begin{align*}
G^\prime(\omega) = E_\infty + E_\alpha \omega^\alpha \cos\left(\alpha \frac{\pi}{2}\right), \qquad
G^{\prime\prime}(\omega) = E_\alpha \omega^\alpha \sin\left(\alpha \frac{\pi}{2}\right),
\end{align*}
which represent, respectively, the {stored and dissipated energies per cycle}. Finally, we define the tangent loss, which represents {the} ratio between {the dissipated/stored energies}, and therefore related to the mechanical damping of the {anomalous medium:}
\begin{equation}\label{Eq:Tangent_loss}
\tan \delta^{loss} = \frac{G^{\prime\prime}(\omega)}{G^\prime(\omega)} = \frac{E_r \omega^\alpha \sin\left(\alpha\frac{\pi}{2}\right)}{1 + E_r\omega^\alpha\cos\left(\alpha\frac{\pi}{2}\right)}
\end{equation}
We set $\omega = \omega_0$ and $E_r = 1$ and demonstrate the results for (\ref{Eq:Tangent_loss}) with varying fractional orders $\alpha$. We present the obtained results in Fig.\ref{Fig:Constitutive_response} \textit{(left)}, where we observe that increasing fractional orders lead to increased dissipation per loading cycle with the increase of the tangent loss, and the hardening part ($\alpha > \alpha_{c}$) is not associated with higher storage in the material. Instead, the increasing dissipation with $\alpha$ suggests an increasing damping of the mechanical structure.

\noindent\textbf{Stress-time response for monotone loads/relaxation}: In this test, we demonstrate how increasing fractional orders for the fractional model leads to increased hardening for sufficiently high strain rates. Therefore, we directly {discretize (\ref{Eq: frac KV}) utilizing an L1-scheme \cite{lin2007finite} in a uniform time-grid and set $E_\infty = 1$, $E_\alpha = 1$. We also assume the following piecewise strain function:} $\varepsilon(t) = (1/24) t$, for $0 \le t < 2.5$ (monotone stress/strain), and $\varepsilon(t) = 1/10$ for $2.5 \le t \le 6$ (stress relaxation). The obtained results are illustrated in Fig.\ref{Fig:Constitutive_response} \textit{(right)}, where we observe that even for relatively low strain rates, there is a ballistic region nearby the initial time where higher fractional orders present higher values of stress, characterizing a {rate-dependent stress-hardening} response. However, due to the dissipative nature of fractional rheological elements, the initially higher-stress material softens {after passing a critical point,} due to its faster relaxation nature.
%

%

\begin{figure}[h]
	\centering
	\begin{subfigure}{0.35\textwidth}
		\centering
		\includegraphics[width=1\linewidth]{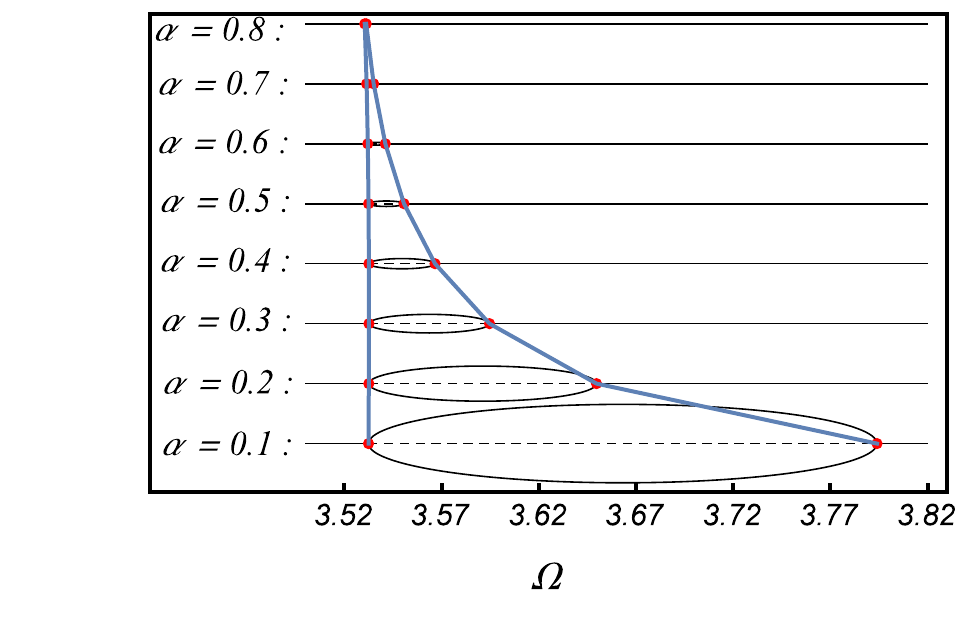}
	\end{subfigure}
	\begin{subfigure}{0.35\textwidth}
		\centering
		\includegraphics[width=1\linewidth]{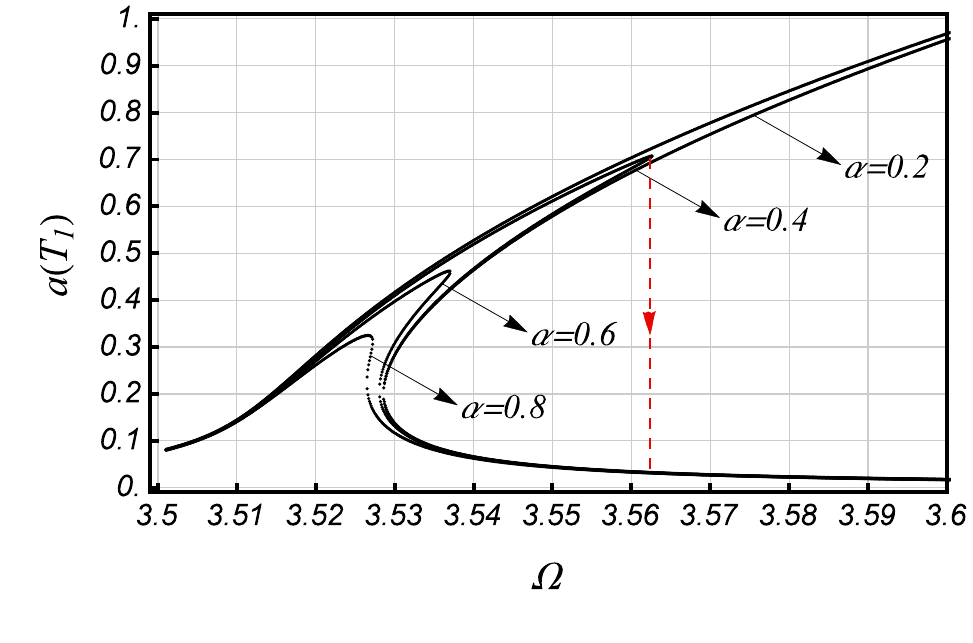}
	\end{subfigure}
	\caption{Primary resonance of the viscoelastic cantilever beam with no lumped mass at the tip. Steady state amplitude \textit{(right)} and its bifurcation diagram \textit{(left)} by changing the detuning parameter $\Delta$ for different values of $\alpha$ and $E_r=0.3, f=1$.   }
	\label{Fig: biff diag}
\end{figure}

\vspace{0.2 in}
\noindent$\bullet$ Primary Resonance Case, $\Omega \approx \omega_0$\\
In the case of primary resonance, the excitation frequency is close to the natural frequency of the system. We let $\Omega = \omega_0 + \epsilon \, \Delta$, where $\Delta$ is called the detuning parameter and thus, write the force function as $\frac{1}{2} F \, e^{i \, \Delta \, T_1} \, e^{i \, \omega_0 \, T_0} + c.c$ . In this case, the force function also contributes to the secular terms. Therefore, we find the governing equations of solution amplitude and phase as
\begin{align}
\label{Eq: primary amp}
\frac{d a}{d T_1}
= &
- E_r \, \omega_0^{\alpha-1} \sin (\alpha \frac{\pi }{2}) 
\left( \frac{1}{2} \, c_l  \, a + \frac{3}{8} \, c_{nl}  \, a^3 \right) + \frac{1}{2} f \, \omega^{-1}_0 \, \sin (\Delta  T_1 - \varphi),
\\ \label{Eq: primary phase}
a \, \frac{d \varphi}{d T_1} =&
\frac{1}{2} c_l \, E_r \, \omega _0^{\alpha-1} \, \cos(\frac{\pi  \alpha }{2})  \, a + \frac{3}{4} c_{nl} \, E_r \, \omega _0^{\alpha-1} \, \cos(\frac{\pi  \alpha }{2}) \, a^3 
+\frac{3}{4} \, \omega _0^{-1} \, k_{nl} \, a^3 
- \frac{1}{2} f \, \omega^{-1}_0 \, \cos (\Delta  T_1 - \varphi) ,
\end{align}
in which the four parameters $\{\alpha,E_r,f,\Delta\}$ mainly change the frequency response of the system. The equations \eqref{Eq: primary amp} and \eqref{Eq: primary phase} can be transformed into an autonomous system, where the $T_1$ does not appear explicitly, by letting $$ \gamma = \Delta \, T_1 - \varphi .$$ The steady state solution occur when $\frac{d a}{d T_1} = \frac{d \varphi}{d T_1} =0$, that gives
\begin{align}
\label{Eq: primary amp ss}
&E_r \, \omega_0^{\alpha-1} \sin (\frac{\pi  \alpha }{2}) 
\left( \frac{c_l}{2} a + \frac{3 c_{nl}}{8}  a^3 \right) 
= \frac{f}{2 \, \omega_0} \sin (\gamma),
\\ \label{Eq: primary phase ss}
& \left(\Delta - \frac{c_l}{2} E_r \, \omega _0^{\alpha-1} \, \cos(\frac{\pi  \alpha }{2}) \right)  a - \frac{3}{4} \left( c_{nl} \, E_r \, \omega _0^{\alpha-1} \, \cos(\frac{\pi  \alpha }{2}) + \omega _0^{-1} \, k_{nl} \right) a^3 
= \frac{f}{2 \, \omega_0} \cos (\gamma) ,
\end{align}
and thus, by squaring and adding these two equations, we get
\begin{align}
\label{Eq: primary ss}
&\left[
\frac{c_l}{2} E_r \, \omega_0^{\alpha-1} \sin (\frac{\pi  \alpha }{2})  \, a
+\frac{3 c_{nl}}{8}  E_r \, \omega_0^{\alpha-1} \sin (\frac{\pi  \alpha }{2}) \, a^3 
\right]^2
+ \left[
\left(\Delta - \frac{c_l}{2} E_r \, \omega _0^{\alpha-1} \, \cos(\frac{\pi  \alpha }{2}) \right)  a \right. \nonumber \\
&\left.- \frac{3}{4} \left( c_{nl} \, E_r \, \omega _0^{\alpha-1} \, \cos(\frac{\pi  \alpha }{2}) + \omega _0^{-1} \, k_{nl} \right) a^3 
\right]^2
= \frac{f^2}{4 \, \omega^2_0}.
\end{align}
This can be written in a simpler way as
\begin{align}
\label{Eq: primary ss - 2}
\left[ A_1 \, a + A_2 \, a^3 \right]^2
+ \left[B_1 \, a + B_2 \, a^3 \right]^2
= C,
\end{align}
where
\begin{alignat*}{4}
&A_1 = \frac{c_l}{2} E_r \, \omega_0^{\alpha-1} \sin (\frac{\pi  \alpha }{2}),  \quad
A_2 = \frac{3 c_{nl}}{8}  E_r \, \omega_0^{\alpha-1} \sin (\frac{\pi  \alpha }{2}), \\
&B_1 = \Delta - \frac{c_l}{2} E_r \, \omega _0^{\alpha-1} \, \cos(\frac{\pi  \alpha }{2}), \quad
B_2 = - \frac{3}{4} \left( c_{nl} \, E_r \, \omega _0^{\alpha-1} \, \cos(\frac{\pi  \alpha }{2}) + \omega _0^{-1} \, k_{nl} \right), \quad C =\frac{f^2}{4 \, \omega^2_0}.
\end{alignat*}
Hence, the steady state response amplitude is the admissible root of 
\begin{align}
\label{Eq: SS amp}
(A_2^2 + B_2^2) a^6 + (2 A_1 A_2 + 2 B_1 B_2)a^4 + (A_1^2 + B_1^2) a^2 - C =0  ,
\end{align}
which is a cubic equation in $a^2$. The discriminant of a cubic equation of the form $a x^3 + b x^2 + c x + d = 0$ is given as $\vartheta = 18 abcd - 4 b^3 d + b^2 c^2 - 4 a c^3 - 27 a^2 d^2$. The cubic equation \eqref{Eq: SS amp} has one real root when $\vartheta<0$ and three distinct real roots when $\vartheta>0$. The main four parameters $\{\alpha,E_r,f,\Delta\}$ dictate the value of coefficients $\{A_1,A_2,B_1,B_2,C\}$, the value of discriminant $\vartheta$, and thus the number of admissible steady state amplitudes. We see that for fixed values of $\{\alpha,E_r,f\}$, by sweeping the detuning parameter $\Delta$ from lower to higher excitation frequency, the stable steady state amplitude bifurcates into two stable branches and one unstable branch, where they converge back to a stable amplitude by further increasing $\Delta$. Fig. \ref{Fig: biff diag} (left) shows the bifurcation diagram by sweeping the detuning parameter $\Delta$ and for different values of $\alpha$ when $E_r=0.3$ and $f=1$. The solid and dashed black lines are the stable and unstable amplitudes, respectively. The blue lines connect the bifurcation points (red dots) for each value of $\alpha$. We see that the bifurcation points are strongly related to the value of $\alpha$, meaning that by introducing extra viscosity to the system, i.e. increasing the value of $\alpha$, the amplitudes bifurcate and then converge back faster. Figure \ref{Fig: biff diag} (right) shows the frequency response of the system, i.e. the magnitude of steady state amplitudes versus excitation frequency. As the excitation frequency is swept to the right, the steady state amplitude increases, reaches a peak value, and then jumps down (see e.g. red dashed line for $\alpha=0.4$). The peak amplitude and the jump magnitude decreases as $\alpha$ is increased. 
%
%

\begin{figure*}[h]
	\centering
	\begin{subfigure}{0.22\textwidth}
		\centering
		\includegraphics[width=1\linewidth]{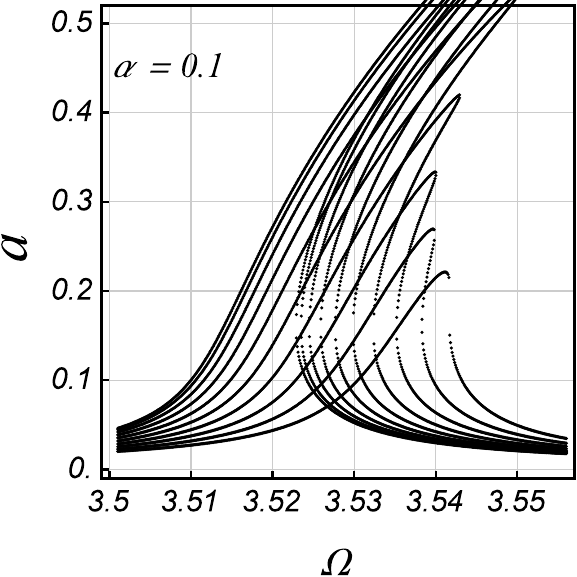}
	\end{subfigure}
	\begin{subfigure}{0.22\textwidth}
		\centering
		\includegraphics[width=1\linewidth]{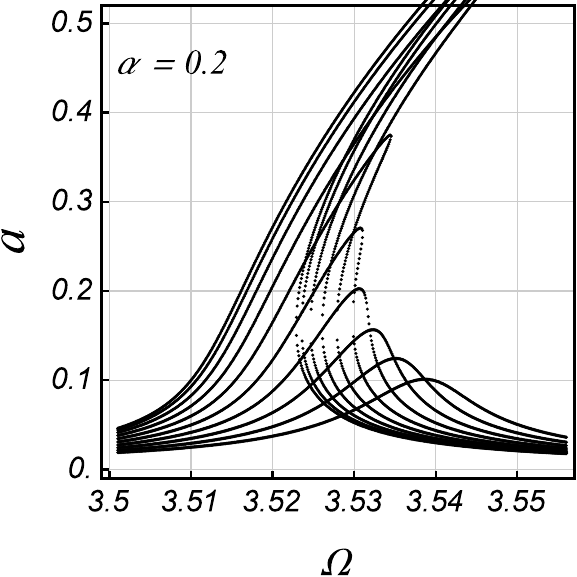}
	\end{subfigure}
	\begin{subfigure}{0.22\textwidth}
		\centering
		\includegraphics[width=1\linewidth]{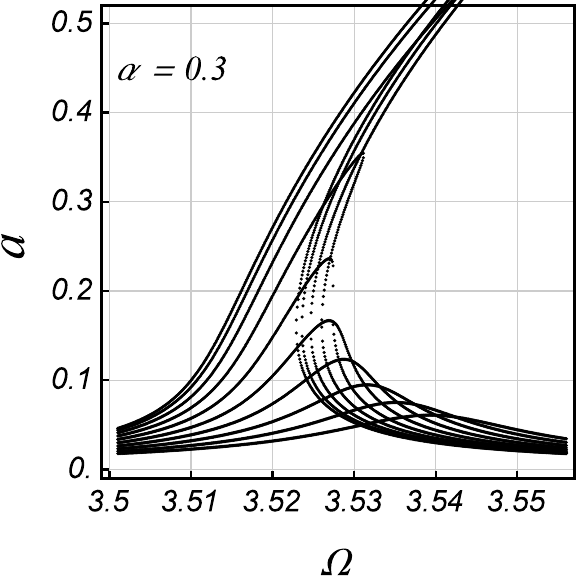}
	\end{subfigure}
	\begin{subfigure}{0.22\textwidth}
		\centering
		\includegraphics[width=1\linewidth]{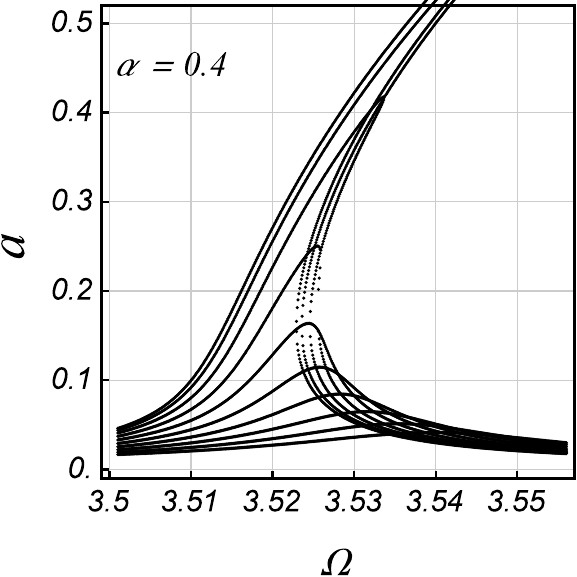}
	\end{subfigure}
	\begin{subfigure}{0.22\textwidth}
		\centering
		\includegraphics[width=1\linewidth]{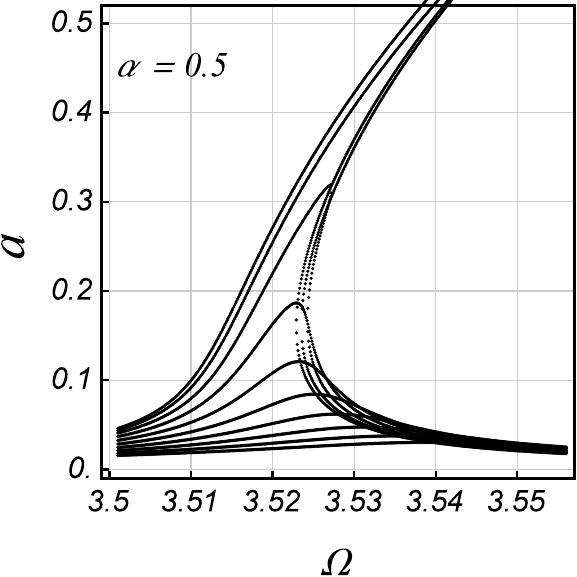}
	\end{subfigure}
	\begin{subfigure}{0.22\textwidth}
		\centering
		\includegraphics[width=1\linewidth]{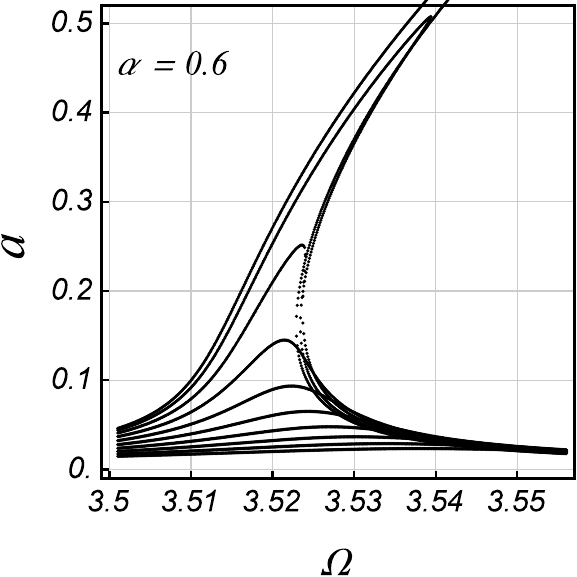}
	\end{subfigure}
	\begin{subfigure}{0.22\textwidth}
		\centering
		\includegraphics[width=1\linewidth]{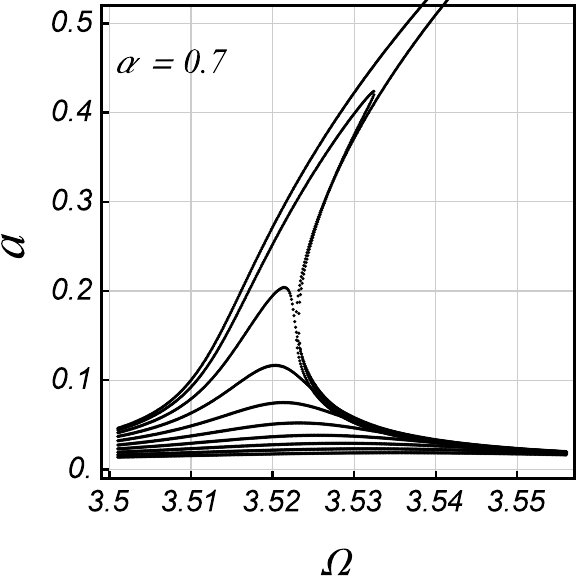}
	\end{subfigure}
	\begin{subfigure}{0.22\textwidth}
		\centering
		\includegraphics[width=1\linewidth]{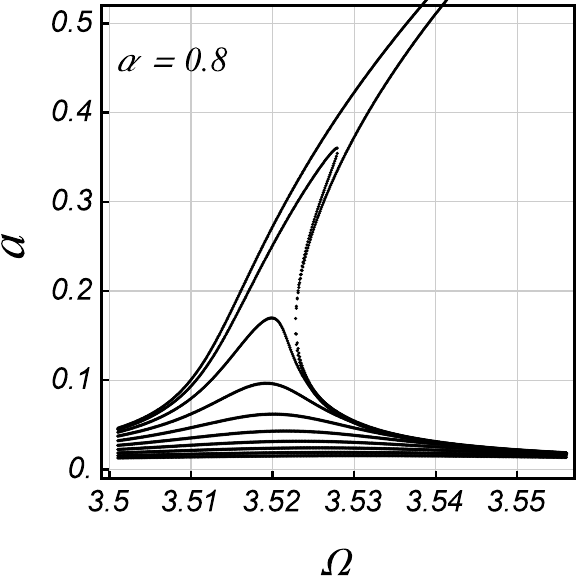}
	\end{subfigure}
	\caption{Frequency-Response curve for the case of primary resonance in the viscoelastic cantilever beam with no lumped mass at the tip. Each sub-figure corresponds to a fixed value of $\alpha$ and $f$ when $E_r = \{0.1, 0.2, \cdots, 1\}$. As effect of fractional element becomes more pronounced, i.e. $\alpha$ and $E_r$ increase, the corresponding amplitude peaks decrease and shift towards the lower frequency range.}
	\label{Fig: amp diag}
\end{figure*}

The coefficient $E_r = \frac{E_{\infty}}{E_{\alpha}}$ is the proportional contribution of fractional and pure elastic element. At a certain value while increasing this parameter, we see that the bifurcation disappears and the frequency response of system slightly changes. Fig. \ref{Fig: amp diag} shows the frequency response of the system for different values of $\{\alpha,E_r\}$ when $f = 0.5$. In each sub-figure, we let $\alpha$ be fixed and then plot the frequency response for $E_r = \{0.1,0.2,\cdots,1\}$; the amplitude peak moves down as $E_r$ is increased. For higher values of $E_r$, we see that as $\alpha$ is increased, the amplitude peaks drift back to the left, showing a softening behavior in the system response.

%

%
\subsubsection{Case 2: Lumped Mass At The Tip}
\label{Sec: case 2}
%
In this case, $M = J = 1$, and thus, given the functions $\phi_1(x)$ in Appendix \ref{Sec: App. Eigenvalue Problem of Linear Model}, the coefficients are computed as $\mathcal{M} = 1 + 70.769 J + 7.2734 M$, $\mathcal{J} = 5008.25$, $\mathcal{K}_l = \mathcal{C}_l = 98.1058$, $\mathcal{M}_b = -0.648623 - 2.69692 M$, and $\mathcal{K}_{nl} = \mathcal{C}_{nl} = 2979.66$. Similar to Case 1, we consider the following cases:

\vspace{0.2 in}
\noindent$\bullet$ Free Vibration, $F = 0$\\
Following the same steps as in Case 1, we see that the equation governing amplitude preserve its structure, but the governing equation of phase contains an extra term accommodating the $m_{nl}$.
\begin{align}
\label{Eq: free amp TipMass}
\frac{d a}{d T_1} 
& = 
- E_r \, \omega_0^{\alpha-1} \sin (\alpha \frac{\pi }{2}) 
\left( \frac{1}{2} \, c_l \, a +  \frac{3}{8} \, c_{nl} \, a^3 \right),
\\ \label{Eq: free phase TipMass}
\frac{d \varphi}{d T_1} 
 =&
\frac{1}{2} c_l \, E_r \, \omega _0^{\alpha-1} \, \cos \left(\frac{\pi  \alpha }{2}\right)  
+ \frac{3}{4} c_{nl} \, E_r \, \omega _0^{\alpha-1} \, \cos \left(\frac{\pi  \alpha }{2}\right) \, a^2 
+\frac{3}{4} \, \omega _0^{-1} \, k_{nl} \, a^2 
-\frac{1}{4} \, m_{nl} \,  \omega_0 \,  a^2.
\end{align}
This extra term does not significantly alter the behavior of phase and the whole system.

\vspace{0.2 in}
\noindent$\bullet$ Primary Resonance Case, $\Omega \approx \omega_0$\\
Similar to the free vibration, we see that the equation governing amplitude preserves its structure while the governing equation of phase contains an extra term accommodating the $m_{nl}$
\begin{align}
\label{Eq: primary amp TipMass}
\frac{d a}{d T_1} =&
- E_r \, \omega_0^{\alpha-1} \sin (\alpha \frac{\pi }{2}) 
\left( \frac{1}{2} \, c_l  \, a + \frac{3}{8} \, c_{nl}  \, a^3 \right) + \frac{1}{2} f \, \omega^{-1}_0 \, \sin (\Delta  T_1 - \varphi),
\\ \label{Eq: primary phase TipMass}
a \, \frac{d \varphi}{d T_1} =&
\frac{1}{2} c_l \, E_r \, \omega _0^{\alpha-1} \, \cos(\frac{\pi  \alpha }{2})  \, a
+ \frac{3}{4} c_{nl} \, E_r \, \omega _0^{\alpha-1} \, \cos(\frac{\pi  \alpha }{2}) \, a^3 +\frac{3}{4} \, \omega _0^{-1} \, k_{nl} \, a^3 
- \frac{1}{2} f \, \omega^{-1}_0 \, \cos (\Delta  T_1 - \varphi) -\frac{1}{4} \, m_{nl} \,  \omega_0 \,  a^3 .
\end{align}
Transforming the equations into an autonomous system by letting $ \gamma = \Delta \, T_1 - \varphi$, we obtain the governing equation of steady state solution as
\begin{align}
\label{Eq: primary ss TipMass}
&\left[
\frac{c_l}{2} E_r \, \omega_0^{\alpha-1} \sin (\frac{\pi  \alpha }{2})  \, a
+\frac{3 c_{nl}}{8}  E_r \, \omega_0^{\alpha-1} \sin (\frac{\pi  \alpha }{2}) \, a^3 
\right]^2
\left[
\left(\Delta - \frac{c_l}{2} E_r \, \omega _0^{\alpha-1} \, \cos(\frac{\pi  \alpha }{2}) \right)  a \right. \nonumber \\
& \left.- \frac{3}{4} \left( c_{nl} \, E_r \, \omega _0^{\alpha-1} \, \cos(\frac{\pi  \alpha }{2}) + \omega _0^{-1} \, k_{nl} + \frac{1}{3} \, m_{nl} \,  \omega_0  \right) a^3 
\right]^2
= \frac{f^2}{4 \, \omega^2_0},
\end{align}
which, similar to Case 1, can be written as $$ (A_2^2 + B_2^2) a^6 + (2 A_1 A_2 + 2 B_1 B_2)a^4 + (A_1^2 + B_1^2) a^2 - C =0 ,$$ where all the $A_1$, $A_2$, $B_1$, and $C$ are the same as in Case 1, but $$B_2 = - \frac{3}{4} \left( c_{nl} \, E_r \, \omega _0^{\alpha-1} \, \cos(\frac{\pi  \alpha }{2}) + \omega _0^{-1} \, k_{nl} + \frac{1}{3} \, m_{nl} \,  \omega_0 \right).$$ The corresponding cubic equation can be solved to obtain the bifurcation diagram and also the frequency response of the system. However, in addition to Case 1, we have an extra parameter $m_{nl}$ which affects the response of the system.

\section{Summary and Discussion}
\label{Sec: Summery}
In this work we investigated the anomalous nonlinear dynamics driven by the application of extraordinary materials. Our anomalous system is represented as a nonlinear fractional {Kelvin-Voigt} viscoelastic cantilever beam.
{A spectral method} was employed for spatial 
discretization of the governing equation of motion, reducing it to a set of 
nonlinear fractional {ODEs}. The corresponding system 
was linearized and the time-fractional integration was {carried out} 
through {a direct L1 finite-difference scheme}, together with a Newmark method. 
{For the nonlinear solution, a method of multiple scale was employed, and the 
time response of the beam subject to a base excitation was obtained}. We performed a set of 
numerical experiments {on the system response under varying fractional orders, representing different stages of material evolution, where we observed:}
\begin{itemize}
    \item Anomalous drift in peak amplitude response according to fractional orders, and the presence of a low-frequency critical point under linear forced vibration.
    \item Short-time and long-time anomalous behaviors under linear free vibration, respectively, for Riemann-Liouville and Caputo definitions.
	\item Super sensitivity of {the amplitude response with respect} to the fractional model parameters at free vibration.
	\item {A critical behavior of the decay rate sensitivity with respect to $\alpha$}, where 
	increasing values of fractional order yielded higher decay rates 
	(softening) before a critical value $\alpha_{cr}$. {Lower decay rates (stress hardening) were observed beyond such critical value.}
	\item A bifurcation behavior under steady-state amplitude at primary 
	resonance case.
\end{itemize}
The choice of a fractional Kelvin-Voigt model in this work allowed us to 
describe a material in the intersection between anomalous and standard 
constitutive behavior, where the contribution of the {SB element yields 
the power-law material response, while the Hookean spring reflects the 
instantaneous response of many engineering materials. In addition, the shifts in amplitude-frequency response with respect to the fractional order motivate future studies on the downscaling of fractional operators to the associated far-from-equilibrium dynamics (polymer caging/reptation, dislocation avalanches) in evolving heterogeneous microstructures \cite{Mashayekhi2019Fractal}}. In terms of modifications of the current model, 
different material distribution functions could be chosen, leading to 
application-based material design for a wide range of 
structural materials and anomalous systems, including microelectromechanical 
systems (MEMS). 
{Finally, regarding numerical discretizations, one could
utilize additional active vibration modes}, as well as faster 
time-fractional integration methods, in order to better capture the fundamental 
dynamics of the presented {system}.

\appendix       
\section{Derivation of Governing Equation Using Extended Hamilton's Principle}
\label{Sec: App. extended Hamilton}
\subsection{Equation of Motion}
We recast the integral \eqref{Eq: work var} as $\delta W = \int_{0}^{L} \int_{A} \sigma \, \delta \varepsilon \, dA \,\, ds$ for the considered cantilever beam, in which the variation of strain is $\delta \varepsilon = -\eta \, \delta \frac{\partial{\psi}}{\partial{s}}$, using \eqref{Eq: strain-curvature relation}. Therefore, by assuming the constitutive equation \eqref{Eq: frac KV}, the variation of total work is expressed as
\small{
\begin{align}
\label{Eq: total work var}
\delta w 
& = 
\int_{0}^{L} \int_{A} \left( -\eta \, E_{\infty} \, \frac{\partial{\psi}}{\partial{s}} -\eta \, E_{\alpha} \, \prescript{RL}{0}{\mathcal D}_{t}^{\alpha} \, \frac{\partial{\psi}}{\partial{s}} \right) \, (-\eta \, \delta \frac{\partial{\psi}}{\partial{s}}) \, dA  \,\, ds
\nonumber\\
& = 
\int_{0}^{L} \left( E_{\infty} \, \left( \int_{A}  \, \eta^2 dA \right) \, \frac{\partial{\psi}}{\partial{s}} 
+ E_{\alpha} \, \left( \int_{A}  \, \eta^2 dA \right) \, \prescript{RL}{0}{\mathcal D}_{t}^{\alpha} \, \frac{\partial{\psi}}{\partial{s}} \right) \, \delta \frac{\partial{\psi}}{\partial{s}} \,\, ds
\nonumber\\
& = 
\int_{0}^{L} \left( E_{\infty} \, I \, \frac{\partial{\psi}}{\partial{s}} + E_{\alpha} \, I \, \prescript{RL}{0}{\mathcal D}_{t}^{\alpha} \, \frac{\partial{\psi}}{\partial{s}} \right) \, \delta \frac{\partial{\psi}}{\partial{s}} \,\, ds
\end{align}}
where $I = \int_{A} \eta^2 \, dA$. By approximation \eqref{Eq: curvature}, we write the variation of curvature as
\begin{align}
\label{Eq: strain variation}
\delta \frac{\partial{\psi}}{\partial{s}} 
= (1 + \frac{1}{2} (\frac{\partial{v}}{\partial{s}})^2) \, \delta \frac{\partial^2{v}}{\partial{s}^2}+ \frac{\partial^2{v}}{\partial{s}^2} \, \frac{\partial{v}}{\partial{s}} \, \delta \frac{\partial{v}}{\partial{s}}.
\end{align}
Therefore, the variation of total energy becomes
\begin{align}
\label{Eq: total work var - 2}
\delta w 
& \!=\!
\int_{0}^{L}  
\left( E_{\infty}  I  \frac{\partial^2{v}}{\partial{s}^2} (1 + \frac{1}{2} (\frac{\partial{v}}{\partial{s}})^2)  
+ E_{\alpha}  I  \prescript{RL}{0}{\mathcal D}_{t}^{\alpha}  \Big[\frac{\partial^2{v}}{\partial{s}^2}(1 + \frac{1}{2} (\frac{\partial{v}}{\partial{s}})^2)\Big]\right) 
(1 + \frac{1}{2} (\frac{\partial{v}}{\partial{s}})^2)  \delta \frac{\partial^2{v}}{\partial{s}^2} \,\, ds 
\nonumber\\
& +
\int_{0}^{L}  
\left( E_{\infty}  I  \frac{\partial^2{v}}{\partial{s}^2} (1 + \frac{1}{2} (\frac{\partial{v}}{\partial{s}})^2) 
+ E_{\alpha}  I  \prescript{RL}{0}{\mathcal D}_{t}^{\alpha}  \Big[\frac{\partial^2{v}}{\partial{s}^2}(1 + \frac{1}{2} (\frac{\partial{v}}{\partial{s}})^2)\Big]\right) 
\frac{\partial^2{v}}{\partial{s}^2}  \frac{\partial{v}}{\partial{s}}  \delta \frac{\partial{v}}{\partial{s}}  \,\, ds 
\end{align}
By expanding the terms and integrating by parts, we have
\begin{align}
\label{Eq: total work var - 3}
\delta w 
& =
\int_{0}^{L}  
\frac{\partial^2}{\partial{s}^2}\left(
\left( E_{\infty}  I  \frac{\partial^2{v}}{\partial{s}^2} (1 + \frac{1}{2} (\frac{\partial{v}}{\partial{s}})^2)  
+ E_{\alpha}  I  \prescript{RL}{0}{\mathcal D}_{t}^{\alpha}  \Big[\frac{\partial^2{v}}{\partial{s}^2}(1 + \frac{1}{2} (\frac{\partial{v}}{\partial{s}})^2)\Big]\right) 
(1 + \frac{1}{2} (\frac{\partial{v}}{\partial{s}})^2) 
\right) 
\delta v \,\, ds 
 \nonumber \\
& -
\int_{0}^{L}  
\frac{\partial}{\partial{s}}\left(
\left( E_{\infty}  I  \frac{\partial^2{v}}{\partial{s}^2} (1 + \frac{1}{2} (\frac{\partial{v}}{\partial{s}})^2) 
+ E_{\alpha} I  \prescript{RL}{0}{\mathcal D}_{t}^{\alpha} \Big[\frac{\partial^2{v}}{\partial{s}^2}(1 + \frac{1}{2} (\frac{\partial{v}}{\partial{s}})^2)\Big]\right) 
\frac{\partial^2{v}}{\partial{s}^2}  \frac{\partial{v}}{\partial{s}} 
\right)
\delta v  \,\, ds 
\nonumber \\
& +
\left( E_{\infty}  I  \frac{\partial^2{v}}{\partial{s}^2} (1 + \frac{1}{2} (\frac{\partial{v}}{\partial{s}})^2)  
+ E_{\alpha}  I  \prescript{RL}{0}{\mathcal D}_{t}^{\alpha}  \Big[\frac{\partial^2{v}}{\partial{s}^2}(1 + \frac{1}{2} (\frac{\partial{v}}{\partial{s}})^2)\Big]\right) 
(1 + \frac{1}{2} (\frac{\partial{v}}{\partial{s}})^2)  \delta \frac{\partial{v}}{\partial{s}} \Bigg|_{0}^{L}
 \nonumber\\
& - \frac{\partial}{\partial{s}}\left(
\left( E_{\infty}  I  \frac{\partial^2{v}}{\partial{s}^2} (1 + \frac{1}{2} (\frac{\partial{v}}{\partial{s}})^2) 
+ E_{\alpha}  I  \prescript{RL}{0}{\mathcal D}_{t}^{\alpha}  \Big[\frac{\partial^2{v}}{\partial{s}^2}(1 + \frac{1}{2} (\frac{\partial{v}}{\partial{s}})^2)\Big]\right) 
(1 + \frac{1}{2} (\frac{\partial{v}}{\partial{s}})^2) \right) \delta v \Bigg|_{0}^{L}
 \nonumber\\
& +
\left( E_{\infty}  I  \frac{\partial^2{v}}{\partial{s}^2} (1 + \frac{1}{2} (\frac{\partial{v}}{\partial{s}})^2)  
+ E_{\alpha}  I  \prescript{RL}{0}{\mathcal D}_{t}^{\alpha} \Big[\frac{\partial^2{v}}{\partial{s}^2}(1 + \frac{1}{2} (\frac{\partial{v}}{\partial{s}})^2)\Big]\right) 
\frac{\partial^2{v}}{\partial{s}^2}  \frac{\partial{v}}{\partial{s}}  \delta v \Bigg|_{0}^{L}
\end{align}
The prescribed geometry boundary conditions at the base of the beam, $s=0$, allow the variation of deflection and its first derivative to be zero at $s=0$, i.e. $\delta v(0,t) = \delta \frac{\partial{v}}{\partial{s}}(0,t) = 0$. Therefore,
\begin{align}
\label{Eq: total work var - 4}
\delta w 
& \!=\!
\!\int_{0}^{L}\!  
\frac{\partial^2}{\partial{s}^2}\left(
\left( E_{\infty}  I  \frac{\partial^2{v}}{\partial{s}^2} (1 + \frac{1}{2} (\frac{\partial{v}}{\partial{s}})^2)  
+ E_{\alpha}  I  \prescript{RL}{0}{\mathcal D}_{t}^{\alpha}  \Big[\frac{\partial^2{v}}{\partial{s}^2}(1 + \frac{1}{2} (\frac{\partial{v}}{\partial{s}})^2)\Big]\right) 
(1 + \frac{1}{2} (\frac{\partial{v}}{\partial{s}})^2) 
\right) 
\delta v \,\, ds 
 \nonumber\\
& -
\int_{0}^{L}  
\frac{\partial}{\partial{s}}\left(
\left( E_{\infty}  I  \frac{\partial^2{v}}{\partial{s}^2} (1 + \frac{1}{2} (\frac{\partial{v}}{\partial{s}})^2)  
+ E_{\alpha}  I  \prescript{RL}{0}{\mathcal D}_{t}^{\alpha}  \Big[\frac{\partial^2{v}}{\partial{s}^2}(1 + \frac{1}{2} (\frac{\partial{v}}{\partial{s}})^2)\Big]\right) 
\frac{\partial^2{v}}{\partial{s}^2}  \frac{\partial{v}}{\partial{s}}  
\right)
\delta v  \,\, ds 
\nonumber\\
& +
\left( E_{\infty}  I  \frac{\partial^2{v}}{\partial{s}^2} (1 + \frac{1}{2} (\frac{\partial{v}}{\partial{s}})^2) \!+\! E_{\alpha}  I  \prescript{RL}{0}{\mathcal D}_{t}^{\alpha}  \Big[\frac{\partial^2{v}}{\partial{s}^2}(1 \!+\! \frac{1}{2} (\frac{\partial{v}}{\partial{s}})^2)\Big]\right) 
(1 \!+\! \frac{1}{2} (\frac{\partial{v}}{\partial{s}})^2) \Bigg|_{s=L}  \delta \frac{\partial{v}}{\partial{s}}(L,t) 
\nonumber\\
& -
\frac{\partial}{\partial{s}}\left(
\left( E_{\infty}  I  \frac{\partial^2{v}}{\partial{s}^2} (1 + \frac{1}{2} (\frac{\partial{v}}{\partial{s}})^2)  
\!+\! E_{\alpha}  I  \prescript{RL}{0}{\mathcal D}_{t}^{\alpha}  \Big[\frac{\partial^2{v}}{\partial{s}^2}(1 \!+\! \frac{1}{2} (\frac{\partial{v}}{\partial{s}})^2)\Big]\right) 
(1 \!+\! \frac{1}{2} (\frac{\partial{v}}{\partial{s}})^2) \right) \Bigg|_{s=L}  \delta v(L,t) 
\nonumber\\
& +
\left( E_{\infty}  I  \frac{\partial^2{v}}{\partial{s}^2} (1 + \frac{1}{2} (\frac{\partial{v}}{\partial{s}})^2)  +
E_{\alpha}  I  \prescript{RL}{0}{\mathcal D}_{t}^{\alpha}  \Big[\frac{\partial^2{v}}{\partial{s}^2}(1 + \frac{1}{2} (\frac{\partial{v}}{\partial{s}})^2\Big]\right) 
\frac{\partial^2{v}}{\partial{s}^2}  \frac{\partial{v}}{\partial{s}} \Bigg|_{s=L}  \delta v(L,t) 
\end{align}
Let $\varrho$ be mass per unit volume of the beam, $M$ and $J$ be the mass and rotatory inertia of the lumped mass at the tip of beam. By considering the displacement and velocity of the beam given in \eqref{Eq: displacement} and \eqref{Eq: velocity}, respectively, the kinetic energy is obtained as
\begin{align}
\label{Eq: kinetic}
T = &\frac{1}{2} \, \int_{0}^{L} \int_{A} \varrho  \, (\frac{\partial{\textbf{r}}}{\partial{t}})^2 \, dA \, ds + \frac{1}{2}M \left( (\frac{\partial{u}}{\partial{t}})^2  + (\frac{\partial{v}}{\partial{t}} + \dot {{v_b}})^2 \right) \Big|_{s=L} + \frac{1}{2} J (\frac{\partial{\psi}}{\partial{t}})^2 \Big|_{s=L} , 
 \nonumber \\
 = &\frac{1}{2} \, \int_{0}^{L} \int_{A} \varrho  \, \Big\{ (\frac{\partial{u}}{\partial{t}} - \eta \, \frac{\partial{\psi}}{\partial{t}} \, \cos(\psi) )^2 
+ (\frac{\partial{v}}{\partial{t}}+ \dot {{v_b}} - \eta \, \frac{\partial{\psi}}{\partial{t}} \, \sin(\psi))^2 \Big\}  \,  dA \, ds 
+ \frac{1}{2}M \left((\frac{\partial{u}}{\partial{t}})^2 + (\frac{\partial{v}}{\partial{t}} + \dot {{v_b}})^2 \right) \Big|_{s=L} + \frac{1}{2} J (\frac{\partial{\psi}}{\partial{t}})^2 \Big|_{s=L} , 
 \nonumber\\
= &
\frac{1}{2} \, \int_{0}^{L} \int_{A} \varrho  \,
\Big\{ (\frac{\partial{u}}{\partial{t}} )^2 - 2 \eta \, \frac{\partial{u}}{\partial{t}}  \, \frac{\partial{\psi}}{\partial{t}}  \, \cos(\psi)  +\eta^2 \, (\frac{\partial{\psi}}{\partial{t}})^2 \, \cos^2(\psi) 
+ (\frac{\partial{v}}{\partial{t}} )^2 + {\dot {{v_b}}}^2 + 2 \frac{\partial{v}}{\partial{t}}  \, \dot {{v_b}}  \nonumber \\
& - 2 \eta \, (\frac{\partial{v}}{\partial{t}} \!+\! \dot {{v_b}}) \frac{\partial{\psi}}{\partial{t}}  \, \sin(\psi) \!+\! \eta^2 (\frac{\partial{\psi}}{\partial{t}})^2 \, \sin^2(\psi) 
\Big\}   dA  ds 
+ \frac{1}{2}M \left((\frac{\partial{u}}{\partial{t}})^2 + (\frac{\partial{v}}{\partial{t}} + \dot {{v_b}})^2 \right) \Big|_{s=L}+ \frac{1}{2} J (\frac{\partial{\psi}}{\partial{t}})^2 \Big|_{s=L} , 
\nonumber\\
 = &
\frac{1}{2} \, \int_{0}^{L} \int_{A} \varrho  \,
\Big\{ (\frac{\partial{u}}{\partial{t}})^2 + (\frac{\partial{v}}{\partial{t}})^2 + {\dot {{v_b}}}^2 + 2 \frac{\partial{v}}{\partial{t}} \, \dot {{v_b}}  - 2 \eta \, \frac{\partial{u}}{\partial{t}}\, \frac{\partial{\psi}}{\partial{t}} \, \cos(\psi) + \eta^2 \, (\frac{\partial{\psi}}{\partial{t}})^2 
- 2 \eta \, (\frac{\partial{v}}{\partial{t}} + \dot {{v_b}}) \frac{\partial{\psi}}{\partial{t}} \, \sin(\psi)  
\Big\}  \,  dA \, ds \nonumber \\
& 
+ \frac{1}{2}M \left((\frac{\partial{u}}{\partial{t}} )^2 + (\frac{\partial{v}}{\partial{t}}  + \dot {{v_b}})^2 \right) \Big|_{s=L} + \frac{1}{2} J (\frac{\partial{\psi}}{\partial{t}})^2 \Big|_{s=L} .
\end{align}
Let
\begin{align*}
\rho = \int_{A} \varrho dA , \quad
\mathcal{J}_1 = \int_{A} \varrho \, \eta \,  dA , \quad
\mathcal{J}_2 = \int_{A} \varrho \, \eta^2 dA.
\end{align*}
$\rho$ is the mass per unit length of the beam, $\mathcal{J}_1$ is the first moment of inertia and is zero because the reference point of coordinate system attached to the cross section coincides with the mass centroid, and $\mathcal{J}_2$ is the second moment of inertia, which is very small for slender beam and can be ignored \cite{hamdan1997large}. Assuming that the velocity along the length of the beam, $\dot u$, is relatively small compared to the lateral velocity $\dot v + \dot {{v_b}}$, the kinetic energy of the beam can be reduced to
\begin{align}
\label{Eq: kinetic - 2}
T = &
\frac{1}{2} \, \rho \int_{0}^{L} ( \frac{\partial{v}}{\partial{t}} + \dot {{v_b}})^2  \,   ds 
+ \frac{1}{2}M (\frac{\partial{v}}{\partial{t}} + \dot {{v_b}})^2  \Big|_{s=L}
+ \frac{1}{2} J (\frac{\partial{\psi}}{\partial{t}})^2 \Big|_{s=L}, 
\end{align}
where its variation can be taken as
\begin{align}
\label{Eq: kinetic variation}
\delta T = & \rho \int_{0}^{L} (\frac{\partial{v}}{\partial{t}} + \dot{{{v_b}}}) \, \delta \frac{\partial{v}}{\partial{t}} \, ds + M (\frac{\partial{v}}{\partial{t}} + \dot{{{v_b}}}) \, \delta \frac{\partial{v}}{\partial{t}} \Big|_{s=L}
+  J \frac{\partial{\psi}}{\partial{t}} \, \delta \frac{\partial{\psi}}{\partial{t}} \Big|_{s=L} ,
\end{align}
in which $\frac{\partial{\psi}}{\partial{t}}$ is given in \eqref{Eq: angular velocity} and $\delta \frac{\partial{\psi}}{\partial{t}}$ can be obtained as 
\begin{align*}
	\delta \frac{\partial{\psi}}{\partial{t}} \simeq  ( 1 + \frac{1}{2} (\frac{\partial{v}}{\partial{s}})^2) \delta \frac{\partial^2{v}}{\partial{t}\partial{s}}+ \frac{\partial{v}}{\partial{s}} \frac{\partial^2{v}}{\partial{t}\partial{s}} \delta \frac{\partial{v}}{\partial{s}}.
\end{align*}
\ Therefore, 
\begin{align}
\label{Eq: kinetic variation - 2}
\delta T 
 \simeq &\rho \int_{0}^{L} (\frac{\partial{v}}{\partial{t}} + \dot{{{v_b}}})  \delta \frac{\partial{v}}{\partial{t}}  ds 
+ M (\frac{\partial{v}}{\partial{t}} + \dot{{{v_b}}})  \delta \frac{\partial{v}}{\partial{t}} \Big|_{s=L}  
\!+\! J \left( \frac{\partial^2{v}}{\partial{t}\partial{s}} ( 1 \!+\! (\frac{\partial{v}}{\partial{s}})^2) \delta \frac{\partial^2{v}}{\partial{t}\partial{s}} \!+\! \frac{\partial{v}}{\partial{s}} (\frac{\partial^2{v}}{\partial{t}\partial{s}})^2 \delta \frac{\partial{v}}{\partial{s}} \right) \Big|_{s=L} .
\end{align}
The time integration of $\delta T$ takes the following form through integration by parts
\begin{align}
\label{Eq: kinetic var - 2}
\int_{t_1}^{t_2} \delta T \, dt = & \int_{t_1}^{t_2} 
\Bigg\{
\rho \int_{0}^{L} (\frac{\partial{v}}{\partial{t}} + \dot{{{v_b}}}) \, \delta \frac{\partial{v}}{\partial{t}} \, ds 
\!+\! M (\frac{\partial{v}}{\partial{t}} + \dot{{{v_b}}}) \, \delta \frac{\partial{v}}{\partial{t}} \Big|_{s=L} 
+ J \left( \frac{\partial^2{v}}{\partial{t}\partial{s}} ( 1 + (\frac{\partial{v}}{\partial{s}})^2) \delta \frac{\partial^2{v}}{\partial{t}\partial{s}} + \frac{\partial{v}}{\partial{s}} (\frac{\partial^2{v}}{\partial{t}\partial{s}})^2 \delta \frac{\partial{v}}{\partial{s}} \right) \Big|_{s=L}
\Bigg\} \, dt
 \nonumber\\
= &
\int_{t_1}^{t_2} \rho \int_{0}^{L} (\frac{\partial{v}}{\partial{t}} \!+\! \dot{{{v_b}}}) \, \delta \frac{\partial{v}}{\partial{t}} \, ds \, dt
\!+\! \!M\! \int_{t_1}^{t_2} (\frac{\partial{v}}{\partial{t}} \!+\! \dot{{{v_b}}}) \, \delta \frac{\partial{v}}{\partial{t}} \Big|_{s=L} \, dt
\!+\! \!J\! \int_{t_1}^{t_2} \left( \frac{\partial^2{v}}{\partial{t}\partial{s}} ( 1 + (\frac{\partial{v}}{\partial{s}})^2) \delta \frac{\partial^2{v}}{\partial{t}\partial{s}} + \frac{\partial{v}}{\partial{s}} (\frac{\partial^2{v}}{\partial{t}\partial{s}})^2 \delta \frac{\partial{v}}{\partial{s}}  \right) \Big|_{s=L} \, dt
 \nonumber\\
= &
\rho \int_{0}^{L} \int_{t_1}^{t_2}  (\frac{\partial{v}}{\partial{t}} + \dot{{{v_b}}}) \, \delta \frac{\partial{v}}{\partial{t}} \, dt \, ds
+ M \int_{t_1}^{t_2} (\frac{\partial{v}}{\partial{t}} + \dot{{{v_b}}}) \, \delta \frac{\partial{v}}{\partial{t}}  \, dt \, \Big|_{s=L}
\!+\! J \int_{t_1}^{t_2} \left( \frac{\partial^2{v}}{\partial{t}\partial{s}} ( 1 + (\frac{\partial{v}}{\partial{s}})^2) \delta \frac{\partial^2{v}}{\partial{t}\partial{s}} + \frac{\partial{v}}{\partial{s}} (\frac{\partial^2{v}}{\partial{t}\partial{s}})^2 \delta \frac{\partial{v}}{\partial{s}} \right)  \, dt \, \Big|_{s=L}
 \nonumber\\
= &
\rho \!\int_{0}^{L}\! \left[ (\frac{\partial{v}}{\partial{t}} \!+\! \dot{{{v_b}}})  \delta{v} \Big|_{t_1}^{t_2} \!-\! \!\int_{t_1}^{t_2}\!  (\frac{\partial^2{v}}{\partial{t}^2} \!+\! \ddot{{{v_b}}})  \delta {v}  dt \right]  ds
\!+\! M (\frac{\partial{v}}{\partial{t}} \!+\! \dot{{{v_b}}})  \delta v \Big|_{s=L} \Big|_{t_1}^{t_2} 
\!-\! M \!\int_{t_1}^{t_2} (\frac{\partial^2{v}}{\partial{t}^2} \!+\! \ddot{{{v_b}}})  \delta v   dt  \Big|_{s=L}
 \nonumber\\
& 
\!+\! J  \frac{\partial^2{v}}{\partial{t}\partial{s}} ( 1 \!+\! (\frac{\partial{v}}{\partial{s}})^2) \delta \frac{\partial{v}}{\partial{s}}  \Big|_{s=L}  \Big|_{t_1}^{t_2}
\!-\! J \int_{t_1}^{t_2} 
\left(  \frac{\partial^3{v}}{\partial{t}^2\partial{s}} ( 1 \!+\! (\frac{\partial{v}}{\partial{s}})^2) \!+\! \frac{\partial{v}}{\partial{s}} (\frac{\partial^2{v}}{\partial{t}\partial{s}})^2 \right)  
\delta \frac{\partial{v}}{\partial{s}}   dt  \Big|_{s=L}
\nonumber\\
\!=\! &
- \!\int_{t_1}^{t_2}\! \Bigg\{
\rho \int_{0}^{L}  (\frac{\partial^2{v}}{\partial{t}^2} \!+\! \ddot{{{v_b}}})  \delta {v}  ds 
\!+\! M (\frac{\partial^2{v}}{\partial{t}^2} \!+\! \ddot{{{v_b}}})  \delta v   \Big|_{s=L}
+ J \left(  \frac{\partial^3{v}}{\partial{t}^2\partial{s}}( 1 \!+\! (\frac{\partial{v}}{\partial{s}})^2) \!+\! \frac{\partial{v}}{\partial{s}} (\frac{\partial^2{v}}{\partial{t}\partial{s}})^2  \right)  \delta \frac{\partial{v}}{\partial{s}}   \Big|_{s=L}
\Bigg\}  dt,
\end{align}
where we consider that $\delta v = \delta \frac{\partial{v}}{\partial{s}} = 0$ at $t=t_1$ and $t=t_2$. Therefore, the extended Hamilton's principle takes the form
\begin{align}
\label{Eq: extended Hamilton}
&
\int_{t_1}^{t_2} \Bigg\{
\int_{0}^{L} \Bigg[
- \rho (\frac{\partial^2{v}}{\partial{t}^2}  + \ddot{{{v_b}}}) 
-\frac{\partial^2}{\partial{s}^2} \left(
\left( E_{\infty} \, I \, \frac{\partial^2{v}}{\partial{s}^2} (1 + \frac{1}{2} (\frac{\partial{v}}{\partial{s}})^2)  
+ E_{\alpha} \, I \, \prescript{RL}{0}{\mathcal D}_{t}^{\alpha} \, \Big[\frac{\partial^2{v}}{\partial{s}^2}(1 + \frac{1}{2} (\frac{\partial{v}}{\partial{s}})^2)\Big]\right) 
(1 + \frac{1}{2} (\frac{\partial{v}}{\partial{s}})^2) 
\right)
 \nonumber\\
&
+\frac{\partial}{\partial{s}}\left( 
\left( E_{\infty} \, I \, \frac{\partial^2{v}}{\partial{s}^2} (1 + \frac{1}{2} (\frac{\partial{v}}{\partial{s}})^2) + E_{\alpha} \, I \, \prescript{RL}{0}{\mathcal D}_{t}^{\alpha} \, \Big[\frac{\partial^2{v}}{\partial{s}^2}(1 + \frac{1}{2} (\frac{\partial{v}}{\partial{s}})^2)\Big]\right) 
\frac{\partial^2{v}}{\partial{s}^2} \, \frac{\partial{v}}{\partial{s}}\, 
\right)
\Bigg] \, \delta {v} \, ds
- M (\frac{\partial^2{v}}{\partial{t}^2} + \ddot{{{v_b}}}) \Big|_{s=L} \, \delta v(L,t) \nonumber\\&
\!-\! J \left(\frac{\partial^3{v}}{\partial{t}^2\partial{s}}(1\!+\! (\frac{\partial{v}}{\partial{s}})^2)\!+\! \frac{\partial{v}}{\partial{s}} (\frac{\partial^2{v}}{\partial{t}\partial{s}})^2 \right) \Big|_{s=L} \,  \delta \frac{\partial{v}}{\partial{s}}(L,t) 
\!-\!
\left( E_{\infty} \, I \, \frac{\partial^2{v}}{\partial{s}^2} (1 \!+\! \frac{1}{2} (\frac{\partial{v}}{\partial{s}})^2) \!+\!
E_{\alpha} I  \prescript{RL}{0}{\mathcal D}_{t}^{\alpha}  \Big[\frac{\partial^2{v}}{\partial{s}^2}(1 \!+\! \frac{1}{2} (\frac{\partial{v}}{\partial{s}})^2)\Big]\right) 
(1 \!+\! \frac{1}{2} (\frac{\partial{v}}{\partial{s}})^2) \Bigg|_{s=L}  \delta \frac{\partial{v}}{\partial{s}}(L,t) 
\nonumber\\
&  +
\frac{\partial}{\partial{s}}\left(
\left( E_{\infty} \, I \, \frac{\partial^2{v}}{\partial{s}^2} (1 + \frac{1}{2} (\frac{\partial{v}}{\partial{s}})^2)  +
E_{\alpha} I  \prescript{RL}{0}{\mathcal D}_{t}^{\alpha}  \Big[\frac{\partial^2{v}}{\partial{s}^2}(1 \!+\! \frac{1}{2} (\frac{\partial{v}}{\partial{s}})^2)\Big]\right) 
(1 \!+\! \frac{1}{2} (\frac{\partial{v}}{\partial{s}})^2) \right) \Bigg|_{s=L}  \delta v(L,t) 
 \nonumber\\
&  -
\left( E_{\infty} \, I \, \frac{\partial^2{v}}{\partial{s}^2} (1 + \frac{1}{2} (\frac{\partial{v}}{\partial{s}})^2)  
+ E_{\alpha}  I  \prescript{RL}{0}{\mathcal D}_{t}^{\alpha} \Big[\frac{\partial^2{v}}{\partial{s}^2}(1 \!+\! \frac{1}{2} (\frac{\partial{v}}{\partial{s}})^2)\Big]\right) 
\frac{\partial^2{v}}{\partial{s}^2} \frac{\partial{v}}{\partial{s}} \Bigg|_{s=L} \delta v(L,t) 
\Bigg\}  dt \!=\! 0
\end{align}
Invoking the arbitrariness of virtual displacement $\delta v$, we obtain the strong form of the equation of motion as:
\begin{align}
\label{Eq: eqn of motion}
& \rho \, \frac{\partial^2{v}}{\partial{t}^2}  
+ 
E_{\infty} \, I \, \frac{\partial^2}{\partial{s}^2}\left(  \frac{\partial^2{v}}{\partial{s}^2} (1 + \frac{1}{2} (\frac{\partial{v}}{\partial{s}})^2)^2 \right) + 
E_{\alpha} \, I \, \frac{\partial^2}{\partial{s}^2}\left(  (1 + \frac{1}{2} (\frac{\partial{v}}{\partial{s}})^2) \, \prescript{RL}{0}{\mathcal D}_{t}^{\alpha} \, \Big[\frac{\partial^2{v}}{\partial{s}^2}(1 + \frac{1}{2} (\frac{\partial{v}}{\partial{s}})^2) \Big]\right)
 \nonumber\\
& - 
E_{\infty} \, I \, \frac{\partial}{\partial{s}}\left(  \frac{\partial{v}}{\partial{s}} \,  (\frac{\partial^2{v}}{\partial{s}^2})^2 (1 + \frac{1}{2} (\frac{\partial{v}}{\partial{s}})^2)  \right)
- 
E_{\alpha} \, I \, \frac{\partial}{\partial{s}} \left(  \frac{\partial{v}}{\partial{s}} \, \frac{\partial^2{v}}{\partial{s}^2} \,  \prescript{RL}{0}{\mathcal D}_{t}^{\alpha} \, \Big[\frac{\partial^2{v}}{\partial{s}^2}(1 + \frac{1}{2} (\frac{\partial{v}}{\partial{s}})^2)\Big] \right) 
= -\rho \,\ddot{{{v_b}}} ,
\end{align}
which is subject to the following natural boundary conditions:
\begin{align}
\label{Eq: natural bc}
& 
J \left(  \frac{\partial^3{v}}{\partial{t}^2\partial{s}}  ( 1 + (\frac{\partial{v}}{\partial{s}})^2) + \frac{\partial{v}}{\partial{s}} (\frac{\partial^2{v}}{\partial{t}\partial{s}})^2   \right)
+ E_{\infty} \, I \, \frac{\partial^2{v}}{\partial{s}^2} (1 + \frac{1}{2} (\frac{\partial{v}}{\partial{s}})^2)^2  
+ E_{\alpha} \, I \, (1 + \frac{1}{2} (\frac{\partial{v}}{\partial{s}})^2) \, \prescript{RL}{0}{\mathcal D}_{t}^{\alpha} \, \Big[\frac{\partial^2{v}}{\partial{s}^2}(1 + \frac{1}{2} (\frac{\partial{v}}{\partial{s}})^2)\Big]
\,\, \Bigg|_{s=L} = 0 ,
 \nonumber\\
& 
M (\frac{\partial^2{v}}{\partial{t}^2}  + \ddot{{{v_b}}})
- \frac{\partial}{\partial{s}}\left(
E_{\infty} \, I \, \frac{\partial^2{v}}{\partial{s}^2} (1 + \frac{1}{2} (\frac{\partial{v}}{\partial{s}})^2)^2  
+ E_{\alpha} \, I \, (1 + \frac{1}{2} (\frac{\partial{v}}{\partial{s}})^2) \, \prescript{RL}{0}{\mathcal D}_{t}^{\alpha} \, \Big[\frac{\partial^2{v}}{\partial{s}^2}(1 + \frac{1}{2} (\frac{\partial{v}}{\partial{s}})^2)\Big] \right)
 \nonumber\\
&
+ \left( E_{\infty} \, I \, \frac{\partial{v}}{\partial{s}} \, (\frac{\partial^2{v}}{\partial{s}^2})^2 (1 + \frac{1}{2} (\frac{\partial{v}}{\partial{s}})^2) 
+ E_{\alpha} \, I \, \frac{\partial{v}}{\partial{s}} \, \frac{\partial^2{v}}{\partial{s}^2}  \, \prescript{RL}{0}{\mathcal D}_{t}^{\alpha} \, \Big[\frac{\partial^2{v}}{\partial{s}^2}(1 + \frac{1}{2} (\frac{\partial{v}}{\partial{s}})^2) \Big]\right) 
\,\, \Bigg|_{s=L}  = 0 .
\end{align}
Following a similar approach as in \eqref{Eq: curvature} in deriving the beam curvature, we obtain the approximations below, where we only consider up to third order terms and remove the higher order terms (HOTs). 
\begin{align*}
&\frac{\partial^2{v}}{\partial{s}^2} (1 + \frac{1}{2} (\frac{\partial{v}}{\partial{s}})^2)^2 
= 
\frac{\partial^2{v}}{\partial{s}^2} + \frac{\partial^2{v}}{\partial{s}^2} (\frac{\partial{v}}{\partial{s}})^2+ \text{HOTs}
\\
&(1 + \frac{1}{2} (\frac{\partial{v}}{\partial{s}})^2) \, \prescript{RL}{0}{\mathcal D}_{t}^{\alpha} \, \Big[\frac{\partial^2{v}}{\partial{s}^2}(1 + \frac{1}{2} (\frac{\partial{v}}{\partial{s}})^2)\Big] 
= 
\prescript{RL}{0}{\mathcal D}_{t}^{\alpha}  \Big[\frac{\partial^2{v}}{\partial{s}^2}(1 \!+\! \frac{1}{2} (\frac{\partial{v}}{\partial{s}})^2)\Big] \!+\!\frac{1}{2} (\frac{\partial{v}}{\partial{s}})^2  \prescript{RL}{0}{\mathcal D}_{t}^{\alpha}  \frac{\partial^2{v}}{\partial{s}^2} \!+\! \text{HOTs}\\
&
\frac{\partial{v}}{\partial{s}} \,  (\frac{\partial^2{v}}{\partial{s}^2})^2 (1 + \frac{1}{2} (\frac{\partial{v}}{\partial{s}})^2) 
= \frac{\partial{v}}{\partial{s}} \,  (\frac{\partial^2{v}}{\partial{s}^2})^2 + \text{HOTs}\\
&
\frac{\partial{v}}{\partial{s}} \, \frac{\partial^2{v}}{\partial{s}^2} \,  \prescript{RL}{0}{\mathcal D}_{t}^{\alpha} \, \Big[\frac{\partial^2{v}}{\partial{s}^2}(1 + \frac{1}{2} (\frac{\partial{v}}{\partial{s}})^2)\Big] = \frac{\partial{v}}{\partial{s}}\, \frac{\partial^2{v}}{\partial{s}^2} \,  \prescript{RL}{0}{\mathcal D}_{t}^{\alpha} \, \frac{\partial^2{v}}{\partial{s}^2} + \text{HOTs}
\end{align*}
Therefore, the strong form can be approximated up to the third order and the problem then reads as: find $v \in V$ such that
\begin{align}
\label{Eq: eqn of motion - 1}
& m \, \frac{\partial^2{v}}{\partial{t}^2}  
+ 
\frac{\partial^2}{\partial{s}^2}\left(  \frac{\partial^2{v}}{\partial{s}^2} + \frac{\partial^2{v}}{\partial{s}^2} (\frac{\partial{v}}{\partial{s}})^2 \right)
- 
\frac{\partial}{\partial{s}}\left(  \frac{\partial{v}}{\partial{s}} \,  (\frac{\partial^2{v}}{\partial{s}^2})^2  \right)
+ 
E_r  \frac{\partial^2}{\partial{s}^2}\left(  \prescript{RL}{0}{\mathcal D}_{t}^{\alpha}  \Big[\frac{\partial^2{v}}{\partial{s}^2}(1 + \frac{1}{2} (\frac{\partial{v}}{\partial{s}})^2)\Big] + \frac{1}{2}  (\frac{\partial{v}}{\partial{s}})^2  \prescript{RL}{0}{\mathcal D}_{t}^{\alpha}  \frac{\partial^2{v}}{\partial{s}^2} \right)\nonumber\\
&- 
E_r  \frac{\partial}{\partial{s}}\left(  \frac{\partial{v}}{\partial{s}} \frac{\partial^2{v}}{\partial{s}^2}   \prescript{RL}{0}{\mathcal D}_{t}^{\alpha} \frac{\partial^2{v}}{\partial{s}^2} \right) 
= -m \ddot{{{v_b}}} ,
\end{align}
By {rearranging}
\begin{align}
\label{Eq: eqn of motion - 2}
& m \, \frac{\partial^2{v}}{\partial{t}^2} 
+ 
\frac{\partial^2}{\partial{s}^2}\left(  
\frac{\partial^2{v}}{\partial{s}^2} + \frac{\partial^2{v}}{\partial{s}^2} (\frac{\partial{v}}{\partial{s}})^2 
+ E_r
\prescript{RL}{0}{\mathcal D}_{t}^{\alpha} \, \Big[\frac{\partial^2{v}}{\partial{s}^2}(1 + \frac{1}{2} (\frac{\partial{v}}{\partial{s}})^2)\Big] 
+  \frac{1}{2}  E_r 
(\frac{\partial{v}}{\partial{s}})^2 \, \prescript{RL}{0}{\mathcal D}_{t}^{\alpha} \, \frac{\partial^2{v}}{\partial{s}^2}
\right)
 \nonumber\\
& 
- 
\frac{\partial}{\partial{s}} \left( 
\frac{\partial{v}}{\partial{s}} \,  (\frac{\partial^2{v}}{\partial{s}^2})^2  
+ E_r
\frac{\partial{v}}{\partial{s}} \, \frac{\partial^2{v}}{\partial{s}^2} \,  \prescript{RL}{0}{\mathcal D}_{t}^{\alpha} \, \frac{\partial^2{v}}{\partial{s}^2}
\right)
= -m \,\ddot{{{v_b}}} ,
\end{align}
subject to the following boundary conditions:
\begin{align}
\label{Eq: bc}
& v  \Big|_{s=0} = \frac{\partial{v}}{\partial{s}}  \Big|_{s=0} = 0 ,
\nonumber\\
& 
\frac{J m}{\rho} \left(  \frac{\partial^3{v}}{\partial{t}^2\partial{s}} ( 1 + (\frac{\partial{v}}{\partial{s}})^2) + \frac{\partial{v}}{\partial{s}} (\frac{\partial^2{v}}{\partial{t}\partial{s}})^2  \right)
\left( 
\frac{\partial^2{v}}{\partial{s}^2} + \frac{\partial^2{v}}{\partial{s}^2} (\frac{\partial{v}}{\partial{s}})^2
+ E_r  \prescript{RL}{0}{\mathcal D}_{t}^{\alpha}  \Big[\frac{\partial^2{v}}{\partial{s}^2}(1 + \frac{1}{2} (\frac{\partial{v}}{\partial{s}})^2)\Big] 
+ \frac{1}{2}  E_r  (\frac{\partial{v}}{\partial{s}})^2  \prescript{RL}{0}{\mathcal D}_{t}^{\alpha}  \frac{\partial^2{v}}{\partial{s}^2}  
\right)  \Bigg|_{s=L} = 0 ,
 \nonumber\\
& 
\frac{M m}{\rho} (\frac{\partial^2{v}}{\partial{t}^2} \!+\! \ddot{{{v_b}}})
\!-\! \frac{\partial}{\partial{s}} \left(
\frac{\partial^2{v}}{\partial{s}^2} \!+\! \frac{\partial^2{v}}{\partial{s}^2} (\frac{\partial{v}}{\partial{s}})^2
\!+\! E_r  \prescript{RL}{0}{\mathcal D}_{t}^{\alpha}  \Big[\frac{\partial^2{v}}{\partial{s}^2}(1 \!+\! \frac{1}{2} (\frac{\partial{v}}{\partial{s}})^2) \Big]
\!+\! \frac{E_r}{2} (\frac{\partial{v}}{\partial{s}})^2  \prescript{RL}{0}{\mathcal D}_{t}^{\alpha}  \frac{\partial^2{v}}{\partial{s}^2}  
\right)
\!+\! \left( 
\frac{\partial{v}}{\partial{s}}  (\frac{\partial^2{v}}{\partial{s}^2})^2  \!+\! E_r  \frac{\partial{v}}{\partial{s}}   \frac{\partial^2{v}}{\partial{s}^2}   \prescript{RL}{0}{\mathcal D}_{t}^{\alpha}  \frac{\partial^2{v}}{\partial{s}^2}	
\right)  \Bigg|_{s=L}  \!=\! 0 ,
\end{align}
where $m = \frac{\rho}{ E_{\infty} \, I}$ and $E_r = \frac{E_{\alpha}}{E_{\infty}}$. 

\bigskip
\subsection{Nondimensionalization}
%
Let the {following} dimensionless variables:
\begin{align}\label{eq:dimenionsless_variables}
&s^{*} \!=\! \frac{s}{L}, \,\,\,
v^{*} \!=\! \frac{v}{L}, \,\,\,
t^{*} \!=\! t \left(\frac{1}{m L^4}\right)^{1/2},\,\,\,
E_r^{*}\! =\! E_r \left(\frac{1}{m L^4}\right)^{\alpha/2},\,\,\,J^{*} \!=\! \frac{J}{\rho L^3}, \,\,\,
M^{*} \!=\! \frac{M}{\rho L},\,\,\,
{{v_b}}^* \!=\!  \frac{{{v_b}}}{L}.
\end{align}
We obtain the following dimensionless equation by substituting the above dimensionless variables. 
\begin{align}
\label{Eq: eqn of motion - dimless - 0}
& m \frac{L}{mL^4} \, \frac{\partial^2 v^*}{\partial {t^*}^2}  
+
\frac{1}{L^2} \frac{\partial^2}{\partial {s^*}^2} 
\Bigg[  
\frac{L}{L^2} \frac{\partial^2 v^*}{\partial {s^*}^2} 
+ \frac{L}{L^2} \frac{\partial^2 v^*}{\partial {s^*}^2} 
( \frac{L}{L} \frac{\partial v^*}{\partial {s^*}} )^2 
+ \frac{E_r^* (mL^4)^{\alpha/2}}{2} \frac{1}{(mL^4)^{\alpha/2}} \frac{L}{L^2} (\frac{L}{L})^2 \prescript{RL}{0}{\mathcal D}_{t^*}^{\alpha} \frac{\partial^2 v^*}{\partial {s^*}^2}(\frac{\partial v^*}{\partial {s^*}})^2 
 \nonumber\\
&
+ E_r^* (mL^4)^{\alpha/2} \frac{1}{(mL^4)^{\alpha/2}} \frac{L}{L^2} \prescript{RL}{0}{\mathcal D}_{t^*}^{\alpha} \frac{\partial^2 v^*}{\partial {s^*}^2}
+ \frac{1}{2}  E_r^* (mL^4)^{\alpha/2} (\frac{L}{L} \frac{\partial v^*}{\partial s^*})^2  \frac{1}{(mL^4)^{\alpha/2}} \frac{L}{L^2} \prescript{RL}{0}{\mathcal D}_{t^*}^{\alpha} \frac{\partial^2 v^*}{\partial {s^*}^2}
\Bigg]
 \nonumber\\
& - 
\frac{1}{L} \frac{\partial}{\partial {s^*}} 
\Bigg[
\frac{L}{L} \frac{\partial v^*}{\partial {s^*}} (\frac{L}{L^2} \frac{\partial^2 v^*}{\partial {s^*}^2})^2 
+ E_r^* (mL^4)^{\alpha/2} \frac{L}{L} \frac{\partial v^*}{\partial s^*} \frac{L}{L^2} \frac{\partial^2 v^*}{\partial {s^*}^2}  \frac{1}{(mL^4)^{\alpha/2}} \frac{L}{L^2} \prescript{RL}{0}{\mathcal D}_{t^*}^{\alpha} \frac{\partial^2 v^*}{\partial {s^*}^2}
\Bigg]
= -m \frac{L}{mL^4} \, \frac{\partial^2 {{v_b}}^*}{\partial {t^*}^2}  ,
\end{align}
which can be simplified to
\begin{align}
\label{Eq: eqn of motion - dimless - 1}
&\frac{\partial^2 v^*}{\partial {t^*}^2}  \!+\! 
\frac{\partial^2}{\partial {s^*}^2} 
\Bigg[  
\frac{\partial^2 v^*}{\partial {s^*}^2} 
\!+\! \frac{\partial^2 v^*}{\partial {s^*}^2} (\frac{\partial v^*}{\partial {s^*}} )^2 
\!+\! \frac{E_r^*}{2} \prescript{RL}{0}{\mathcal D}_{t^*}^{\alpha} \frac{\partial^2 v^*}{\partial {s^*}^2}(\frac{\partial v^*}{\partial {s^*}})^2 
\!+\! E_r^* \prescript{RL}{0}{\mathcal D}_{t^*}^{\alpha} \frac{\partial^2 v^*}{\partial {s^*}^2}
\!+\! \frac{1}{2}  E_r^* (\frac{\partial v^*}{\partial s^*})^2 \prescript{RL}{0}{\mathcal D}_{t^*}^{\alpha} \frac{\partial^2 v^*}{\partial {s^*}^2}
\Bigg]
 \nonumber\\
&
\!-\!
\frac{\partial}{\partial {s^*}} 
\Bigg[
\frac{\partial v^*}{\partial {s^*}} (\frac{\partial^2 v^*}{\partial {s^*}^2})^2 
\!+\! E_r^* \frac{\partial v^*}{\partial s^*} \frac{\partial^2 v^*}{\partial {s^*}^2}  \prescript{RL}{0}{\mathcal D}_{t^*}^{\alpha} \frac{\partial^2 v^*}{\partial {s^*}^2}
\Bigg]
= -\frac{\partial^2 {{v_b}}^*}{\partial {t^*}^2}  ,
\end{align}
The dimensionless boundary conditions are also obtained by substituting dimensionless variables in \eqref{Eq: bc}. We can show similarly that they preserve their structure as: 
\begin{align*}
& v^{*} \, \Big|_{s^{*}=0} = \frac{\partial v^{*}}{\partial s^{*}} \, \Big|_{s^{*}=0} = 0 ,
 \nonumber\\
& 
\frac{J^* \rho L^3 m}{\rho} \frac{1}{mL^4}
\Bigg[ 
\frac{\partial^3 v^{*}}{\partial t^{*} \partial^2 s^{*}} \left( 1 \!+\! \left(\frac{\partial v^{*}}{\partial s^{*}}\right)^2 \right) 
\!+\!\frac{\partial v^{*}}{\partial s^{*}} \left(\frac{\partial^2 v^{*}}{\partial t^{*} \partial s^{*}}\right)^2  
\Bigg]
\!+\!
\frac{1}{L} \Bigg[  
\frac{\partial^2 v^*}{\partial {s^*}^2} 
\!+\! \frac{\partial^2 v^*}{\partial {s^*}^2} (\frac{\partial v^*}{\partial {s^*}} )^2 
\!+\! \frac{E_r^*}{2} \prescript{RL}{0}{\mathcal D}_{t^*}^{\alpha} \frac{\partial^2 v^*}{\partial {s^*}^2}(\frac{\partial v^*}{\partial {s^*}})^2 
 \nonumber\\
&
\!+\! E_r^* \prescript{RL}{0}{\mathcal D}_{t^*}^{\alpha} \frac{\partial^2 v^*}{\partial {s^*}^2}
\!\!+ \frac{1}{2}  E_r^* (\frac{\partial v^*}{\partial s^*})^2 \prescript{RL}{0}{\mathcal D}_{t^*}^{\alpha} \frac{\partial^2 v^*}{\partial {s^*}^2}
\Bigg]
\Bigg|_{s^{*}=1} \!=\! 0 ,
 \nonumber\\
& 
\frac{M^* \rho L m}{\rho} \frac{L}{mL^4} \left( \frac{\partial^2 v^{*}}{\partial^2 t^{*}}  + \frac{\partial^2 {{v_b}}^{*}}{\partial^2 t^{*}}  \right)
\!-\!\frac{1}{L^2} \frac{\partial v^{*}}{\partial s^{*}} \Bigg[  
\frac{\partial^2 v^*}{\partial {s^*}^2} 
\!+\! \frac{\partial^2 v^*}{\partial {s^*}^2} (\frac{\partial v^*}{\partial {s^*}} )^2 
\!+\! \frac{E_r^*}{2} \prescript{RL}{0}{\mathcal D}_{t^*}^{\alpha} \frac{\partial^2 v^*}{\partial {s^*}^2}(\frac{\partial v^*}{\partial {s^*}})^2 
\!+\! E_r^* \prescript{RL}{0}{\mathcal D}_{t^*}^{\alpha} \frac{\partial^2 v^*}{\partial {s^*}^2}
 \nonumber\\
& 
\!+\! \frac{1}{2}  E_r^* (\frac{\partial v^*}{\partial s^*})^2 \prescript{RL}{0}{\mathcal D}_{t^*}^{\alpha} \frac{\partial^2 v^*}{\partial {s^*}^2}
\Bigg]
\!+\!\frac{1}{L^2} 
\Bigg[
\frac{\partial v^*}{\partial {s^*}} (\frac{\partial^2 v^*}{\partial {s^*}^2})^2 
\!+\! E_r^* \frac{\partial v^*}{\partial s^*} \frac{\partial^2 v^*}{\partial {s^*}^2}  \prescript{RL}{0}{\mathcal D}_{t^*}^{\alpha} \frac{\partial^2 v^*}{\partial {s^*}^2}
\Bigg] \Bigg|_{s^{*}=1}  \!=\! 0 ,
\end{align*}
Therefore, the dimensionless equation of motion becomes (after dropping $^*$ for the sake of simplicity)
\begin{align}
\label{Eq: eqn of motion - dimless - 2}
&\frac{\partial^2{v}}{\partial{t}^2} 
 \!+\! 
\frac{\partial^2}{\partial{s}^2} \left(  
\frac{\partial^2{v}}{\partial{s}^2} \!+\! \frac{\partial^2{v}}{\partial{s}^2} (\frac{\partial{v}}{\partial{s}})^2
\!+\! E_r \prescript{RL}{0}{\mathcal D}_{t}^{\alpha} \, \Big[\frac{\partial^2{v}}{\partial{s}^2}(1 + \frac{1}{2} (\frac{\partial{v}}{\partial{s}})^2)\Big] 
\!+\! \frac{1}{2}  E_r (\frac{\partial{v}}{\partial{s}})^2 \, \prescript{RL}{0}{\mathcal D}_{t}^{\alpha} \, \frac{\partial^2{v}}{\partial{s}^2}
\right)
\!-\! 
\frac{\partial}{\partial{s}} \left( 
\frac{\partial{v}}{\partial{s}}  \,  (\frac{\partial^2{v}}{\partial{s}^2})^2  
\!+\! E_r \frac{\partial{v}}{\partial{s}}  \, \frac{\partial^2{v}}{\partial{s}^2} \,  \prescript{RL}{0}{\mathcal D}_{t}^{\alpha} \, \frac{\partial^2{v}}{\partial{s}^2}
\right)
\!=\! \!-\!\ddot{{{v_b}}} ,
\end{align}
which is subject to the following dimensionless boundary conditions
\begin{align}
\label{Eq: bc - dimless}
& v \, \Big|_{s=0} = \frac{\partial{v}}{\partial{s}} \, \Big|_{s=0} = 0 ,
\nonumber\\
& 
J \left(  \frac{\partial^3{v}}{\partial{t}^2\partial{s}}  ( 1 + (\frac{\partial{v}}{\partial{s}})^2) + \frac{\partial{v}}{\partial{s}}  (\frac{\partial^2{v}}{\partial{s}\partial{t}})^2   \right)
+
\left( 
\frac{\partial^2{v}}{\partial{s}^2} + \frac{\partial^2{v}}{\partial{s}^2} (\frac{\partial{v}}{\partial{s}})^2
+ E_r \, \prescript{RL}{0}{\mathcal D}_{t}^{\alpha} \, \Big[\frac{\partial^2{v}}{\partial{s}^2}(1 + \frac{1}{2} (\frac{\partial{v}}{\partial{s}})^2)\Big] 
+ \frac{1}{2}  E_r \, (\frac{\partial{v}}{\partial{s}})^2 \, \prescript{RL}{0}{\mathcal D}_{t}^{\alpha} \, \frac{\partial^2{v}}{\partial{s}^2}  
\right) \,\, \Bigg|_{s=1} = 0 ,
 \nonumber\\
& 
M (\frac{\partial^2{v}}{\partial{t}^2} \!+\! \ddot{{{v_b}}})
\!-\! \frac{\partial}{\partial{s}} \left(
\frac{\partial^2{v}}{\partial{s}^2} \!+\! \frac{\partial^2{v}}{\partial{s}^2} (\frac{\partial{v}}{\partial{s}})^2
\!+\! E_r \prescript{RL}{0}{\mathcal D}_{t}^{\alpha} \, \Big[\frac{\partial^2{v}}{\partial{s}^2}(1 \!+\! \frac{1}{2} (\frac{\partial{v}}{\partial{s}})^2) \Big]
\!+\! \frac{1}{2}  E_r (\frac{\partial{v}}{\partial{s}})^2 \prescript{RL}{0}{\mathcal D}_{t}^{\alpha} \frac{\partial^2{v}}{\partial{s}^2}  
\right)
\!+\! \left( 
\frac{\partial{v}}{\partial{s}} \, (\frac{\partial^2{v}}{\partial{s}^2})^2  \!+\! E_r \, \frac{\partial{v}}{\partial{s}} \, \frac{\partial^2{v}}{\partial{s}^2}\prescript{RL}{0}{\mathcal D}_{t}^{\alpha} \frac{\partial^2{v}}{\partial{s}^2}	
\right) \Bigg|_{s=1} \!=\! 0.
\end{align}
%
\section{Single Mode Decomposition}
\label{Sec: App. Single Mode Decomposition} 
%
{In order to demonstrate that} the single mode decomposition satisfies the weak form solution and boundary conditions, we consider the case of no lumped mass at {the tip, \textit{i.e.},} $M=J=0$, and check if the proposed approximate solution {solves the weak form and corresponding} boundary conditions. First we substitute the boundary conditions in the weak formulation and then we use {single-mode} approximation to recover equation \eqref{Eq: weak form - discrete 2}. {Therefore, start integrating Equation \eqref{Eq: weak form - 1} by parts as follows:}
\begin{align}
\label{Eq: weak form single mode - 1}
& \int_{0}^{1} \frac{\partial^2{v}}{\partial{t}^2}  \tilde{v}  ds  
\!+\! 
\frac{\partial}{\partial{s}}\left(  
\frac{\partial^2{v}}{\partial{s}^2} \!+\! \frac{\partial^2{v}}{\partial{s}^2} (\frac{\partial{v}}{\partial{s}})^2 
\!+\! E_r \prescript{RL}{0}{\mathcal D}_{t}^{\alpha}  \Big[\frac{\partial^2{v}}{\partial{s}^2}(1 \!+\! \frac{1}{2} (\frac{\partial{v}}{\partial{s}})^2) \Big] \!+\! \frac{1}{2}  E_r  (\frac{\partial{v}}{\partial{s}})^2  \prescript{RL}{0}{\mathcal D}_{t}^{\alpha}  \frac{\partial^2{v}}{\partial{s}^2} \right) \tilde{v}  \Bigg|_{0}^{1} \nonumber \\
& \!-\! \left(  
\frac{\partial^2{v}}{\partial{s}^2} \!+\! \frac{\partial^2{v}}{\partial{s}^2} (\frac{\partial{v}}{\partial{s}})^2 
\!+\! E_r \prescript{RL}{0}{\mathcal D}_{t}^{\alpha}  \Big[\frac{\partial^2{v}}{\partial{s}^2}(1 \!+\! \frac{1}{2} (\frac{\partial{v}}{\partial{s}})^2) \Big]
+ \frac{1}{2}  E_r 
(\frac{\partial{v}}{\partial{s}})^2  \prescript{RL}{0}{\mathcal D}_{t}^{\alpha}  \frac{\partial^2{v}}{\partial{s}^2}
\right)
\frac{\partial\tilde{v}}{\partial{s}}  \Bigg|_{0}^{1}
+ \int_{0}^{1} \left(  
\frac{\partial^2{v}}{\partial{s}^2} + \frac{\partial^2{v}}{\partial{s}^2} (\frac{\partial{v}}{\partial{s}})^2 
+ E_r
\prescript{RL}{0}{\mathcal D}_{t}^{\alpha}  \Big[\frac{\partial^2{v}}{\partial{s}^2}(1 \!+\! \frac{1}{2} (\frac{\partial{v}}{\partial{s}})^2) \Big]
\right. \nonumber\\
&\left.
\!+\! \frac{1}{2}  E_r 
(\frac{\partial{v}}{\partial{s}})^2  \prescript{RL}{0}{\mathcal D}_{t}^{\alpha}  \frac{\partial^2{v}}{\partial{s}^2}
\right)
\frac{\partial^2{\tilde{v}}}{\partial{s}^2}  ds
- \left( 
\frac{\partial{v}}{\partial{s}}   (\frac{\partial^2{v}}{\partial{s}^2})^2  
+ E_r
\frac{\partial{v}}{\partial{s}}  \frac{\partial^2{v}}{\partial{s}^2}   \prescript{RL}{0}{\mathcal D}_{t}^{\alpha}  \frac{\partial^2{v}}{\partial{s}^2}
\right)
\tilde{v}  \Bigg|_{0}^{1}
+\int_{0}^{1}\left( 
\frac{\partial{v}}{\partial{s}}   (\frac{\partial^2{v}}{\partial{s}^2})^2  
+ E_r
\frac{\partial{v}}{\partial{s}}  \frac{\partial^2{v}}{\partial{s}^2}   \prescript{RL}{0}{\mathcal D}_{t}^{\alpha}  \frac{\partial^2{v}}{\partial{s}^2}
\right)
\frac{\partial{\tilde{v}}}{\partial{s}}  ds
= f(t).
\end{align}
When $M=J=0$ {the boundary conditions in equation \eqref{Eq: bc - dimless - body} are given by:}
\begin{align}
\label{Eq: weak form single mode - 2}
& v \, \Big|_{s=0} = \frac{\partial{v}}{\partial{s}} \, \Big|_{s=0} = 0 ,
\nonumber\\
& 
\left( 
\frac{\partial^2{v}}{\partial{s}^2}+  \frac{\partial^2{v}}{\partial{s}^2} ( \frac{\partial{v}}{\partial{s}})^2
+ E_r \, \prescript{RL}{0}{\mathcal D}_{t}^{\alpha} \,  \Big[\frac{\partial^2{v}}{\partial{s}^2}(1 + \frac{1}{2}  (\frac{\partial{v}}{\partial{s}})^2) \Big]
+ \frac{1}{2}  E_r \,  (\frac{\partial{v}}{\partial{s}})^2 \, \prescript{RL}{0}{\mathcal D}_{t}^{\alpha} \,  \frac{\partial^2{v}}{\partial^2{s}}
\right) \,\, \Bigg|_{s=1} = 0 ,
\nonumber\\
& 
- \frac{\partial{v}}{\partial{s}}\left(
\frac{\partial^2{v}}{\partial{s}^2}+ \frac{\partial^2{v}}{\partial{s}^2} (\frac{\partial{v}}{\partial{s}})^2
+ E_r \, \prescript{RL}{0}{\mathcal D}_{t}^{\alpha} \, \Big[\frac{\partial^2{v}}{\partial{s}^2}(1 + \frac{1}{2} (\frac{\partial{v}}{\partial{s}})^2)\Big]
+ \frac{1}{2}  E_r \, (\frac{\partial{v}}{\partial{s}})^2 \, \prescript{RL}{0}{\mathcal D}_{t}^{\alpha} \, \frac{\partial^2{v}}{\partial{s}^2}
\right)
+ \left( 
\frac{\partial{v}}{\partial{s}} \, (\frac{\partial^2{v}}{\partial{s}^2})^2  + E_r \, \frac{\partial{v}}{\partial{s}} \, \frac{\partial^2{v}}{\partial{s}^2}  \, \prescript{RL}{0}{\mathcal D}_{t}^{\alpha} \, \frac{\partial^2{v}}{\partial{s}^2} 
\right) \,\, \Bigg|_{s=1}  = 0 ,
\end{align}
By substituting \eqref{Eq: weak form single mode - 2} in \eqref{Eq: weak form single mode - 1} we {obtain:}
\begin{align}
\label{Eq: weak form single mode - 3}
&
\int_{0}^{1} \frac{\partial^2{v}}{\partial{t}^2} \, \tilde{v} \, ds  
+ 
\
\int_{0}^{1} \left(  
\frac{\partial^2{v}}{\partial{s}^2}+ \frac{\partial^2{v}}{\partial{s}^2} (\frac{\partial{v}}{\partial{s}})^2 
+ E_r
\prescript{RL}{0}{\mathcal D}_{t}^{\alpha} \, \Big[\frac{\partial^2{v}}{\partial{s}^2}(1 + \frac{1}{2} (\frac{\partial{v}}{\partial{s}})^2)\Big] \right.
\left. + \frac{1}{2}  E_r 
(\frac{\partial{v}}{\partial{s}})^2 \, \prescript{RL}{0}{\mathcal D}_{t}^{\alpha} \, \frac{\partial^2{v}}{\partial{s}^2}
\right)
\, \frac{\partial^2{\tilde{v}}}{\partial{s}^2} \, ds
\nonumber\\
& 
+\int_{0}^{1}\left( 
\frac{\partial{v}}{\partial{s}} \,  (\frac{\partial^2{v}}{\partial{s}^2})^2  
+ E_r
\frac{\partial{v}}{\partial{s}} \, \frac{\partial^2{v}}{\partial{s}^2} \,  \prescript{RL}{0}{\mathcal D}_{t}^{\alpha} \, \frac{\partial^2{v}}{\partial{s}^2}
\right)
\, \frac{\partial^2{\tilde{v}}}{\partial{s}^2}\tilde{v}^{\prime} \, ds
= f(t) ,
\end{align}
The modal discretization {utilized} in \eqref{Eq: assumed mode} can be simplified as  $v(s,t)=q(t)\phi(s)$, where we choose $\phi(0)=\phi^{\prime}(0)=0$, and $\phi^{\prime\prime}(0)=\phi^{\prime\prime\prime}(0)=0$. {Furthermore, setting} $ \tilde{v}=\phi(s)$, \eqref{Eq: weak form single mode - 3} is {given by:}
\begin{align}
\label{Eq: weak form single mode - 4}
&
\ddot{q}\int_{0}^{1} \phi \phi  ds  
\!+\! \int_{0}^{1} \left(  
q \phi^{\prime\prime} \!+\! q^{3}\phi^{\prime\prime}
{\phi^{\prime}}^2 \right)\phi^{\prime\prime}  ds
\!+\! \int_{0}^{1} E_r \phi^{\prime\prime} \prescript{RL}{0}{\mathcal D}_{t}^{\alpha}  q\phi^{\prime\prime} ds
\!+\!\int_{0}^{1}
\frac{1}{2}  E_r 
\phi^{\prime\prime}(\phi^{\prime})^2  \prescript{RL}{0}{\mathcal D}_{t}^{\alpha}  q^3 \phi^{\prime\prime}
ds
 \nonumber\\
&
\!+\! \int_{0}^{1}
\frac{1}{2}  E_r 
{\phi^{\prime}}^2 \phi^{\prime\prime} q^{2}  \prescript{RL}{0}{\mathcal D}_{t}^{\alpha}  q \phi^{\prime\prime}
ds
\!+\!\int_{0}^{1}\left( \phi^{\prime}   {\phi^{\prime\prime}}^2q^3 \!+\!  E_r
\phi^{\prime}\phi^{\prime\prime}\phi^{\prime\prime} q^2    \prescript{RL}{0}{\mathcal D}_{t}^{\alpha}  q 
\right)\phi^{\prime}
ds
\!=\! f(t) ,
\end{align}
considering \eqref{Eq: coeff unimodal} for the case {without a lumped mass, we have,}
\begin{align}
\label{Eq: weak form single mode - 5}
& \mathcal{M} = \int_{0}^{1} \phi^2 \, ds, \,  \quad\mathcal{K}_l = \mathcal{C}_l = \int_{0}^{1}  {\phi^{\prime\prime}}^2 \,\, ds,\, \quad
\mathcal{K}_{nl} = \mathcal{C}_{nl} = \int_{0}^{1}  {\phi^{\prime}}^2 \, {\phi^{\prime\prime}}^2 \,\, ds,
\quad  \mathcal{M}_b = \int_{0}^{1} \phi \, ds  .
\nonumber\\
\end{align}
{Substituting the above relationships into \eqref{Eq: weak form single mode - 4}, we obtain} \eqref{Eq: weak form - discrete 2}:
\begin{align}
\label{Eq: weak form single mode - 6}
&\mathcal{M} \, \ddot{q} 
+ \mathcal{K}_l \, q + E_r \, \mathcal{C}_l \, \prescript{RL}{0}{\mathcal D}_{t}^{\alpha} q  
+ 2 \mathcal{K}_{nl} \, q^3 + \frac{E_r \, \mathcal{C}_{nl}}{2} \, \left( \prescript{RL}{0}{\mathcal D}_{t}^{\alpha} q^3 + 3 \, q^2 \, \prescript{RL}{0}{\mathcal D}_{t}^{\alpha} q  \right)
= -\mathcal{M}_b \, \ddot{{{v_b}}}.
\end{align}
\bigskip
\section{Deriving the Linearized Equation {of Motion}}
\label{Sec: App. Linearization}

{Since our problem considers geometrical nonlinearities, we let the following kinematic linearizations under the assumption small motions, respective, for the rotation angle (\ref{Eq: rotation}), angular velocity (\ref{Eq: angular velocity}), and curvature (\ref{Eq: curvature}):}
\begin{equation}\label{eq:linear_approximations}
\psi \simeq \frac {\partial{v}}{\partial{s}}, \qquad \frac {\partial{\psi}}{\partial{t}} \simeq \frac {\partial^2{v}}{\partial{t}\partial{s}},\qquad \frac {\partial{\psi}}{\partial{s}} \simeq \frac {\partial^2{v}}{\partial{s}^2},
\end{equation}
{With their corresponding variations given by:}
{
\begin{equation}\label{eq:linear_variations}
\delta \psi \simeq \delta \frac {\partial{v}}{\partial{s}}, \qquad \delta \frac {\partial{\psi}}{\partial{t}} \simeq \delta \frac {\partial^2{v}}{\partial{t}\partial{s}},\qquad \delta \frac{\partial{\psi}}{\partial{s}} \simeq \delta \frac {\partial^2{v}}{\partial{s}^2}
\end{equation}}
{Similar to \ref{Sec: App. extended Hamilton}, the variation of total work is expressed as:}
\begin{align}
\label{Eq: total work var linearized}
\delta w 
= &
\int_{0}^{L} \int_{A} \left( -\eta \, E_{\infty} \, \frac{\partial{\psi}}{\partial{s}} -\eta \, E_{\alpha} \, \prescript{RL}{0}{\mathcal D}_{t}^{\alpha} \, \frac{\partial{\psi}}{\partial{s}} \right) \, (-\eta \, \delta \frac{\partial{\psi}}{\partial{s}}) \, dA  \,\, ds
 \nonumber\\
 = &
\int_{0}^{L} \left( E_{\infty} \, \left( \int_{A}  \, \eta^2 dA \right) \, \frac{\partial{\psi}}{\partial{s}} 
+ E_{\alpha} \, \left( \int_{A}  \, \eta^2 dA \right) \, \prescript{RL}{0}{\mathcal D}_{t}^{\alpha} \, \frac{\partial{\psi}}{\partial{s}} \right) \, \delta \frac{\partial{\psi}}{\partial{s}} \,\, ds
\nonumber\\
= &
\int_{0}^{L} \left( E_{\infty} \, I \, \frac{\partial{\psi}}{\partial{s}} + E_{\alpha} \, I \, \prescript{RL}{0}{\mathcal D}_{t}^{\alpha} \, \frac{\partial{\psi}}{\partial{s}} \right) \, \delta \frac{\partial{\psi}}{\partial{s}} \,\, ds.
\end{align}
{Employing approximation \eqref{eq:linear_variations} for the variation of curvature, the variation of total energy becomes:}
\begin{align}
\label{Eq: Total Energy Linearized}
\delta w 
=
\int_{0}^{L}  
\left( E_{\infty} \, I \, \frac{\partial^2{v}}{\partial{s}^2}  + E_{\alpha} \, I \, \prescript{RL}{0}{\mathcal D}_{t}^{\alpha} \, \Big[\frac{\partial^2{v}}{\partial{s}^2}\Big]\right)  \, \delta \frac{\partial^2{v}}{\partial{s}^2} \,\, ds.
\end{align}
{Expanding} the terms and integrating by parts, we have
\begin{align}
\label{Eq: Total Energy Linearized Expantion}
\delta w 
\!=\!
\int_{0}^{L}  
\frac{\partial^2}{\partial{s}^2}
\left( E_{\infty} \, I \, \frac{\partial^2{v}}{\partial{s}^2} \!+\! E_{\alpha} \, I \, \prescript{RL}{0}{\mathcal D}_{t}^{\alpha} \Big[\frac{\partial^2{v}}{\partial{s}^2}\Big]\right)  
\delta v \, ds 
\!+\!
\left( E_{\infty} \, I \, \frac{\partial^2{v}}{\partial{s}^2} \!+\! E_{\alpha} \, I \, \prescript{RL}{0}{\mathcal D}_{t}^{\alpha} \Big[\frac{\partial^2{v}}{\partial{s}^2}\Big]\right) 
 \delta \frac{\partial{v}}{\partial{s}} \Bigg|_{0}^{L}
\!-\!
\frac{\partial}{\partial{s}}
\left( E_{\infty} \, I \, \frac{\partial^2{v}}{\partial{s}^2} \!+\! E_{\alpha} \, I \, \prescript{RL}{0}{\mathcal D}_{t}^{\alpha} \, \Big[\frac{\partial^2{v}}{\partial{s}^2}\Big]\right)\, \delta v \Bigg|_{0}^{L}
\end{align}
{Employing the boundary conditions} $\delta v(0,t) = \delta \frac{\partial{v}}{\partial{s}}(0,t) = 0$ {into \eqref{Eq: Total Energy Linearized Expantion}, we obtain:}
\begin{align}
\label{Eq: total work var 4}
\delta w 
& \!=\!
\int_{0}^{L}  
\frac{\partial^2}{\partial{s}^2}
\left( E_{\infty} I \frac{\partial^2{v}}{\partial{s}^2} \!+\! E_{\alpha} I \prescript{RL}{0}{\mathcal D}_{t}^{\alpha}  \Big[\frac{\partial^2{v}}{\partial{s}^2}\Big]\right)
\delta v ds 
\!+\!
\left( E_{\infty} I \frac{\partial^2{v}}{\partial{s}^2} \!+\! E_{\alpha} I \prescript{RL}{0}{\mathcal D}_{t}^{\alpha}  \Big[\frac{\partial^2{v}}{\partial{s}^2}\Big]\right)\Bigg|_{s=L} \delta \frac{\partial{v}}{\partial{s}}(L,t) 
\nonumber \\
& \!-\!
\frac{\partial}{\partial{s}}\left(
\left( E_{\infty} I \frac{\partial^2{v}}{\partial{s}^2} \!+\! E_{\alpha} I \prescript{RL}{0}{\mathcal D}_{t}^{\alpha}  \Big[\frac{\partial^2{v}}{\partial{s}^2}\Big]\right)\right) \Bigg|_{s=L} \delta v(L,t).
\end{align}
{From \ref{Sec: App. extended Hamilton}, the kinetic energy of the beam is given to:}
\begin{align}
\label{Eq: Kinetic Linearized - 2}
T =& 
\frac{1}{2} \, \rho \int_{0}^{L} ( \frac{\partial{v}}{\partial{t}} + \dot {{v_b}})^2  \,   ds 
+ \frac{1}{2}M (\frac{\partial{v}}{\partial{t}} + \dot {{v_b}})^2  \Big|_{s=L}
+ \frac{1}{2} J (\frac{\partial{\psi}}{\partial{t}})^2 \Big|_{s=L}, 
\end{align}
where its variation can be taken as
\begin{align}
\label{Eq: Kinetic Variation Linearized}
\delta T 
= &\rho \int_{0}^{L} (\frac{\partial{v}}{\partial{t}} + \dot{{{v_b}}}) \, \delta \frac{\partial{v}}{\partial{t}} \, ds + M (\frac{\partial{v}}{\partial{t}} + \dot{{{v_b}}}) \, \delta \frac{\partial{v}}{\partial{t}} \Big|_{s=L} 
+  J \frac{\partial{\psi}}{\partial{t}} \, \delta \frac{\partial{\psi}}{\partial{t}} \Big|_{s=L} ,
\end{align}
{Employing the approximations (\ref{eq:linear_approximations}) and (\ref{eq:linear_variations}), to the above equation, we obtain:}
{
\begin{align}
\label{Eq: Kinetic Variation Linearized - 2}
\delta T 
 \simeq\,&\rho \int_{0}^{L} (\frac{\partial{v}}{\partial{t}} \!+\! \dot{{{v_b}}})  \delta \frac{\partial{v}}{\partial{t}}  ds 
\!+\! M (\frac{\partial{v}}{\partial{t}} \!+\! \dot{{{v_b}}})  \delta \frac{\partial{v}}{\partial{t}} \Big|_{s=L} 
\!+\! J  \frac{\partial^2{v}}{\partial{t}\partial{s}} \delta \frac{\partial^2{v}}{\partial{t}\partial{s}} \Big|_{s=L} .
\end{align}}
The time integration of $\delta T$ takes the following form through integration by parts
{
\begin{align}
\label{Eq: Kinetic Var Linearized- 2}
\int_{t_1}^{t_2} \delta T \, dt 
= & \int_{t_1}^{t_2} 
\Bigg\{
\rho \int_{0}^{L} \left(\frac{\partial{v}}{\partial{t}} + \dot{{{v_b}}}\right) \, \delta \frac{\partial{v}}{\partial{t}} \, ds 
+ M \left(\frac{\partial{v}}{\partial{t}} + \dot{{{v_b}}}\right) \, \delta \frac{\partial{v}}{\partial{t}} \Big|_{s=L} 
+ J \left(\frac{\partial^2{v}}{\partial{t}\partial{s}}\right) \delta \frac{\partial^2{v}}{\partial{t}\partial{s}} \Big|_{s=L}
\Bigg\} \, dt
 \nonumber\\
= &
\int_{t_1}^{t_2} \rho \int_{0}^{L} \left(\frac{\partial{v}}{\partial{t}} \!+\! \dot{{{v_b}}}\right) \, \delta \frac{\partial{v}}{\partial{t}} \, ds \, dt
\!+\! M \int_{t_1}^{t_2} \left(\frac{\partial{v}}{\partial{t}} \!+\! \dot{{{v_b}}}\right) \, \delta \frac{\partial{v}}{\partial{t}} \Big|_{s=L} \, dt
\!+\! J \int_{t_1}^{t_2} \left(\frac{\partial^2{v}}{\partial{t}\partial{s}}\right) \delta \frac{\partial^2{v}}{\partial{t}\partial{s}}\Big|_{s=L} \, dt
 \nonumber\\
= &
\rho \int_{0}^{L} \int_{t_1}^{t_2}  \left(\frac{\partial{v}}{\partial{t}} \!+\! \dot{{{v_b}}}\right) \, \delta \frac{\partial{v}}{\partial{t}} \, dt \, ds
\!+\! M \int_{t_1}^{t_2} \left(\frac{\partial{v}}{\partial{t}} \!+\! \dot{{{v_b}}}\right) \, \delta \frac{\partial{v}}{\partial{t}}  \, dt \, \Big|_{s=L}
\!+\! J \int_{t_1}^{t_2} \left(\frac{\partial^2{v}}{\partial{t}\partial{s}}\right) \delta \frac{\partial^2{v}}{\partial{t}\partial{s}}  \, dt \, \Big|_{s=L}
\nonumber\\
= &
\rho \int_{0}^{L} \left[ \left(\frac{\partial{v}}{\partial{t}} \!+\! \dot{{{v_b}}}\right)  \delta{v} \Big|_{t_1}^{t_2} \!-\! \int_{t_1}^{t_2}  \left(\frac{\partial^2{v}}{\partial{t}^2} \!+\! \ddot{{{v_b}}}\right)  \delta {v}  dt \right] ds
\!+\! M \left(\frac{\partial{v}}{\partial{t}} \!+\! \dot{{{v_b}}}\right)  \delta v \Big|_{s=L} \Big|_{t_1}^{t_2} 
\!-\! M \int_{t_1}^{t_2} \left(\frac{\partial^2{v}}{\partial{t}^2} \!+\! \ddot{{{v_b}}}\right)  \delta v   dt  \Big|_{s=L}
\nonumber\\
& 
+ J  \left(\frac{\partial^2{v}}{\partial{t}\partial{s}}\right) \delta \frac{\partial{v}}{\partial{s}}  \Big|_{s=L}  \Big|_{t_1}^{t_2}
- J \int_{t_1}^{t_2} 
\left(  \frac{\partial^3{v}}{\partial{t}^2\partial{s}} \right) \delta \frac{\partial{v}}{\partial{s}}   dt  \Big|_{s=L}
 \nonumber\\
= &
- \!\int_{t_1}^{t_2}\! \Bigg\{
\rho \!\int_{0}^{L}\!  (\frac{\partial^2{v}}{\partial{t}^2} \!+\! \ddot{{{v_b}}})  \delta {v}  ds 
\!+\! M (\frac{\partial^2{v}}{\partial{t}^2} \!+\! \ddot{{{v_b}}})  \delta v   \Big|_{s=L}
\!+\! J \left(  \frac{\partial^3{v}}{\partial{t}^2\partial{s}} \right)  \delta \frac{\partial{v}}{\partial{s}}   \Big|_{s=L}
\Bigg\}  dt,
\end{align}}
where we consider that $\delta v = \delta \frac{\partial{v}}{\partial{s}} = 0$ at $t=t_1$ and $t=t_2$. Therefore, the extended Hamilton's principle takes the form
{
\begin{align}
\label{Eq: Extended Hamilton Linearized}
&
\int_{t_1}^{t_2} \Bigg\{
\int_{0}^{L} \Bigg[
- \rho (\frac{\partial^2{v}}{\partial{t}^2} + \ddot{{{v_b}}})
 -\frac{\partial^2}{\partial{s}^2}
\left( E_{\infty} \, I \, \frac{\partial^2{v}}{\partial{s}^2} + E_{\alpha} \, I \, \prescript{RL}{0}{\mathcal D}_{t}^{\alpha} \, \Big[\frac{\partial^2{v}}{\partial{s}^2}\Big]\right) \Bigg]\delta v \, ds
 - \left[ M (\frac{\partial^2{v}}{\partial{t}^2} \!+\! \ddot{{{v_b}}}) 
  \!-\! \frac{\partial}{\partial{s}}
\left( E_{\infty} \, I \, \frac{\partial^2{v}}{\partial{s}^2} + E_{\alpha} \, I \, \prescript{RL}{0}{\mathcal D}_{t}^{\alpha} \, \Big[\frac{\partial^2{v}}{\partial{s}^2}\Big]\right) \right] \delta v \Big|_{s=L} \nonumber \\
& - \left[  J \left(  \frac{\partial^3{v}}{\partial{t}^2\partial{s}} \right) -
\left( E_{\infty} \, I \, \frac{\partial^2{v}}{\partial{s}^2} + E_{\alpha} \, I \, \prescript{RL}{0}{\mathcal D}_{t}^{\alpha} \, \Big[\frac{\partial^2{v}}{\partial{s}^2}\Big]\right) \right] \delta \frac{\partial{v}}{\partial{s}}\Bigg|_{s=L} 
\,\, \Bigg\} \, dt = 0.
\end{align}}
Invoking the arbitrariness of virtual displacement $\delta v$, we obtain the strong form of the equation of motion as:
\begin{align}
\label{Eq: Eqn of Motion Linearized}
& \rho \, \frac{\partial^2{v}}{\partial{t}^2}  
+ 
E_{\infty} \, I \, \frac{\partial^2}{\partial{s}^2}\left(  \frac{\partial^2{v}}{\partial{s}^2}  \right)
+ 
E_{\alpha} \, I \, \frac{\partial^2}{\partial{s}^2}\left(   \, \prescript{RL}{0}{\mathcal D}_{t}^{\alpha} \, \Big[\frac{\partial^2{v}}{\partial{s}^2} \Big]\right)= { -\rho \ddot{v}_b }
\end{align}
which is subject to the following natural boundary conditions:
{ 
\begin{align}
\label{Eq: Natural BC Linearized}
& \Bigg\{ J   \frac{\partial^3{v}}{\partial{t}^2\partial{s}}  - E_{\infty} \, I \, \frac{\partial^2{v}}{\partial{s}^2}  - E_{\alpha} \, I \,  \, \prescript{RL}{0}{\mathcal D}_{t}^{\alpha} \, \Big[\frac{\partial^2{v}}{\partial{s}^2}\Big]\, \Bigg\}
\,\,  \Bigg|_{s=L} = 0 ,
\nonumber\\
& 
\Bigg\{ M \left(\frac{\partial^2{v}}{\partial{t}^2} + { \ddot{v}_b } \right)
- \frac{\partial}{\partial{s}}\left(
E_{\infty} \, I \, \frac{\partial^2{v}}{\partial{s}^2}   + E_{\alpha} \, I \, \, \prescript{RL}{0}{\mathcal D}_{t}^{\alpha} \, \Big[\frac{\partial^2{v}}{\partial{s}^2}\Big] \right) \Bigg\}  \, \Bigg|_{s=L}  = 0. 
\end{align}}
Therefore, the strong form reads as: find $v \in V$ such that
\begin{align}
\label{Eq: eqn of motion - 2 Linearized}
& m \, \frac{\partial^2{v}}{\partial{t}^2} 
+ 
\frac{\partial^2}{\partial{s}^2}\left(  
\frac{\partial^2{v}}{\partial{s}^2} 
+ E_r \prescript{RL}{0}{\mathcal D}_{t}^{\alpha} \, \Big[\frac{\partial^2{v}}{\partial{s}^2}\Big] 
\right)
= { -m \ddot{v}_b } ,
\end{align}
subject to the following boundary conditions:
\begin{align}
\label{Eq: bc Linearized}
& v \, \Big|_{s=0} = \frac{\partial{v}}{\partial{s}} \, \Big|_{s=0} = 0 ,
\\ \nonumber
& 
{ \frac{J m}{\rho}   \frac{\partial^3{v}}{\partial{t}^2\partial{s}} -} \left( 
\frac{\partial^2{v}}{\partial{s}^2}
+ E_r \, \prescript{RL}{0}{\mathcal D}_{t}^{\alpha} \, \Big[\frac{\partial^2{v}}{\partial{s}^2}\Big] 
\,
\right) \,\, \Bigg|_{s=L} = 0 ,
\\ \nonumber
& 
\frac{M m}{\rho} \left(\frac{\partial^2{v}}{\partial{t}^2} + { \ddot{v}_b }\right)
- \frac{\partial}{\partial{s}} \left(
\frac{\partial^2{v}}{\partial{s}^2} 
+ E_r \, \prescript{RL}{0}{\mathcal D}_{t}^{\alpha} \, \Big[\frac{\partial^2{v}}{\partial{s}^2} \Big]
\right)\,\, \Bigg|_{s=L}  = 0 ,
\end{align}
where $m = \frac{\rho}{ E_{\infty} \, I}$ and $E_r = \frac{E_{\alpha}}{E_{\infty}}$. 

\subsection{Nondimensionalization of Linearized Equation of Motion}

{Employing the dimensionless variables defined by (\ref{eq:dimenionsless_variables}) in a similar fashion as \ref{Sec: App. extended Hamilton}, and dropping the superscript ${}^*$ for simplicity, our dimensionless linearized equation of motion becomes:}
\begin{align}
\label{Eq: eqn of motion - dimless - 2 Linearized}
\frac{\partial^2{v}}{\partial{t}^2} 
& + 
\frac{\partial^2}{\partial{s}^2} \left(  
\frac{\partial^2{v}}{\partial{s}^2}
+ E_r \prescript{RL}{0}{\mathcal D}_{t}^{\alpha} \, \Big[\frac{\partial^2{v}}{\partial{s}^2}\Big] 
\right)
= { -\ddot{v}_b},
\end{align}
which is subject to the following dimensionless boundary conditions
\begin{align}
\label{Eq: bc - dimless Linearized}
& v \, \Big|_{s=0} = \frac{\partial{v}}{\partial{s}} \, \Big|_{s=0} = 0 ,
 \nonumber\\
& 
{ J   \frac{\partial^3{v}}{\partial{t}^2\partial{s}} -} \left( 
\frac{\partial^2{v}}{\partial{s}^2}
+ E_r \, \prescript{RL}{0}{\mathcal D}_{t}^{\alpha} \, \Big[\frac{\partial^2{v}}{\partial{s}^2}\Big]  
\right) \,\, \Bigg|_{s=1} = 0 ,
\nonumber\\
& 
{ M \left(\frac{\partial^2{v}}{\partial{t}^2} + \ddot{v}_b \right) }- \frac{\partial}{\partial{s}} \left(
\frac{\partial^2{v}}{\partial{s}^2}
+ E_r \, \prescript{RL}{0}{\mathcal D}_{t}^{\alpha} \, \Big[\frac{\partial^2{v}}{\partial{s}^2} \Big]
\right)
\,\, \Bigg|_{s=1}  = 0 ,
\end{align}
{To obtain the corresponding weak form, we multiply both sides of (\ref{Eq: eqn of motion - dimless - 2 Linearized}) by proper test functions $\tilde{v}(s) \in \tilde{V} $ and integrate the result over $\Omega_s = [0,1]$. Therefore:}
\begin{align}
\label{Eq: weak form - 1 Linearized}
& \int_{0}^{1} \frac{\partial^2{v}}{\partial{t}^2} \, \tilde{v} \, ds  
+ 
\int_{0}^{1} 
\frac{\partial^2}{\partial{s}^2}\left(  
\frac{\partial^2{v}}{\partial{s}^2}
+ E_r
\prescript{RL}{0}{\mathcal D}_{t}^{\alpha} \, \Big[\frac{\partial^2{v}}{\partial{s}^2}\Big] 
\right)
\, \tilde{v} \, ds
= { -\int_0^1 \frac{\partial^2 v_b }{\partial t^2}\tilde{v}\,ds}.
\end{align}
{Integrating be above expression by parts and employing the corresponding boundary conditions (\ref{Eq: bc - dimless Linearized}) with $M = J = 0$, we obtain:}
\begin{align}
\label{Eq: weak form - 3 linearized }
& \int_{0}^{1} { \frac{\partial^2 v}{\partial {t}^2}} \, \tilde{v} \, ds 
+
\int_{0}^{1} \frac{\partial^2{v}}{\partial{s}^2}\, \frac{\partial^2{\tilde{v}}}{\partial{s}^2} \, ds
+ E_r 
\int_{0}^{1} \prescript{RL}{0}{\mathcal D}_{t}^{\alpha} \, \Big[\frac{\partial^2{v}}{\partial{s}^2}\Big]\,\, \frac{\partial^2{\tilde{v}}}{\partial{s}^2} \, ds = { -\int_0^1 \frac{\partial^2 v_b }{\partial t^2}\tilde{v}\,ds}.
\end{align}
Using \eqref{Eq: assumed mode} and \eqref{Eq: Solution/Test Space}, the problem  \eqref{Eq: weak form - 3 linearized } reads: find $v_N \in V_N$ such that
\begin{align}
\label{Eq: weak form - discrete Linearized}
& \int_{0}^{1} { \frac{\partial^2 v_N}{\partial {t}^2}} \, \tilde{v}_N \, ds 
+ 
\int_{0}^{1} \frac{\partial^2{{v}_N}}{\partial{s}^2} \, \frac{\partial^2{\tilde{v}_N}}{\partial{s}^2} \, ds
+ E_r 
\int_{0}^{1} \prescript{RL}{0}{\mathcal D}_{t}^{\alpha} \, \Big[\frac{\partial^2{{v}_N}}{\partial{s}^2}\Big] \,\,  \frac{\partial^2{\tilde{v}_N}}{\partial{s}^2} \, ds,
= { -\int_0^1 \frac{\partial^2 v_b }{\partial t^2}\tilde{v}\,ds}.
\end{align}
for all $ \tilde{v}_N \in \tilde{V}_N$. {Substituting the single-mode approximation in \eqref{Eq: weak form - discrete Linearized} in a similar fashion as Section \ref{subsection}, we obtain the following unimodal governing equation of motion:}
\begin{align}
\label{Eq: weak form - discrete 2 Linearized}
\mathcal{M} \, \ddot{q}
+ \mathcal{K}_l \, q + E_r \, \mathcal{C}_l \, \prescript{RL}{0}{\mathcal D}_{t}^{\alpha} q  
= { -\mathcal{M}_b \ddot{v}_b},
\end{align}
{where the coefficients $\mathcal{M}$, $\mathcal{K}_l$, $\mathcal{C}_l$ and $\mathcal{M}_b$ are given by (\ref{Eq: weak form single mode - 5}). Finally, multiplying both sides of \eqref{Eq: weak form - discrete 2 Linearized} by $1/\mathcal{M}$, we obtain:}
\begin{align}
\label{Eq: weak form linear - discrete 2 Linearized}
\ddot{q} + E_r \, c_l \, \prescript{RL}{0}{\mathcal D}_{t}^{\alpha} q  + k_l \, q  = { -m_b \ddot{v}_b},
\end{align}
{with the coefficients $c_l = \frac{\mathcal{C}_l}{\mathcal{M}}$, $k_l = \frac{\mathcal{K}_l}{\mathcal{M}}$  and $m_n = \frac{\mathcal{M}_b}{\mathcal{M}}$. We note that although we previously assumed $M=J=0$ for the boundary conditions, we push the effects of the lumped mass with $M=J=1$ through the choice of our spatial eigenfunctions $\phi(s)$ \textit{(see \ref{Sec: App. Eigenvalue Problem of Linear Model})}.}
\section{Eigenvalue Problem of Linear Model}
\label{Sec: App. Eigenvalue Problem of Linear Model}
\begin{figure}[h]
	\centering
	\begin{subfigure}{0.24\textwidth}
		\centering
		\includegraphics[width=1\linewidth]{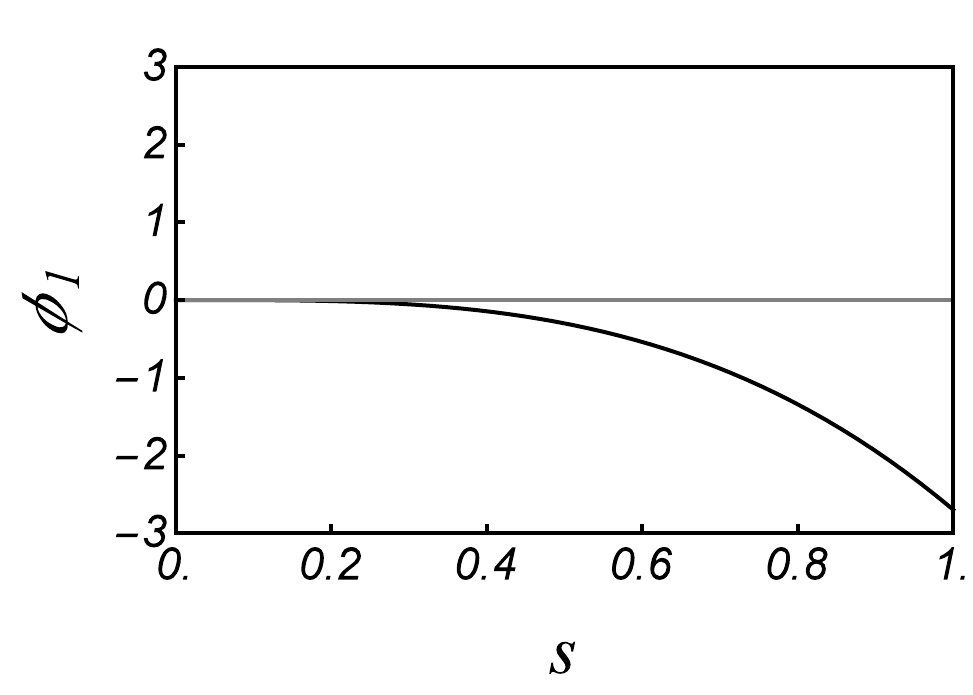}
	\end{subfigure}
	\begin{subfigure}{0.24\textwidth}
		\centering
		\includegraphics[width=1\linewidth]{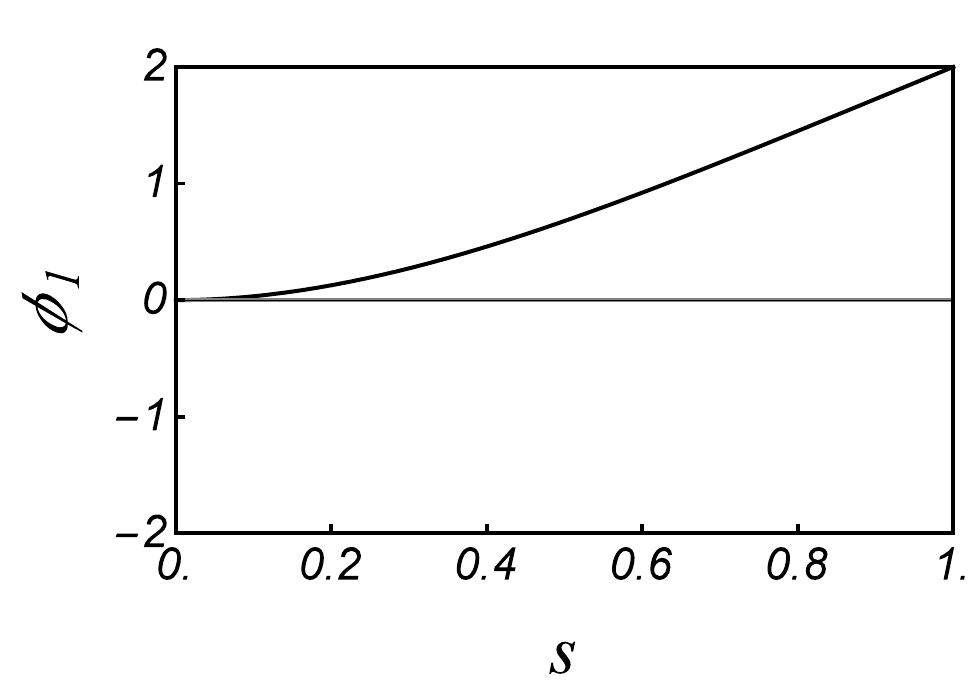}
	\end{subfigure}

	\caption{\footnotesize{Top: The first eigenfunctions, $X_1(s)$, of the undamped linear counterpart of our model. It is used as the spatial functions in the single mode approximation. Bottom: The first eigenfunctions, $X_1(s)$, of the undamped linear counterpart of our model with no lumped mass at the tip. It is used as the spatial functions in the single mode approximation.}}
	\label{Fig: Mode Shape Tip Mass}
\end{figure}

%
The assumed modes $\phi_i(s)$ in discretization \eqref{Eq: assumed mode} are obtained by solving the corresponding eigenvalue problem of free vibration of undamped linear counterparts to our model. Thus, the dimensionless linearized undamped equation of motion takes the form
\begin{align}
\label{Eq: Linear Euler Bernouli Beam}
& \frac{\partial^2}{\partial t^2}  v(s,t)
+ \frac{\partial^4}{\partial s^4} v(s,t)
= 0.
\end{align}
subject to linearized boundary conditions:
\begin{align}
\label{Eq: Linear BC}
& v(0,t) = 0 , &&  \frac{\partial^2{v}}{\partial{s}^2}(1,t) = -J \, \ddot{v}^{\prime}(1,t), 
 \nonumber\\
& v^{\prime}(0,t) = 0,
&& v^{\prime\prime\prime}(1,t) = M \, \ddot{v}(1,t) ,
\end{align}
where $ \dot{( \,\,\,)} = \frac{d }{dt}$ and $ ( \,\,\,)^{'} = \frac{d }{ds}$. We derive the corresponding eigenvalue problem by applying the separation of variables, i.e. $v(x,t) = X(s) T(t)$ to \eqref{Eq: Linear Euler Bernouli Beam}. Therefore, 
\begin{align}
\label{Eq: separ of vara}
\ddot{T}(t) X(s) + T(t) X^{''''}(s) = 0, \qquad 
\frac{\ddot{T}(t)}{T(t)}  + \frac{X^{''''}(s)}{X(s)} = 0, \qquad
\frac{\ddot{T}(t)}{T(t)} = -  \frac{X^{''''}(s)}{X(s)} = \lambda,
\end{align}
which gives the following equations
\begin{align}
\label{Eq: eigenProb 1}
&\ddot{T}(t) + \omega^2 T(t) = 0 , 
\\
\label{Eq: eigenProb 2}
&X^{\prime\prime\prime\prime}(s) - \beta^4 X(s) = 0, 
\end{align}
where $\beta^4 = \omega^2$ and the boundary conditions are
\begin{align*}
%
&X(0) = 0, && X^{\prime\prime}(1) =  J \, \omega^2 \, X^{\prime}(1), \\ \nonumber
&X^{\prime}(0) = 0,
&&X^{\prime\prime\prime}(1) = - M \, \omega^2 \, X(1).
\end{align*}
the solution to \eqref{Eq: eigenProb 2} is of the form $X(s) = A \sin(\beta s) + B \cos(\beta s) + C \sinh(\beta s) + D \cosh(\beta s)$, where $C = - A$ and $D = -B$, using the boundary conditions at $s=0$. Therefore, $$X(s) = A \left( \sin(\beta s) - \sinh(\beta s) \right) + B \left( \cos(\beta s) - \cosh(\beta s) \right).$$ Applying the first bondary condition at $s = 1$, i.e. $X^{\prime\prime}(1) =  J \, \omega^2 \, X^{\prime}(1)$ gives 
$$B = - \frac{ \sin(\beta) + \sinh(\beta) + J \beta^3 (\cos(\beta) - \cosh(\beta))}{ \cos(\beta) + \cosh(\beta)- J \beta^3 (\sin(\beta) - \sinh(\beta))} A, $$
that results in
\begin{align*}
&X(s) = A \left[ \left( \sin(\beta s) - \sinh(\beta s) \right) 
\!-\! \frac{ \sin(\beta) \!+\! \sinh(\beta) \!+\! J \beta^3 (\cos(\beta) \!-\! \cosh(\beta))}{ \cos(\beta) \!+\! \cosh(\beta)\!-\! J \beta^3 (\sin(\beta) \!-\! \sinh(\beta))} \left( \cos(\beta s) 
\!-\! \cosh(\beta s) \right) \right].
\end{align*}
Finally, using the second boundary condition at $s=1$ gives the {following transcendental equation for the case where $M=J=1$:}
\begin{align}
\label{Eq: freq}
& \!-\!\left(1 \!+\! \beta ^4 \!+\! \cos(\beta ) \cosh(\beta ) \right)
\!+\!\beta \left( \sin (\beta ) \cosh (\beta)\!-\!\cos (\beta ) \sinh (\beta ) \right) 
\!+\!\beta ^3 \left( \sin (\beta ) \cosh (\beta )\!-\!\sinh (\beta ) \cosh (\beta ) \right) \nonumber \\
& \!+\! \beta ^4 \left( \sin (\beta ) \sinh (\beta ) \!+\! \cos (\beta ) \cosh (\beta ) \right) 
\!=\! 0.
\end{align}
The first eigenvalue is computed as $\beta_1^2 = \omega_1 = 1.38569$, which results to the following first normalized eigenfunction, given in Fig. \ref{Fig: Mode Shape Tip Mass} (left). 
\begin{equation}\label{eq:eig_mass}
{ \phi(s) = 	X_1(s) = 5.50054 \sin (\beta_1 s)-0.215842 \cos (\beta_1 s)
	-5.50054 \sinh (\beta_1 s)+0.215842 \cosh (\beta_1 s), \quad \beta_1^2 = 1.38569.}
\end{equation}

We note that \eqref{Eq: freq} reduces to $1 + \cos(\beta) \cosh(\beta) = 0 $ for the case that there is no lumped mass at the tip of beam; this in fact gives the natural frequencies of a linear cantilever beam. In this case, the first eigenvalue is computed as $\beta_1^2 = \omega_1 = 3.51602$, which results to the following first normalized eigenfunction, given in Fig. \ref{Fig: Mode Shape Tip Mass} (right). 

\begin{equation}\label{eq:eig_no_mass}
{ \phi(s) = X_1(s) = 0.734096 \sin (\beta_1 s) - \cos (\beta_1 s)-0.734096 \sinh (\beta_1 s)
+\cosh (\beta_1 s), \quad \beta^2_1 = 3.51602.}
\end{equation}

\section*{Acknowledgments}
This work was supported by the ARO Young Investigator Program Award (W911NF-19-1-0444), and the National Science Foundation Award (DMS-1923201), also partially by MURI/ARO (W911NF-15-1-0562) and the AFOSR Young Investigator Program Award (FA9550-17-1-0150).

           
\bibliographystyle{elsarticle-num-names}
\bibliography{references}

\begin{thebibliography}{73}
\expandafter\ifx\csname natexlab\endcsname\relax\def\natexlab#1{#1}\fi
\providecommand{\url}[1]{\texttt{#1}}
\providecommand{\href}[2]{#2}
\providecommand{\path}[1]{#1}
\providecommand{\DOIprefix}{doi:}
\providecommand{\ArXivprefix}{arXiv:}
\providecommand{\URLprefix}{URL: }
\providecommand{\Pubmedprefix}{pmid:}
\providecommand{\doi}[1]{\href{http://dx.doi.org/#1}{\path{#1}}}
\providecommand{\Pubmed}[1]{\href{pmid:#1}{\path{#1}}}
\providecommand{\bibinfo}[2]{#2}
\ifx\xfnm\relax \def\xfnm[#1]{\unskip,\space#1}\fi
\bibitem[{Sagaut and Cambon(2018)}]{Sagaut2018}
\bibinfo{author}{P.~Sagaut}, \bibinfo{author}{C.~Cambon},
  \bibinfo{title}{Homogeneous Turbulence Dynamics},
  \bibinfo{publisher}{Springer}, \bibinfo{year}{2018}.
\bibitem[{Akhavan-Safaei et~al.(2020)Akhavan-Safaei, Seyedi, and
  Zayernouri}]{AkhavanSafaei2020}
\bibinfo{author}{A.~Akhavan-Safaei}, \bibinfo{author}{S.~Seyedi},
  \bibinfo{author}{M.~Zayernouri},
\newblock \bibinfo{title}{Anomalous features in internal cylinder flow
  instabilities subject to uncertain rotational effects},
\newblock \bibinfo{journal}{Physics of Fluids} \bibinfo{volume}{32}
  (\bibinfo{year}{2020}) \bibinfo{pages}{094107}.
\bibitem[{Habtour et~al.(2016)Habtour, Cole, Riddick, Weiss, Robeson,
  Sridharan, and Dasgupta}]{habtour2016detection}
\bibinfo{author}{E.~Habtour}, \bibinfo{author}{D.~P. Cole},
  \bibinfo{author}{J.~C. Riddick}, \bibinfo{author}{V.~Weiss},
  \bibinfo{author}{M.~Robeson}, \bibinfo{author}{R.~Sridharan},
  \bibinfo{author}{A.~Dasgupta},
\newblock \bibinfo{title}{Detection of fatigue damage precursor using a
  nonlinear vibration approach},
\newblock \bibinfo{journal}{Structural Control and Health Monitoring}
  \bibinfo{volume}{23} (\bibinfo{year}{2016}) \bibinfo{pages}{1442--1463}.
\bibitem[{Kapnistos et~al.(2008)Kapnistos, Lang, Vlassopoulos,
  Pyckhout-Hintzen, Richter, Cho, Chang, and Rubinstein}]{Kapnistos2008}
\bibinfo{author}{M.~Kapnistos}, \bibinfo{author}{M.~Lang},
  \bibinfo{author}{D.~Vlassopoulos}, \bibinfo{author}{D.~Pyckhout-Hintzen},
  \bibinfo{author}{D.~Richter}, \bibinfo{author}{D.~Cho},
  \bibinfo{author}{T.~Chang}, \bibinfo{author}{M.~Rubinstein},
\newblock \bibinfo{title}{Unexpected power-law stress relaxation of entangled
  ring polymers},
\newblock \bibinfo{journal}{Nature Materials} \bibinfo{volume}{7}
  (\bibinfo{year}{2008}) \bibinfo{pages}{997--1002}.
\bibitem[{McKinley and Jaishankar(2013)}]{McKinley2013}
\bibinfo{author}{G.~McKinley}, \bibinfo{author}{A.~Jaishankar},
  \bibinfo{title}{Critical gels, scott blair and the fractional calculus of
  soft squishy materials}, \bibinfo{howpublished}{Presentation},
  \bibinfo{year}{2013}.
\bibitem[{Metzler and Klafter(2000)}]{Metzler2000}
\bibinfo{author}{R.~Metzler}, \bibinfo{author}{J.~Klafter},
\newblock \bibinfo{title}{The random walk's guide to anomalous diffusion: a
  fractional dynamics approach},
\newblock \bibinfo{journal}{Physics Reports} \bibinfo{volume}{339}
  (\bibinfo{year}{2000}) \bibinfo{pages}{1 -- 77}.
\bibitem[{Nnetu et~al.(2013)Nnetu, Knorr, Pawlizak, Fuhs, and
  K\"{a}s}]{Nnetu2013}
\bibinfo{author}{K.~Nnetu}, \bibinfo{author}{M.~Knorr},
  \bibinfo{author}{S.~Pawlizak}, \bibinfo{author}{T.~Fuhs},
  \bibinfo{author}{J.~K\"{a}s},
\newblock \bibinfo{title}{Slow and anomalous dynamics of an mcf-10a epithelial
  cell monolayer},
\newblock \bibinfo{journal}{Soft Matter} \bibinfo{volume}{39}
  (\bibinfo{year}{2013}).
\bibitem[{Wong et~al.(2004)Wong, Gardel, Reichman, Weeks, Valentine, Bausch,
  and Weitz}]{Wong2004}
\bibinfo{author}{I.~Wong}, \bibinfo{author}{M.~Gardel},
  \bibinfo{author}{D.~Reichman}, \bibinfo{author}{E.~Weeks},
  \bibinfo{author}{M.~Valentine}, \bibinfo{author}{A.~Bausch},
  \bibinfo{author}{D.~Weitz},
\newblock \bibinfo{title}{Anomalous diffusion probes microstructure dynamics of
  entangled f-actin networks},
\newblock \bibinfo{journal}{Physical Review Letters} \bibinfo{volume}{92}
  (\bibinfo{year}{2004}).
\bibitem[{Bonadkar et~al.(2016)Bonadkar, Gerum, Kuhn, Sporer, Lippert,
  Schneider, Aifantis, and Fabry}]{Bonadkar2016}
\bibinfo{author}{N.~Bonadkar}, \bibinfo{author}{R.~Gerum},
  \bibinfo{author}{M.~Kuhn}, \bibinfo{author}{M.~Sporer},
  \bibinfo{author}{A.~Lippert}, \bibinfo{author}{W.~Schneider},
  \bibinfo{author}{K.~Aifantis}, \bibinfo{author}{B.~Fabry},
\newblock \bibinfo{title}{Mechanical plasticity of cells},
\newblock \bibinfo{journal}{Nat. Mater.} \bibinfo{volume}{15}
  (\bibinfo{year}{2016}) \bibinfo{pages}{1090 -- 1094}.
\bibitem[{Richeton et~al.(2005)Richeton, Weiss, and Louchet}]{Richeton2005}
\bibinfo{author}{T.~Richeton}, \bibinfo{author}{J.~Weiss},
  \bibinfo{author}{F.~Louchet},
\newblock \bibinfo{title}{Breakdown of avalanche critical behaviour in
  polycrystalline plasticity},
\newblock \bibinfo{journal}{Nat. Mater.} \bibinfo{volume}{4}
  (\bibinfo{year}{2005}) \bibinfo{pages}{465 -- 469}.
\bibitem[{Christensen(2012)}]{christensen2012theory}
\bibinfo{author}{R.~Christensen}, \bibinfo{title}{Theory of viscoelasticity: an
  introduction}, \bibinfo{publisher}{Elsevier}, \bibinfo{year}{2012}.
\bibitem[{Pipkin(2012)}]{pipkin2012lectures}
\bibinfo{author}{A.~Pipkin}, \bibinfo{title}{Lectures on viscoelasticity
  theory}, volume~\bibinfo{volume}{7}, \bibinfo{publisher}{Springer Science \&
  Business Media}, \bibinfo{year}{2012}.
\bibitem[{Aster et~al.(2019)Aster, Borchers, and Thurber}]{aster2018parameter}
\bibinfo{author}{R.~Aster}, \bibinfo{author}{B.~Borchers},
  \bibinfo{author}{C.~Thurber}, \bibinfo{title}{Parameter estimation and
  inverse problems (third edition)}, \bibinfo{publisher}{Elsevier},
  \bibinfo{year}{2019}.
\bibitem[{Bagley(1989)}]{bagley1989power}
\bibinfo{author}{R.~Bagley},
\newblock \bibinfo{title}{Power law and fractional calculus model of
  viscoelasticity},
\newblock \bibinfo{journal}{AIAA journal} \bibinfo{volume}{27}
  (\bibinfo{year}{1989}) \bibinfo{pages}{1412--1417}.
\bibitem[{Jaishankar and McKinley(2013)}]{Jaishankar2013}
\bibinfo{author}{A.~Jaishankar}, \bibinfo{author}{G.~McKinley},
\newblock \bibinfo{title}{Power-law rheology in the bulk and at the interface:
  quasi-properties and fractional constitutive equations},
\newblock \bibinfo{journal}{Proc R Soc A 469: 20120284}
  (\bibinfo{year}{2013}).
\bibitem[{Nutting(1921)}]{nutting1921new}
\bibinfo{author}{P.~Nutting},
\newblock \bibinfo{title}{A new general law of deformation},
\newblock \bibinfo{journal}{Journal of the Franklin Institute}
  \bibinfo{volume}{191} (\bibinfo{year}{1921}) \bibinfo{pages}{679--685}.
\bibitem[{Gemant(1936)}]{gemant1936method}
\bibinfo{author}{A.~Gemant},
\newblock \bibinfo{title}{A method of analyzing experimental results obtained
  from elasto-viscous bodies},
\newblock \bibinfo{journal}{Physics} \bibinfo{volume}{7} (\bibinfo{year}{1936})
  \bibinfo{pages}{311--317}.
\bibitem[{Bagley and Torvik(1983)}]{bagley1983theoretical}
\bibinfo{author}{R.~Bagley}, \bibinfo{author}{P.~Torvik},
\newblock \bibinfo{title}{A theoretical basis for the application of fractional
  calculus to viscoelasticity},
\newblock \bibinfo{journal}{Journal of Rheology} \bibinfo{volume}{27}
  (\bibinfo{year}{1983}) \bibinfo{pages}{201--210}.
\bibitem[{Naghibolhosseini(2015)}]{naghibolhosseini2015estimation}
\bibinfo{author}{M.~Naghibolhosseini}, \bibinfo{title}{Estimation of
  outer-middle ear transmission using \uppercase{DPOAE}s and fractional-order
  modeling of human middle ear}, Ph.D. thesis, City University of New York,
  NY., \bibinfo{year}{2015}.
\bibitem[{Naghibolhosseini and Long(2018)}]{naghibolhosseini2018fractional}
\bibinfo{author}{M.~Naghibolhosseini}, \bibinfo{author}{G.~Long},
\newblock \bibinfo{title}{Fractional-order modelling and simulation of human
  ear},
\newblock \bibinfo{journal}{International Journal of Computer Mathematics}
  \bibinfo{volume}{95} (\bibinfo{year}{2018}) \bibinfo{pages}{1257--1273}.
\bibitem[{Suzuki et~al.(2016)Suzuki, Zayernouri, Bittencourt, and
  Karniadakis}]{suzuki2016fractional}
\bibinfo{author}{J.~Suzuki}, \bibinfo{author}{M.~Zayernouri},
  \bibinfo{author}{M.~Bittencourt}, \bibinfo{author}{G.~Karniadakis},
\newblock \bibinfo{title}{Fractional-order uniaxial visco-elasto-plastic models
  for structural analysis},
\newblock \bibinfo{journal}{Computer Methods in Applied Mechanics and
  Engineering} \bibinfo{volume}{308} (\bibinfo{year}{2016})
  \bibinfo{pages}{443--467}.
\bibitem[{Shitikova et~al.(2017)Shitikova, Rossikhin, and
  Kandu}]{shitikova2017interaction}
\bibinfo{author}{M.~Shitikova}, \bibinfo{author}{Y.~Rossikhin},
  \bibinfo{author}{V.~Kandu},
\newblock \bibinfo{title}{Interaction of internal and external resonances
  during force driven vibrations of a nonlinear thin plate embedded into a
  fractional derivative medium},
\newblock \bibinfo{journal}{Procedia engineering} \bibinfo{volume}{199}
  (\bibinfo{year}{2017}) \bibinfo{pages}{832--837}.
\bibitem[{Rossikhin and Shitikova(1997)}]{rossikhin1997applications}
\bibinfo{author}{Y.~Rossikhin}, \bibinfo{author}{M.~Shitikova},
\newblock \bibinfo{title}{Applications of fractional calculus to dynamic
  problems of linear and nonlinear hereditary mechanics of solids},
\newblock \bibinfo{journal}{Applied Mechanics Reviews} \bibinfo{volume}{50}
  (\bibinfo{year}{1997}) \bibinfo{pages}{15--67}.
\bibitem[{Samiee et~al.(2020)Samiee, Akhavan-Safaei, and
  Zayernouri}]{Samiee2019turbulence}
\bibinfo{author}{M.~Samiee}, \bibinfo{author}{A.~Akhavan-Safaei},
  \bibinfo{author}{M.~Zayernouri},
\newblock \bibinfo{title}{A fractional subgrid-scale model for turbulent flows:
  Theoretical formulation and a priori study},
\newblock \bibinfo{journal}{Physics of Fluids} \bibinfo{volume}{32}
  (\bibinfo{year}{2020}) \bibinfo{pages}{055102}.
\bibitem[{Lorenzo and Hartley(2002)}]{lorenzo2002variable}
\bibinfo{author}{C.~Lorenzo}, \bibinfo{author}{T.~Hartley},
\newblock \bibinfo{title}{Variable order and distributed order fractional
  operators},
\newblock \bibinfo{journal}{Nonlinear dynamics} \bibinfo{volume}{29}
  (\bibinfo{year}{2002}) \bibinfo{pages}{57--98}.
\bibitem[{Atanackovic et~al.(2009)Atanackovic, Oparnica, and
  Pilipovi{\'c}}]{Atanackovic2009distributional}
\bibinfo{author}{T.~Atanackovic}, \bibinfo{author}{L.~Oparnica},
  \bibinfo{author}{S.~Pilipovi{\'c}},
\newblock \bibinfo{title}{Distributional framework for solving fractional
  differential equations},
\newblock \bibinfo{journal}{Integral Transforms and Special Functions}
  \bibinfo{volume}{20} (\bibinfo{year}{2009}) \bibinfo{pages}{215--222}.
\bibitem[{Caputo(1995)}]{caputo1995m}
\bibinfo{author}{M.~Caputo},
\newblock \bibinfo{title}{Mean fractional-order-derivatives differential
  equations and filters},
\newblock \bibinfo{journal}{Annali dell'Universit\`{a} di Ferrara}
  \bibinfo{volume}{41} (\bibinfo{year}{1995}) \bibinfo{pages}{73--84}.
\bibitem[{Caputo(2001)}]{caputo2001m}
\bibinfo{author}{M.~Caputo},
\newblock \bibinfo{title}{Distributed order differential equation modelling
  dielectric induction and diffusion},
\newblock \bibinfo{journal}{Fract. Calc. Appl. Anal.} \bibinfo{volume}{4}
  (\bibinfo{year}{2001}) \bibinfo{pages}{421--442}.
\bibitem[{Chechkin et~al.(2002)Chechkin, Gorenflo, and Sokolov}]{chechkin2002}
\bibinfo{author}{A.~Chechkin}, \bibinfo{author}{R.~Gorenflo},
  \bibinfo{author}{I.~Sokolov},
\newblock \bibinfo{title}{Retarding subdiffusion and accelerating
  superdiffusion governed by distributed-order fractional diffusion equations},
\newblock \bibinfo{journal}{Physical Review E} \bibinfo{volume}{66}
  (\bibinfo{year}{2002}) \bibinfo{pages}{046129}.
\bibitem[{Chechkin et~al.(2008)Chechkin, Gonchar, Gorenflo, Korabel, and
  Sokolov}]{chechkin2008generalized}
\bibinfo{author}{A.~Chechkin}, \bibinfo{author}{V.~Gonchar},
  \bibinfo{author}{R.~Gorenflo}, \bibinfo{author}{N.~Korabel},
  \bibinfo{author}{I.~Sokolov},
\newblock \bibinfo{title}{Generalized fractional diffusion equations for
  accelerating subdiffusion and truncated l{\'e}vy flights},
\newblock \bibinfo{journal}{Physical Review E} \bibinfo{volume}{78}
  (\bibinfo{year}{2008}) \bibinfo{pages}{021111}.
\bibitem[{Li et~al.(2011)Li, Sheng, and Chen}]{li2011distributed}
\bibinfo{author}{Y.~Li}, \bibinfo{author}{H.~Sheng}, \bibinfo{author}{Y.~Chen},
\newblock \bibinfo{title}{On distributed order integrator/differentiator},
\newblock \bibinfo{journal}{Signal Processing} \bibinfo{volume}{91}
  (\bibinfo{year}{2011}) \bibinfo{pages}{1079--1084}.
\bibitem[{Li and Chen(2014)}]{li2014lyapunov}
\bibinfo{author}{Y.~Li}, \bibinfo{author}{Y.~Chen},
\newblock \bibinfo{title}{Lyapunov stability of fractional-order nonlinear
  systems: A distributed-order approach},
\newblock in: \bibinfo{booktitle}{ICFDA'14 International Conference on
  Fractional Differentiation and Its Applications 2014},
  \bibinfo{organization}{IEEE}, \bibinfo{year}{2014}, pp.
  \bibinfo{pages}{1--6}.
\bibitem[{Duan and Baleanu(2018)}]{duan2018steady}
\bibinfo{author}{J.~Duan}, \bibinfo{author}{D.~Baleanu},
\newblock \bibinfo{title}{Steady periodic response for a vibration system with
  distributed order derivatives to periodic excitation},
\newblock \bibinfo{journal}{Journal of Vibration and Control}
  \bibinfo{volume}{24} (\bibinfo{year}{2018}) \bibinfo{pages}{3124--3131}.
\bibitem[{Bagley and Torvik(2000)}]{bagley2000}
\bibinfo{author}{R.~Bagley}, \bibinfo{author}{P.~Torvik},
\newblock \bibinfo{title}{On the existence of the order domain and the solution
  of distributed order equations-part i},
\newblock \bibinfo{journal}{International Journal of Applied Mathematics}
  \bibinfo{volume}{2} (\bibinfo{year}{2000}) \bibinfo{pages}{865--882}.
\bibitem[{Kharazmi et~al.(2017)Kharazmi, Zayernouri, and
  Karniadakis}]{kharazmi2017petrov}
\bibinfo{author}{E.~Kharazmi}, \bibinfo{author}{M.~Zayernouri},
  \bibinfo{author}{G.~Karniadakis},
\newblock \bibinfo{title}{\uppercase{P}etrov--\uppercase{G}alerkin and spectral
  collocation methods for distributed order differential equations},
\newblock \bibinfo{journal}{SIAM Journal on Scientific Computing}
  \bibinfo{volume}{39} (\bibinfo{year}{2017}) \bibinfo{pages}{A1003--A1037}.
\bibitem[{Kharazmi and Zayernouri(2018)}]{kharazmi2018fractional}
\bibinfo{author}{E.~Kharazmi}, \bibinfo{author}{M.~Zayernouri},
\newblock \bibinfo{title}{Fractional pseudo-spectral methods for
  distributed-order fractional pdes},
\newblock \bibinfo{journal}{International Journal of Computer Mathematics}
  (\bibinfo{year}{2018}) \bibinfo{pages}{1--22}.
\bibitem[{{\L}ab{\k{e}}dzki et~al.(2018){\L}ab{\k{e}}dzki, Pawlikowski, and
  Radowicz}]{labkedzki2018transverse}
\bibinfo{author}{P.~{\L}ab{\k{e}}dzki}, \bibinfo{author}{R.~Pawlikowski},
  \bibinfo{author}{A.~Radowicz},
\newblock \bibinfo{title}{Transverse vibration of a cantilever beam under base
  excitation using fractional rheological model},
\newblock in: \bibinfo{booktitle}{AIP Conference Proceedings}, volume
  \bibinfo{volume}{2029}, \bibinfo{organization}{AIP Publishing},
  \bibinfo{year}{2018}, p. \bibinfo{pages}{020034}.
\bibitem[{Ansari et~al.(2016)Ansari, Oskouie, and Gholami}]{Ansari2016Nanobeam}
\bibinfo{author}{R.~Ansari}, \bibinfo{author}{M.~F. Oskouie},
  \bibinfo{author}{R.~Gholami},
\newblock \bibinfo{title}{Size-dependent geometrically nonlinear free vibration
  analysis of fractional viscoelastic nanobeams based on the nonlocal
  elasticity theory},
\newblock \bibinfo{journal}{Physica E: Low-dimensional Systems and
  Nanostructures} \bibinfo{volume}{75} (\bibinfo{year}{2016})
  \bibinfo{pages}{266 -- 271}.
\bibitem[{Faraji~Oskouie et~al.(2017)Faraji~Oskouie, Ansari, and
  Sadeghi}]{FarajiOskouie2017}
\bibinfo{author}{M.~Faraji~Oskouie}, \bibinfo{author}{R.~Ansari},
  \bibinfo{author}{F.~Sadeghi},
\newblock \bibinfo{title}{Nonlinear vibration analysis of fractional
  viscoelastic euler---bernoulli nanobeams based on the surface stress theory},
\newblock \bibinfo{journal}{Acta Mechanica Solida Sinica} \bibinfo{volume}{30}
  (\bibinfo{year}{2017}) \bibinfo{pages}{416--424}.
\bibitem[{Eyebe et~al.(2017)Eyebe, Betchewe, Mohamadou, and Kofane}]{Eyebe2018}
\bibinfo{author}{G.~Eyebe}, \bibinfo{author}{G.~Betchewe},
  \bibinfo{author}{A.~Mohamadou}, \bibinfo{author}{T.~Kofane},
\newblock \bibinfo{title}{Nonlinear vibration of a nonlocal nanobeam resting on
  fractional-order viscoelastic pasternak foundations},
\newblock \bibinfo{journal}{Fractal and Fractional} \bibinfo{volume}{2}
  (\bibinfo{year}{2017}).
\bibitem[{Lewandowski and Wielentejczyk(2017)}]{lewandowski2017nonlinear}
\bibinfo{author}{R.~Lewandowski}, \bibinfo{author}{P.~Wielentejczyk},
\newblock \bibinfo{title}{Nonlinear vibration of viscoelastic beams described
  using fractional order derivatives},
\newblock \bibinfo{journal}{Journal of Sound and Vibration}
  \bibinfo{volume}{399} (\bibinfo{year}{2017}) \bibinfo{pages}{228--243}.
\bibitem[{Samiee et~al.(2019{\natexlab{a}})Samiee, Zayernouri, and
  Meerschaert}]{samiee2019unified1}
\bibinfo{author}{M.~Samiee}, \bibinfo{author}{M.~Zayernouri},
  \bibinfo{author}{M.~Meerschaert},
\newblock \bibinfo{title}{A unified spectral method for fpdes with two-sided
  derivatives; part i: a fast solver},
\newblock \bibinfo{journal}{Journal of Computational Physics}
  \bibinfo{volume}{385} (\bibinfo{year}{2019}{\natexlab{a}})
  \bibinfo{pages}{225--243}.
\bibitem[{Samiee et~al.(2019{\natexlab{b}})Samiee, Zayernouri, and
  Meerschaert}]{samiee2019unified2}
\bibinfo{author}{M.~Samiee}, \bibinfo{author}{M.~Zayernouri},
  \bibinfo{author}{M.~Meerschaert},
\newblock \bibinfo{title}{A unified spectral method for fpdes with two-sided
  derivatives; part ii: Stability, and error analysis},
\newblock \bibinfo{journal}{Journal of Computational Physics}
  \bibinfo{volume}{385} (\bibinfo{year}{2019}{\natexlab{b}})
  \bibinfo{pages}{244--261}.
\bibitem[{Samiee et~al.(2018)Samiee, Kharazmi, Zayernouri, and
  Meerschaert}]{samiee2018petrov}
\bibinfo{author}{M.~Samiee}, \bibinfo{author}{E.~Kharazmi},
  \bibinfo{author}{M.~Zayernouri}, \bibinfo{author}{M.~Meerschaert},
\newblock \bibinfo{title}{Petrov-galerkin method for fully distributed-order
  fractional partial differential equations},
\newblock \bibinfo{journal}{arXiv preprint arXiv:1805.08242}
  (\bibinfo{year}{2018}).
\bibitem[{Lubich(1986)}]{lubich1986discretized}
\bibinfo{author}{C.~Lubich},
\newblock \bibinfo{title}{Discretized fractional calculus},
\newblock \bibinfo{journal}{SIAM Journal on Mathematical Analysis}
  \bibinfo{volume}{17} (\bibinfo{year}{1986}) \bibinfo{pages}{704--719}.
\bibitem[{Zayernouri et~al.(2015)Zayernouri, Ainsworth, and
  Karniadakis}]{zayernouri2015tempered}
\bibinfo{author}{M.~Zayernouri}, \bibinfo{author}{M.~Ainsworth},
  \bibinfo{author}{G.~Karniadakis},
\newblock \bibinfo{title}{Tempered fractional sturm--liouville eigenproblems},
\newblock \bibinfo{journal}{SIAM Journal on Scientific Computing}
  \bibinfo{volume}{37} (\bibinfo{year}{2015}) \bibinfo{pages}{A1777--A1800}.
\bibitem[{Suzuki and Zayernouri(2020)}]{suzuki2018automated}
\bibinfo{author}{J.~Suzuki}, \bibinfo{author}{M.~Zayernouri},
\newblock \bibinfo{title}{A self-singularity-capturing scheme for fractional
  differential equations},
\newblock \bibinfo{journal}{International Journal of Computer Mathematics}
  (\bibinfo{year}{2020}) \bibinfo{pages}{1--28}.
\bibitem[{Zayernouri and Matzavinos(2016)}]{zayernouri2016fractionalAdams}
\bibinfo{author}{M.~Zayernouri}, \bibinfo{author}{A.~Matzavinos},
\newblock \bibinfo{title}{Fractional adams--bashforth/moulton methods: an
  application to the fractional keller--segel chemotaxis system},
\newblock \bibinfo{journal}{Journal of Computational Physics}
  \bibinfo{volume}{317} (\bibinfo{year}{2016}) \bibinfo{pages}{1--14}.
\bibitem[{Lin and Xu(2007)}]{lin2007finite}
\bibinfo{author}{Y.~Lin}, \bibinfo{author}{C.~Xu},
\newblock \bibinfo{title}{Finite difference/spectral approximations for the
  time-fractional diffusion equation},
\newblock \bibinfo{journal}{Journal of Computational Physics}
  \bibinfo{volume}{225} (\bibinfo{year}{2007}) \bibinfo{pages}{1533--1552}.
\bibitem[{Zhou et~al.(2020)Zhou, Suzuki, Zhang, and Zayernouri}]{Zhou2019IMEX}
\bibinfo{author}{Y.~Zhou}, \bibinfo{author}{J.~Suzuki},
  \bibinfo{author}{C.~Zhang}, \bibinfo{author}{M.~Zayernouri},
\newblock \bibinfo{title}{Implicit-explicit time integration of nonlinear
  fractional differential equations},
\newblock \bibinfo{journal}{Applied Numerical Mathematics}
  \bibinfo{volume}{156} (\bibinfo{year}{2020}) \bibinfo{pages}{555--583}.
\bibitem[{Mashayekhi et~al.(2019)Mashayekhi, Hussaini, and
  Oates}]{Mashayekhi2019Fractal}
\bibinfo{author}{S.~Mashayekhi}, \bibinfo{author}{Y.~Hussaini},
  \bibinfo{author}{W.~Oates},
\newblock \bibinfo{title}{A physical interpretation of fractional
  viscoelasticity based on the fractal structure of media: Theory and
  experimental validation},
\newblock \bibinfo{journal}{J. Mech. Phys. Solids} \bibinfo{volume}{128}
  (\bibinfo{year}{2019}) \bibinfo{pages}{137--150}.
\bibitem[{Mainardi(2010)}]{mainardi2010fractional}
\bibinfo{author}{F.~Mainardi}, \bibinfo{title}{Fractional calculus and waves in
  linear viscoelasticity: an introduction to mathematical models},
  \bibinfo{publisher}{World Scientific}, \bibinfo{year}{2010}.
\bibitem[{Rogosin and Mainardi(2014)}]{rogosin2014george}
\bibinfo{author}{S.~Rogosin}, \bibinfo{author}{F.~Mainardi},
\newblock \bibinfo{title}{George william scott blair--the pioneer of factional
  calculus in rheology},
\newblock \bibinfo{journal}{arXiv preprint arXiv:1404.3295}
  (\bibinfo{year}{2014}).
\bibitem[{Mainardi and Gorenflo(2008)}]{mainardi2008time}
\bibinfo{author}{F.~Mainardi}, \bibinfo{author}{R.~Gorenflo},
\newblock \bibinfo{title}{Time-fractional derivatives in relaxation processes:
  a tutorial survey},
\newblock \bibinfo{journal}{arXiv preprint arXiv:0801.4914}
  (\bibinfo{year}{2008}).
\bibitem[{Mainardi and Spada(2011)}]{mainardi2011creep}
\bibinfo{author}{F.~Mainardi}, \bibinfo{author}{G.~Spada},
\newblock \bibinfo{title}{Creep, relaxation and viscosity properties for basic
  fractional models in rheology},
\newblock \bibinfo{journal}{The European Physical Journal Special Topics}
  \bibinfo{volume}{193} (\bibinfo{year}{2011}) \bibinfo{pages}{133--160}.
\bibitem[{Meirovitch(2010)}]{meirovitch2010fundamentals}
\bibinfo{author}{L.~Meirovitch}, \bibinfo{title}{Fundamentals of vibrations},
  \bibinfo{publisher}{Waveland Press}, \bibinfo{year}{2010}.
\bibitem[{Bonet and Wood(1997)}]{bonet1997nonlinear}
\bibinfo{author}{J.~Bonet}, \bibinfo{author}{R.~Wood},
  \bibinfo{title}{Nonlinear continuum mechanics for finite element analysis},
  \bibinfo{publisher}{Cambridge university press}, \bibinfo{year}{1997}.
\bibitem[{Lion(1997)}]{lion1997thermodynamics}
\bibinfo{author}{A.~Lion},
\newblock \bibinfo{title}{On the thermodynamics of fractional damping
  elements},
\newblock \bibinfo{journal}{Continuum Mechanics and Thermodynamics}
  \bibinfo{volume}{9} (\bibinfo{year}{1997}) \bibinfo{pages}{83--96}.
\bibitem[{Tadmor(2012)}]{tadmor2012review}
\bibinfo{author}{E.~Tadmor},
\newblock \bibinfo{title}{A review of numerical methods for nonlinear partial
  differential equations},
\newblock \bibinfo{journal}{Bulletin of the American Mathematical Society}
  \bibinfo{volume}{49} (\bibinfo{year}{2012}) \bibinfo{pages}{507--554}.
\bibitem[{Azrar et~al.(1999)Azrar, Benamar, and White}]{azrar1999semi}
\bibinfo{author}{L.~Azrar}, \bibinfo{author}{R.~Benamar},
  \bibinfo{author}{R.~White},
\newblock \bibinfo{title}{Semi-analytical approach to the non-linear dynamic
  response problem of s--s and c--c beams at large vibration amplitudes part i:
  general theory and application to the single mode approach to free and forced
  vibration analysis},
\newblock \bibinfo{journal}{Journal of sound and vibration}
  \bibinfo{volume}{224} (\bibinfo{year}{1999}) \bibinfo{pages}{183--207}.
\bibitem[{Tseng and Dugundji(1971)}]{tseng1971nonlinear}
\bibinfo{author}{W.-Y. Tseng}, \bibinfo{author}{J.~Dugundji},
\newblock \bibinfo{title}{Nonlinear vibrations of a buckled beam under harmonic
  excitation}  (\bibinfo{year}{1971}).
\bibitem[{Loutridis et~al.(2005)Loutridis, Douka, and
  Hadjileontiadis}]{loutridis2005forced}
\bibinfo{author}{S.~Loutridis}, \bibinfo{author}{E.~Douka},
  \bibinfo{author}{L.~Hadjileontiadis},
\newblock \bibinfo{title}{Forced vibration behaviour and crack detection of
  cracked beams using instantaneous frequency},
\newblock \bibinfo{journal}{Ndt \& E International} \bibinfo{volume}{38}
  (\bibinfo{year}{2005}) \bibinfo{pages}{411--419}.
\bibitem[{Hamdan and Dado(1997)}]{hamdan1997large}
\bibinfo{author}{M.~Hamdan}, \bibinfo{author}{M.~Dado},
\newblock \bibinfo{title}{Large amplitude free vibrations of a uniform
  cantilever beam carrying an intermediate lumped mass and rotary inertia},
\newblock \bibinfo{journal}{Journal of Sound and Vibration}
  \bibinfo{volume}{206} (\bibinfo{year}{1997}) \bibinfo{pages}{151--168}.
\bibitem[{Lestari and Hanagud(2001)}]{lestari2001nonlinear}
\bibinfo{author}{W.~Lestari}, \bibinfo{author}{S.~Hanagud},
\newblock \bibinfo{title}{Nonlinear vibration of buckled beams: some exact
  solutions},
\newblock \bibinfo{journal}{International Journal of Solids and Structures}
  \bibinfo{volume}{38} (\bibinfo{year}{2001}) \bibinfo{pages}{4741--4757}.
\bibitem[{Eisley(1964)}]{eisley1964nonlinear}
\bibinfo{author}{J.~G. Eisley},
\newblock \bibinfo{title}{Nonlinear vibration of beams and rectangular plates},
\newblock \bibinfo{journal}{Zeitschrift f{\"u}r angewandte Mathematik und
  Physik ZAMP} \bibinfo{volume}{15} (\bibinfo{year}{1964})
  \bibinfo{pages}{167--175}.
\bibitem[{Hsu(1960)}]{hsu1960application}
\bibinfo{author}{C.~Hsu},
\newblock \bibinfo{title}{On the application of elliptic functions in
  non-linear forced oscillations},
\newblock \bibinfo{journal}{Quarterly of Applied Mathematics}
  \bibinfo{volume}{17} (\bibinfo{year}{1960}) \bibinfo{pages}{393--407}.
\bibitem[{Pillai and Rao(1992)}]{pillai1992nonlinear}
\bibinfo{author}{S.~Pillai}, \bibinfo{author}{B.~N. Rao},
\newblock \bibinfo{title}{On nonlinear free vibrations of simply supported
  uniform beams},
\newblock \bibinfo{journal}{Journal of sound and vibration}
  \bibinfo{volume}{159} (\bibinfo{year}{1992}) \bibinfo{pages}{527--531}.
\bibitem[{Evensen(1968)}]{evensen1968nonlinear}
\bibinfo{author}{D.~A. Evensen},
\newblock \bibinfo{title}{Nonlinear vibrations of beams with various boundary
  conditions.},
\newblock \bibinfo{journal}{AIAA journal} \bibinfo{volume}{6}
  (\bibinfo{year}{1968}) \bibinfo{pages}{370--372}.
\bibitem[{Svenkeson et~al.(2016)Svenkeson, Glaz, Stanton, and
  West}]{svenkeson2016spectral}
\bibinfo{author}{A.~Svenkeson}, \bibinfo{author}{B.~Glaz},
  \bibinfo{author}{S.~Stanton}, \bibinfo{author}{B.~West},
\newblock \bibinfo{title}{Spectral decomposition of nonlinear systems with
  memory},
\newblock \bibinfo{journal}{Physical Review E} \bibinfo{volume}{93}
  (\bibinfo{year}{2016}) \bibinfo{pages}{022211}.
\bibitem[{Shoshani et~al.(2017)Shoshani, Shaw, and
  Dykman}]{shoshani2017anomalous}
\bibinfo{author}{O.~Shoshani}, \bibinfo{author}{S.~Shaw},
  \bibinfo{author}{M.~Dykman},
\newblock \bibinfo{title}{Anomalous decay of nanomechanical modes going through
  nonlinear resonance},
\newblock \bibinfo{journal}{Scientific reports} \bibinfo{volume}{7}
  (\bibinfo{year}{2017}) \bibinfo{pages}{18091}.
\bibitem[{Nayfeh and Mook(2008)}]{nayfeh2008nonlinear}
\bibinfo{author}{A.~Nayfeh}, \bibinfo{author}{D.~Mook},
  \bibinfo{title}{Nonlinear oscillations}, \bibinfo{publisher}{John Wiley \&
  Sons}, \bibinfo{year}{2008}.
\bibitem[{Rossikhin and Shitikova(2010)}]{rossikhin2010application}
\bibinfo{author}{Y.~A. Rossikhin}, \bibinfo{author}{M.~Shitikova},
\newblock \bibinfo{title}{Application of fractional calculus for dynamic
  problems of solid mechanics: novel trends and recent results},
\newblock \bibinfo{journal}{Applied Mechanics Reviews} \bibinfo{volume}{63}
  (\bibinfo{year}{2010}) \bibinfo{pages}{010801}.
\bibitem[{Samko et~al.(1993)Samko, Kilbas, and Marichev}]{samko1993fractional}
\bibinfo{author}{S.~Samko}, \bibinfo{author}{A.~Kilbas},
  \bibinfo{author}{O.~Marichev}, \bibinfo{title}{Fractional integrals and
  derivatives: theory and applications}, \bibinfo{publisher}{CRC},
  \bibinfo{year}{1993}.

\end{thebibliography}

\end{document}